\documentclass[12pt]{elsarticle}
\usepackage{amsmath,amssymb}
\usepackage{eepic}
\usepackage{hyperref}
\usepackage{graphicx}
\usepackage{enumitem}\biboptions{sort&compress}
\usepackage[left=2.0cm,right=2.0cm,top=1cm,bottom=2cm,bindingoffset=0cm]{geometry}
\usepackage{color}
\renewcommand{\alarm}[1]{\textcolor{cyan}{#1}} 

\newtheorem{lemma}{Lemma}
\newtheorem{theorem}{Theorem}
\newtheorem{proposition}{Proposition}
\newtheorem{corollary}{Corollary}
\newtheorem{definition}{Definition}
\newtheorem{property}{Property}
\newtheorem{remark}{Remark}

\usepackage{listings, xcolor}
\newcommand{\beginMatlab}[1]
{
	\definecolor{dkgreen}{rgb}{0,0.6,0}
	\definecolor{gray}{rgb}{0.5,0.5,0.5}
	\lstloadlanguages{Matlab}%
	\lstset{
	language=Matlab,
	keywordstyle=[1]\color{blue}\bfseries, 
	keywordstyle=[2]\color{dkgreen},   
	keywordstyle=[3]\color{black},      
	escapebegin=\color{dkgreen},
	keywords={abs,collect,simple,exp,log,eigen,break,case,catch,continue,else,elseif,end,for,
	function,global,if,otherwise,persistent,return,switch,try,while,diff,simplify,
	subs,int2str,syms,clear,all, factor,hold, on, off, plot,solve, title,xlabel, ylabel,
	grid,eval},
	caption = {#1},
	basicstyle=\ttfamily,
	morekeywords={},
	commentstyle=\color{dkgreen},
	stringstyle=\color{gray},
	numbers=left, 
	firstnumber=1, 
	numberstyle=\tiny\color{blue}, 
	numbersep=10pt,
	backgroundcolor=\color{white},
	frame=shadowbox,
	tabsize=4,
	showspaces=false,
	showstringspaces=false,
	numbers=left,
	mathescape
 	}
}

\newenvironment{example}[1][Example.]{\begin{trivlist}
\item[\hskip \labelsep {\bfseries #1}]}{\end{trivlist}}

\journal{arXiv, {\rm DRAFT}}

\DeclareMathOperator{\LCE}{\nu^{L}}
\DeclareMathOperator{\LEs}{\nu}
\DeclareMathOperator{\LEo}{\widetilde{\nu}}

\begin{document}

\begin{frontmatter}

\title{
  A short survey on Lyapunov dimension for \\
  finite dimensional dynamical systems in Euclidean space.
}

 \author{N.V. Kuznetsov\corref{cor}\,}
 \ead{nkuznetsov239@gmail.com}
 \author{G.A.~Leonov}

 \address[spbu]{Faculty of Mathematics and Mechanics,
 St.~Petersburg State University, Russia}
 \address[ras]{
 Institute for Problems in Mechanical Engineering of the Russian Academy of Sciences,
 Russia}
 \address[fin]{Department of Mathematical Information Technology,
 University of Jyv\"{a}skyl\"{a}, Finland}

\begin{abstract}
  Nowadays there are a number of surveys and theoretical works devoted
  to the Lyapunov exponents and Lyapunov dimension,
  however most of them are devoted to infinite dimensional systems
  or rely on special ergodic properties of the system.
  At the same time
  the provided illustrative examples are often finite dimensional systems
  and the rigorous proof of their ergodic properties
  can be a difficult task.
  Also the Lyapunov exponents and Lyapunov dimension
  have become so widespread and common
  that they are often used without references to the rigorous definitions or pioneering works.

  The survey is devoted to the finite dimensional dynamical systems in Euclidean space
  and its aim is to explain, in a simple but rigorous way,
  the connection between the key works in the area:
  by Kaplan and Yorke (the concept of Lyapunov dimension, 1979),
  Douady and Oesterl\'{e}
  (estimation of Hausdorff dimension via the Lyapunov dimension of maps, 1980),
  Constantin, Eden, Foias, and Temam
  (estimation of Hausdorff dimension via the Lyapunov exponents and dimension of dynamical systems, 1985-90),
  Leonov (estimation of the Lyapunov dimension via the direct Lyapunov method, 1991),
  and numerical methods for the computation of Lyapunov exponents and Lyapunov dimension.

  In this survey a concise overview of the classical results is presented,
  various definitions of Lyapunov exponents and Lyapunov dimension are discussed.
  An effective analytical method for the estimation of Lyapunov dimension is presented,
  its application to the self-excited and hidden attractors of well-known dynamical systems
  is demonstrated, and analytical formulas of exact Lyapunov dimension are obtained.
 \end{abstract}


\end{frontmatter}
\setcounter{secnumdepth}{2}
\setcounter{tocdepth}{2}
{\small\tableofcontents}

\section{Introduction: Hausdorff dimension}
The theory of topological dimension \cite{HurewiczW-1941,Kuratowski-1966},
developed in the first half of the 20th century, is of little use in
giving the scale of dimensional characteristics of attractors. The
point is that the topological dimension can take integer values
only. Hence the scale of dimensional characteristics compiled in
this manner turns out to be quite poor.
For the analysis of attractors, the Hausdorff dimension of a set is
much better. This dimensional characteristic can take any
nonnegative value (not greater than the topological dimension of the space),
and it coincides with the topological dimension for such typical objects in Euclidean space
as a smooth curve, a smooth surface, or a countable set of points.

Let us give the definition of the Hausdorff dimension
and its upper estimations based on the Lyapunov exponents
following mainly
(\cite{KaplanY-1979,DouadyO-1980,Ledrappier-1981,ConstantinFT-1985,Smith-1986,EdenFT-1991,Hunt-1996,Temam-1997,IlyashenkoW-1999-AMS,BoichenkoLR-2005,BarreiraG-2011,Leonov-2012-PMM,Chueshov-2015}).

Consider a set $K \subseteq \mathbb{R}^n$ and numbers $d \ge 0$, $\varepsilon>0$.
We cover $K$ by a countable set of balls $B_{r_j}$ of radius $r_j<\varepsilon$, and define
\[
  \mu_H(K,d,\varepsilon) :=\inf\bigg\{\sum\limits_{j\geq 1} r_j^d|\ r_j \leq \varepsilon, K \subset \bigcup_{j\geq 1}B_{r_j} \bigg\},
\]
where the infimum is taken over all such countable $\varepsilon$-coverings $K$.
It is obvious that $\mu_H(K,d,\varepsilon)$ does not decrease with decreasing $\varepsilon$.
Therefore there exists a limit (perhaps infinite), namely
$$
  \mu_H(K,d)=\lim\limits_{\varepsilon\to0+0}\mu_H(K,d,\varepsilon)=\sup_{\varepsilon>0}\mu_H(K,d,\varepsilon).
$$

\begin{definition}
 The function $\mu_H(\cdot,d)$ is called the \emph{Hausdorff $d$-measure} on $\mathbb{R}^n$.
\end{definition}

For a certain set $K$, the function $\mu_H(K,\cdot)$ has the following property.
It is possible to find $d_{cr}=d_{cr}(K)\in[0,n]$ such that
\[
  \begin{array}{l}
  \mu_H(K,d)=\infty,\ \forall d\, \in (0,d_{cr}); \quad 
  \mu_H(K,d)=0,\ \forall\,d>d_{cr}.
  \end{array}
\]
We have $d_{cr}(\mathbb{R}^n)=n$.
\begin{definition}
The \emph{Hausdorff dimension} of the set $K$ is defined as
$$
  \dim_{\rm H} K :=d_{cr}(K)=\inf \{d \geq 0 \mid \mu_H(K,d)=0 \}.
$$
\end{definition}

%
%

\section{Singular value function and invariant sets of maps and dynamical systems}
In the seminal paper \cite{DouadyO-1980}
Douady and Oesterl\'{e}
showed how to obtain an upper estimate of the Hausdorff dimension of set $K$.
To demonstrate their approach,
let us consider some definitions and auxiliary results.

Let $U$ be an open subset of $\mathbb{R}^n$ and $\varphi: U \to \mathbb{R}^n$
be a continuously differentiable map.
With respect to the canonical basis in $\mathbb{R}^n$ the function $\varphi(u)$
has the $n\times n$ Jacobian matrix
\[
  D\varphi(u) =   D_u\varphi(u)
  =\left(
   \frac{\partial\varphi_i(u)}{\partial u_j}
  \right)_{n\times n},
  \quad u \in U.
\]
Let $\sigma_i(u)=\sigma_i(D\varphi(u))$, $i=1,2,...,n$, be the singular values of $D\varphi(u)$
(i.e. $\sigma_i(u)\geq0$ and $\sigma_i(u)^2$ are the eigenvalues of the symmetric matrix $D\varphi(u)^*D\varphi(u)$ with respect to their algebraic multiplicity)
ordered so that
$\sigma_1(u)\geq \cdots \geq \sigma_n(u) \geq 0$ for any $u \in K$.
If $\sigma_n(u) > 0$, then the unit ball $B$ is transformed by $D\varphi(u)$ into the ellipsoid $D\varphi(u)B$
and the lengths of its principal semiaxes coincide with the singular values.
\begin{definition}
The singular value function of $D\varphi(u)$ of order $d \in [0,n]$
at $u \in U$ is defined as
\begin{equation}\label{defomega} 
  \omega_d(D\varphi(u)) := \left\{
  \begin{aligned}
    & 1, && d=0, \\
    & \sigma_1(u)\cdots\sigma_{d}(u), && d \in \{1,\ldots,n\}, \\ 
    & \sigma_1(u)\cdots\sigma_{\lfloor d \rfloor}(u)\sigma_{\lfloor d \rfloor+1}(u)^{d-\lfloor d \rfloor}, && d \in (0,1)\cup\ldots\cup(n-2,n-1),
  \end{aligned}
  \right.
\end{equation}
where ${\lfloor d \rfloor}$ is the largest integer less or equal to $d$.
\end{definition}

Remark that $|\det D\varphi(u)|=\omega_n(D\varphi(u))$.

Similarly, introducing the singular value function for arbitrary quadratic matrices,
by the Horn inequality \cite{HornJ-1994-book}
for any two $n\times n$ matrices $A$ and $B$ and any $d \in [0,n]$
we have (see, e.g. \cite[p.28]{BoichenkoLR-2005})
\begin{lemma}\label{thm:horn}
\begin{equation}\label{horn}
  \omega_d(AC) \leq \omega_d(A)\omega_d(C).
\end{equation}
\end{lemma}

\begin{definition}\label{def:inv_set-map}
  A set $K \subset U \subseteq \mathbb{R}^n$
  with respect to the map $\varphi$
  is said to be:

  1) \emph{positively invariant}  if
  \(
    \varphi(K) \subset K,
  \)

  2) \emph{invariant} if
  \(
    \varphi (K) = K,
  \)

  3) and \emph{negatively invariant} if
 \(
    \varphi(K) \supset K,
 \)
 \noindent
 where $\varphi(K) = \left\{ \varphi(u) ~ | ~ {u} \in K\right \}$.
\end{definition}

Consider an autonomous differential equation
\begin{equation} \label{eq:ode}
 \dot u = f(u),
\end{equation}
where $f : U \subseteq \mathbb{R}^n \to \mathbb{R}^n$ is a continuously differentiable vector-function.
Suppose that any solution $u(t,u_0)$ of \eqref{eq:ode} such that $u(0,u_0)=u_0 \in U$
exists for $t \in [0,\infty)$, is unique, and stays in $U$.
Then the evolutionary operator $\varphi^t(u_0) := u(t,u_0)$
is continuously differentiable and satisfies the semigroup property:
\begin{equation}\label{group_prop}
    \varphi^{t+s}(u_0) = \varphi^{t}(\varphi^{s}(u_0)), \ \varphi^0(u_0)=u_0
    \quad \forall ~ t, s \geq 0, \ \forall u_0 \in U.
\end{equation}
Thus $\{\varphi^t\}_{t\geq0}$ is a smooth \emph{dynamical system}
in the phase space $(U, ||~\cdot~||)$:
$\big(\{\varphi^t\}_{t\geq0},(U \subseteq \mathbb{R}^n,||\cdot||) \big)$.
Here $||u|| = \sqrt{u_1^2 + \cdots + u_n^2}$
is Euclidean norm of the vector
${u} = (u_1, \ldots, u_n) \in \mathbb{R}^n$.
Similarly, we can consider a dynamical system generated by the difference equation
\begin{equation}\label{eq:odife}
   u(t+1) = \varphi(u(t)), \quad t = 0,1,..\,,
\end{equation}
where $\varphi : U \subseteq \mathbb{R}^n \to U$ is a continuously differentiable vector-function.
Here
$\varphi^t(u) = \underbrace{(\varphi\circ \varphi \circ \cdots \varphi)(u)}\limits_{t-{\rm times}}$, $\varphi^0(u)=u$,
and the existence and uniqueness (in the forward-time direction)
take place for all $t\geq0$.
Further $\{\varphi^t\}_{t\geq0}$ denotes a smooth \emph{dynamical system} with continuous or discrete time.

\begin{definition}\label{def:inv_set}
  A set $K \subset U \subseteq \mathbb{R}^n$
  with respect to the dynamical system $\{\varphi^t\}_{t\geq0}$
  is said to be positively invariant, invariant or negatively invariant
  if the corresponding property takes place
  with respect to the map $\varphi^{t}$ for all $t>0$.

%
%
\end{definition}

Consider the linearizations of systems (\ref{eq:ode}) and \eqref{eq:odife}
along the solution $\varphi^t(u)$:
\begin{equation} \label{sflct}
  \dot y = J(\varphi^t(u))y,  \quad J(u) = Df(u),\\
\end{equation}
\begin{equation} \label{sfldt}
  y(t+1) = J(\varphi^t(u))y(t), \quad J(u) = D\varphi(u),
\end{equation}
where $J(u)$ is the $n\times n$ Jacobian matrix, all elements of which are continuous functions of $u$.
Consider the fundamental matrix
\begin{equation}\label{nfm}
  D\varphi^t(u)=\big(y^1(t),...,y^n(t)\big),
  \quad D\varphi^0(u) = I,
\end{equation}
which consists of linearly independent solutions $\{y^i(t)\}_{i=1}^{n}$ of
the linearized system.
An important cocycle property of fundamental matrix \eqref{nfm} is as follows
\begin{equation}\label{cocycle}
  D\varphi^{t+s}(u) = D\varphi^t\big(\varphi^s(u)\big)D\varphi^{s}(u),
  \ \forall t,s \geq 0, \ \forall u \in U \subseteq \mathbb{R}^n.
\end{equation}
Consider the singular values of the matrix $D\varphi^t(u)$ sorted by descending
for each $t \in [0,+\infty)$ and $u \in U \subseteq \mathbb{R}^n$:
\begin{equation}
  \sigma_i(t,u) := \sigma_i(D\varphi^t(u)), \
  \sigma_1(t,u) \geq ... \geq \sigma_n(t,u) \geq 0
  \quad \forall t\geq 0, u \in U \subseteq \mathbb{R}^n.
\end{equation}
Similar to \eqref{defomega}, we introduce
the singular value function of $D\varphi^t(u)$ of order $d$: $\omega_d(D\varphi^t(u))$.

For a fixed $t \geq 0$ one can consider the map defined by
the evolutionary operator $\varphi^t(u)$: $\varphi^t : U \subseteq \mathbb{R}^n \to U$.

Further we need the following auxiliary statements.

\begin{lemma}\label{thm:wleftcont}
  From formula \eqref{defomega}
  it follows that for any $u \in U$ and $t\geq0$
  the function $d \mapsto \omega_{d}(D\varphi^t(u))$ is a left-continuous function.
\end{lemma}

\begin{lemma}\label{thm:wcontu}
  For any $d \in [0,n]$ and $t\geq0$
  the function $u \mapsto \omega_{d}(D\varphi^t(u))$
  is continuous on $U$ (see, e.g. \cite[p.554]{Gelfert-2003}).
  Therefore for a compact set $K \subset U$ and $t\geq0$ we have
  \begin{equation}\label{wsupmax}
    \sup_{u\in K}\omega_d(D\varphi^{t}(u))=\max_{u\in K}\omega_d(D\varphi^{t}(u)).
\end{equation}
\end{lemma}
{\bf Proof.}
It follows from the continuity of the functions $u \mapsto \sigma_i(D\varphi(u))\ i=1,2,...,n$ on $U$. 
$\blacksquare$

Next, unless otherwise stated,
the invariance of the set $K \subset U \subseteq \mathbb{R}^n$
is considered with respect to the dynamical system
$\big(\{\varphi^t\}_{t\geq0},(U\subseteq \mathbb{R}^n,||\cdot||) \big)$: $\varphi^t(K)\!=\!K,\ \forall t\!\geq\!0.$

\begin{lemma}\label{thm:wsubexp}
  For a compact invariant set $K$
  and any $d \in [0,n]$,
  the function $t \mapsto \max\limits_{u \in K}\, \omega_d(D\varphi^t(u))$
  is sub-exponential, i.e.
\begin{equation}\label{wsubexp}
  \max_{u\in K}\omega_d(D\varphi^{t+s}(u))
  \leq
  \max_{u\in K}\omega_d(D\varphi^{t}(u))
  \max_{u\in K}\omega_d(D\varphi^{s}(u)) \quad \forall t,s\geq0;
\end{equation}
If $\max_{u\in K}\omega_d(D\varphi^{t}(u))>0$ for $t\geq0$,
then $\ln\max_{u\in K}\omega_d(D\varphi^{t+s}(u))$ is subadditive, i.e.
\[
 \ln\max_{u\in K}\omega_d(D\varphi^{t+s}(u))
 \leq
   \ln\max_{u\in K}\omega_d(D\varphi^{t}(u))
   +
   \ln\max_{u\in K}\omega_d(D\varphi^{s}(u)).
\]
\end{lemma}
{\bf Proof.}
By \eqref{cocycle} and \eqref{horn} we get
\[
\begin{aligned}
  & \max_{u\in K}\omega_d(D\varphi^{t+s}(u))
  =
  \max_{u\in K}\big(\omega_d (D\varphi^{t}(\varphi^{s}(u))D\varphi^{s}(u))\big)
  \leq
  \\
  &\leq
  \max_{u\in K}\omega_d(D\varphi^{t}(\varphi^{s}(u)))
  \max_{u\in K}\omega_d(D\varphi^{s}(u))
  \leq
  \max_{u\in K}\omega_d(D\varphi^{t}(u))
  \max_{u\in K}\omega_d(D\varphi^{s}(u)).
\end{aligned}
\]
$\blacksquare$
\begin{corollary}\label{thm:weqpoint}
For an equilibrium point $u_{eq}\equiv\varphi^t(u_{eq})$ we have
\begin{equation}\label{weqpoint}
  \omega_{d}(D\varphi^{t}(u_{eq}))
  = \big(\omega_{d}(D\varphi(u_{eq}))\big)^t, \quad t\geq0.
\end{equation}
\end{corollary}
\begin{corollary}\label{thm:winfliminf}
Remark that for a compact invariant set $K$
\[
  \inf_{t > 0}\max\limits_{u \in K}\omega_{d}(D\varphi^t(u))<1
  \Leftrightarrow
  \liminf_{t \to +\infty}\max\limits_{u \in K}\omega_{d}(D\varphi^t(u))<1.
\]
In this case\footnote{
Considering additional properties
of the dynamical system and the singular value function,
one could get $\lim_{t \to +\infty}$ instead of $\liminf_{t \to +\infty}$,
but we do not need it for our further consideration.
}
\begin{equation}\label{winfliminf}
  \inf_{t > 0}\max\limits_{u \in K}\omega_{d}(D\varphi^t(u))
  =
  \liminf_{t \to +\infty}\max\limits_{u \in K}\omega_{d}(D\varphi^t(u))=0.
\end{equation}
\end{corollary}
{\bf Proof.}
Let
\(
  \inf_{t>0}\max\limits_{u \in K}\omega_{d}(D\varphi^t(u))=M<1.
\)
There are $\delta>0$ and $t_0=t_0(\delta)$ such that
$\max\limits_{u \in K}\omega_{d}(D\varphi^{t_0}(u))\leq 1-\delta$.
Thus by \eqref{wsubexp} we have for $u \in K$ and $n \geq 0$
\[
 0 \leq
 \omega_{d}(D\varphi^{nt_0}(u))
 \leq \max\limits_{u \in K}\omega_{d}(D\varphi^{nt_0}(u))\leq
  (\max\limits_{u \in K}\omega_{d}(D\varphi^{t_0}(u)))^n \leq (1-\delta)^n \to_{n\to+\infty} 0
\]
and therefore
$
M= \liminf_{n \to +\infty}\omega_{d}(D\varphi^{nt_0}(u))=
\liminf_{t \to +\infty}\max\limits_{u \in K}\omega_{d}(D\varphi^t(u))=0$.
The same is true if
we consider
$\liminf_{t \to +\infty}\max\limits_{u \in K}\omega_{d}(D\varphi^t(u))<1$  first.
$\blacksquare$
\begin{corollary}\label{thm:liminf}
If for a fixed $t>0$ and $d \in [0,n]$
we have $\max_{u\in K}\omega_d(D\varphi^{t}(u))<1$, then
\[
  \liminf_{t\to+\infty}\max_{u\in K}\omega_d(D\varphi^{t}(u)) = \liminf_{t\to+\infty}\omega_d(D\varphi^{t}(u))= 0, \quad u \in K.
\]
\end{corollary}

\begin{lemma} \cite[p.33]{ConstantinFT-1985},\cite[pp.359-360]{Temam-1997}
   From the sub-exponential behavior of singular value function (see \eqref{wsubexp})
   on a compact invariant set $K$ it follows that
  \begin{equation}\label{1twinflim}
   \inf\limits_{t > 0}
   \big(\max\limits_{u \in K}\omega_{d}(D\varphi^t(u))\big)^{1/t} =
   \lim\limits_{t \to +\infty} 
   \big(\max\limits_{u \in K}\omega_{d}(D\varphi^t(u))\big)^{1/t}, \ t>0.
\end{equation}
\end{lemma}
{\bf Proof.}
  The proof of this result follows from Fekete's lemma 
  for the subadditive functions \cite[pp.463-464]{Kuczma-2009}
  \footnote{
    If $f : \mathbb{R}^n \to \mathbb{R}$
    is a measurable subadditive function, then for every $u \in \mathbb{R}^n$
    there exists the limit $\lim\limits_{t\to\infty}\frac{f(tx)}{t}$.
  }.
$\blacksquare$

\begin{corollary}
  If $\omega_{d}(D\varphi^t(u)>0$, then
  \begin{equation}\label{1tlnwinflim}
    \inf\limits_{t > 0}\max_{u \in K}
    \frac{1}{t} \ln\big(\omega_{d}(D\varphi^t(u))\big) =
    \lim\limits_{t \to +\infty}\max_{u \in K}
    \frac{1}{t} \ln
    \big(\omega_{d}(D\varphi^t(u))\big).
 \end{equation}
\end{corollary}

\medskip

For a compact set $K$, $t>0$, $u \in K$, and $d \in [0,n]$
we consider two scalar functions $g_d(t,u)$ and $f_d(t,u)$.
Suppose that $g_0(t,u) = f_0(t,u)\equiv c$,
therefore the following expressions
\[
  d_g^{+}(t,u) = \sup\{d \in [0,n]: g_d(t,u)\geq c \},
  \quad
  d_f^{+}(t,u) = \sup\{d \in [0,n]: f_d(t,u)\geq c \}
\]
are well defined. Also we consider
\[
  d_g^{-}(t,u) = \inf\{d \in [0,n]: g_d(t,u)< c \},
  \quad
  d_f^{-}(t,u) = \inf\{d \in [0,n]: f_d(t,u)< c \}.
\]
Here and further if the infimum on the empty set is considered,
then we assume that the infimum is equal $n$.
Define
\[
  d_f^{+}(t)=\sup\{d \in [0,n]: \sup_{u \in K}f_d(t,u) \geq c \},
  \quad
  d_f^{-}(t)=\inf\{d \in [0,n]: \sup_{u \in K}f_d(t,u)<c \}.
\]

\begin{lemma}\label{thm:fg}
We have the following properties

\begin{description}
\item[P1] 
  If for fixed $t>0$ and $u \in K$
  the implication $\big(g_d(t,u)<c \Rightarrow f_d(t,u)<c\big)$ holds
  $\forall d \in [0,n]$, then
  \begin{equation}\label{wp1}
    \inf\{d \in [0,n]: f_d(t,u)<c \}
    \leq
    \inf\{d \in [0,n]: g_d(t,u)<c \};
  \end{equation}
\item[P2] 
    If for fixed $t>0$ and $u \in K$
    the inequality $f_d(t,u) \leq g_d(t,u)$
    holds $\forall d \in [0,n]$,
    then
  \begin{equation}\label{wp2}
      \inf\{d \in [0,n]: f_d(t,u)<c \}
      \leq
      \inf\{d \in [0,n]: g_d(t,u)<c \};
  \end{equation}
  \begin{equation}\label{wp2sup}
      \sup\{d \in [0,n]: f_d(t,u) \geq c \}
      \leq
      \sup\{d \in [0,n]: g_d(t,u) \geq c \};
  \end{equation}
\item[P3] 
  If
  \begin{equation}\label{wp3cond2}
    \sup\{d \in [0,n]: f_d(t,u)\geq c \}
    = \inf\{d \in [0,n]: f_d(t,u) <c \},
  \end{equation}
  then
  \begin{equation}\label{wp3}
      \sup\{d \in [0,n]: \sup_{u \in K}f_d(t,u) \geq c \}
      =
      \sup_{u \in K}\sup\{d \in [0,n]: f_d(t,u)\geq c \}
      ;
  \end{equation}
\item[P4] 
    If for fixed $t>0$ the equality
    \begin{equation}\label{wp4cond1}
      \sup_{u \in K}f_d(t,u) = \max_{u \in K}f_d(t,u) \quad \forall d \in [0,n]
    \end{equation}
    is valid and \eqref{wp3cond2} holds,
    then
  \begin{equation}\label{wp4}
      \inf\{d \in [0,n]: \sup_{u \in K}f_d(t,u)<c \}
      =
      \sup_{u \in K}\inf\{d \in [0,n]: f_d(t,u)<c \}
      ;
  \end{equation}
 \item[P5]
   \begin{equation}\label{wp5}
     \inf\{d \in [0,n]: \inf_{t>0} f_d(t,u)<c \} =
     \inf_{t>0}\inf\{d \in [0,n]: f_d(t,u)<c \}.
   \end{equation}
\end{description}
\end{lemma}
{\bf Proof.}

{\bf (P1), (P2)}:
Since in (P1) and (P2)
the set of possible $d$, considered in the left-hand side of expression,
involves the set of possible $d$, considered in the right-hand side of expression,
we have the corresponding inequalities for the infimums of the sets.
Similarly we get relation for supremums.

{\bf (P3)}:
Since $f_d(t,u) \leq \sup_{u \in K} f_d(t,u)$,
by \eqref{wp2sup} in (P2)
we have $\sup_{u \in K}d_f^{+}(t,u) \leq d_f^{+}(t)$.

\noindent
Let $\sup_{u \in K}d_f^{+}(t,u) < d_f^{+}(t)$
$\Rightarrow$ by \eqref{wp3cond2}
$\exists d_0 \in \big(\sup_{u \in K}d_f^{+}(t,u), d_f^{+}(t)\big):
f_{d'}(t,u)<c$   $\forall d' \in [d_0,n]\ \forall u \in K$
$\Rightarrow$ $d_f^{+}(t) \leq d_0$.
Thus we get the contradiction.

{\bf (P4)}:
Since $\sup_{u \in K}f_d(t,u)<c$ implies $f_d(t,u)<c$ for all $u \in K$,
by (P1) we have $d_f^{-}(t,u) \leq d_f^{-}(t)$ for all $u \in K$
$\Rightarrow$ $\sup_{u \in K}d_f^{-}(t,u) \leq d_f^{-}(t)$.

Let $\sup_{u \in K}d_f^{-}(t,u) < d_f^{-}(t)$.
Then $\exists d_0 \in \big(\sup_{u \in K}d_f^{-}(t,u),\,d_f^{-}(t)\big)$.
Since $d_0 < d_f^{-}(t)$, we have $\sup_{u \in K} f_{d_0}(t,u)\geq c.$
Therefore, from condition \eqref{wp4cond1}, $\exists u_0: f_{d_0}(t,u_0)\geq c$.
Finally, according to condition \eqref{wp3cond2},
we have $d_0 \leq \sup\{d \in [0,n]: f_d(t,u_0) \geq c\}=\inf\{d \in [0,n]: f_d(t,u_0)<c \}
=d_f^{-}(t,u_0)
\leq \sup_{u \in K}d_f^{-}(t,u)$.
Thus we get the contradiction.

{\bf (P5)}:
Since $\inf_{t>0} f_d(t,u) \leq f_d(t,u)$, by \eqref{wp2} from (P2)
we have $d_f^{-}(u) \leq d_f^{-}(t,u)$
and, thus, $d_f^{-}(u) \leq \inf_{t>0}d_f^{-}(t,u)$.

\noindent
Let $d_f^{-}(u) < \inf_{t>0}d_f^{-}(t,u)$
$\Rightarrow$
$\exists d_0 \in \big[d_f^{-}(u), \inf_{t>0}d_f(t,u)\big): \inf_{t>0}f_{d_0}(t,u)<c$
$\Rightarrow$
$\exists t_0: f_{d_0}(t_0,u)<c$
$\Rightarrow$
$d_0\geq d_f^{-}(t_0,u) \geq \inf_{t>0}d_f^{-}(t,u)$.
Thus we get the contradiction.
$\blacksquare$


\begin{theorem}\label{thm:ucrit}
\cite[p.147, eq.3.21]{EdenFT-1991},\cite[p.112, eq.4.19]{Eden-1990}
Let $\omega_{d}(D\varphi^t(u))>0$.
For a compact invariant set $K$
and $d \in [0,n]$ there is a point $u^{cr}=u^{cr}(d) \in K$ (it may be not unique)
such that
\begin{equation}\label{finiteLEgglobalLE}
    \frac{1}{t}\ln\omega_{d}(D\varphi^t(u^{cr}(d)))
    \geq
    \limsup\limits_{t \to +\infty}\sup\limits_{u \in K}
    \frac{1}{t}\ln\omega_{d}(D\varphi^t(u))
    \quad \forall t>0.
\end{equation}
\end{theorem}
Relation \eqref{finiteLEgglobalLE} is presented in
\cite[p.147, eq.3.21]{EdenFT-1991},\cite[pp.114, eq.5.6]{Eden-1990}
and its proof is based on the theory of positive operators \cite{ChoquetF-1975}
(see also \cite{GundlachS-2000}).
\begin{corollary} (see, e.g. \cite[pp.113-114]{Eden-1990})
\begin{equation}\label{suplim=limsup}
\begin{aligned}
&
    \sup\limits_{u \in K}
    \limsup\limits_{t \to +\infty}
    \frac{1}{t}\ln\omega_{d}(D\varphi^t(u))
    &=&
    \lim\limits_{t \to +\infty}\frac{1}{t}\ln\omega_{d}(D\varphi^t(u^{cr}(d)))=
    \\
    && =&
    \max\limits_{u \in K}
    \limsup\limits_{t \to +\infty}
    \frac{1}{t}\ln\omega_{d}(D\varphi^t(u))
    =
    \limsup\limits_{t \to +\infty}\sup\limits_{u \in K}
    \frac{1}{t}\ln\omega_{d}(D\varphi^t(u)).
\end{aligned}
\end{equation}
\end{corollary}
{\bf Proof.}
It is easy to check that (see, e.g. \cite[p.31]{ConstantinFT-1985})
\begin{equation}\label{suplim-limsup}
  \sup\limits_{u \in K}\limsup\limits_{t \to +\infty}\frac{1}{t}\ln\omega_{d}(D\varphi^t(u))
  \leq
  \limsup\limits_{t \to +\infty}\sup\limits_{u \in K}\frac{1}{t}\ln\omega_{d}(D\varphi^t(u)).
\end{equation}
Thus, taking into account \eqref{finiteLEgglobalLE}, we get \eqref{suplim=limsup}.
$\blacksquare$

\section{Lyapunov dimension of maps}
The concept of the Lyapunov dimension had been suggested in the seminal paper
by Kaplan and Yorke \cite{KaplanY-1979} and
later it was rigorously developed in a number of papers
(see, e.g. \cite{FredericksonKYY-1983,ConstantinFT-1985}).

The following two definitions are inspirited by Douady--Oesterl\'{e} \cite{DouadyO-1980}.
\begin{definition}\label{def:DOlocal-map}
The (local) Lyapunov dimension\footnote{This is not a dimension
in a rigorous sense (see, e.g. \cite{HurewiczW-1941,Kuratowski-1966,AleksandrovP-1973})}
of a continuously differentiable map $\varphi: U \subseteq \mathbb{R}^n \to \mathbb{R}^n$
at the point $u \in U$ is defined as
\[
  d_{\rm L}(\varphi,u) := \sup\{d \in [0,n]: \omega_{d}(D\varphi(u)) \geq 1 \}.
\]
\end{definition}
For any $u \in U$ this value is well-defined
since $\omega_{0}(D\varphi(u)) \equiv 1$.

By Lemma~\ref{thm:wleftcont} we get
\begin{equation}\label{supDmax}
  d_{\rm L}(\varphi,u) =
  \max\{d \in [0,n]:\ \omega_{d}(D\varphi(u)) \geq 1 \}.
\end{equation}
Additionally, since the singular values in \eqref{defomega} are ordered by decreasing,
we have
\begin{equation}\label{maplocalDOinf}
  d_{\rm L}(\varphi,u)=\max\{d \in [0,n]:\ \omega_{d}(D\varphi(u)) \geq 1 \} =
  \inf\{d \in [0,n]:\ \omega_{d}(D\varphi(u)) < 1 \}
\end{equation}
if the infimum exists (i.e. there exists
$d \in (0,n]$ such that $\omega_{d}(D\varphi(u))< 1$).
Here and further in the similar constructions
if the infimum does not exist,
we assume that the infimum and considered dimension
are taken equal to $n$.

\begin{definition}\label{def:DO-map}
  The Lyapunov dimension of
  a continuously differentiable map $\varphi: U \subseteq \mathbb{R}^n \to \mathbb{R}^n$
  of the compact set $K \subset U \subseteq \mathbb{R}^n$
  is defined as
 \[
    d_{\rm L}(\varphi,K) := \sup\limits_{u \in K} d_{\rm L}(\varphi,u)
    = \sup\limits_{u \in K} \sup\{d \in [0,n]: \omega_{d}(D\varphi(u)) \geq 1 \}.
 \]
\end{definition}


Remark that by Lemma~\ref{thm:fg}\,(property~\eqref{wp3}) and Lemma~\ref{thm:wcontu}
we have
\begin{equation}\label{supUsup=supsupU}
\begin{aligned}
    & d_{\rm L}(\varphi,K)
    = \sup\limits_{u \in K} \sup\{d \in [0,n]: \omega_{d}(D\varphi(u)) \geq 1 \} =
    \sup\{d \in [0,n]: \max\limits_{u \in K}\omega_{d}(D\varphi(u)) \geq 1 \}.
\end{aligned}
\end{equation}
Additionally, by \eqref{maplocalDOinf}
and Lemma~\ref{thm:fg}\,(property~\eqref{wp4}), we have
\begin{equation}\label{mapDOsupinf}
  d_{\rm L}(\varphi,K)=
  \sup\limits_{u \in K}
  \inf\{d \in [0,n]:\ \omega_{d}(D\varphi(u)) < 1 \} =
  \inf\{d \in [0,n]:\ \max\limits_{u \in K}\omega_{d}(D\varphi(u)) < 1 \}
\end{equation}
if the infimum exists (i.e. there exists $d \in (0,n]$ such that
$\max\limits_{u \in K}\omega_{d}(D\varphi(u)) < 1$).

\begin{theorem}\label{thm:DO} (Douady--Oesterl\'{e}, \cite[p.1135]{DouadyO-1980};
see also \cite[p.369]{Temam-1997},\cite[p.239]{Smith-1986},\cite[p.332]{BoichenkoLR-2005})
If the continuously differentiable map $\varphi: U \subseteq \mathbb{R}^n \to \mathbb{R}^n$
has a negatively invariant or invariant compact set $K \subset U$, i.e.
 \[
   \varphi(K) \supseteq K,
 \]
 then
 \[
  \begin{aligned}
   & \dim_{\rm H} K \leq d_{\rm L}(\varphi,K).
   \end{aligned}
 \]
\end{theorem}

Remark that under the assumptions of Theorem~\ref{thm:DO}
if $\omega_{d}(D\varphi(u))<1$ for some $d \leq 1$,
then $\dim_{\rm H} K = 0$ (see, e.g. \cite[p.371]{Temam-1997}).
Thus, taking into account Lemma~\ref{thm:wcontu}, we have
\begin{lemma}(see, e.g. \cite[p.554]{Gelfert-2003})
  The functions $u \mapsto d_{\rm L}(\varphi,u)$ is continuous on $U$
  except at a point $u$, which satisfies $\sigma_1(D\varphi(u)) = 1$,
  where it is still upper semi-continuous.
\end{lemma}
\begin{corollary}
By the Weierstrass extreme value theorem for the upper semi-continuous functions, 
there exists a critical point $u_{\rm L}$ (it may be not unique) such that 
\begin{equation}\label{DOTmax}
   \sup_{u \in K}d_{\rm L}(\varphi,u)
   =
   \max_{u \in K}d_{\rm L}(\varphi,u)
   =
   d_{\rm L}(\varphi,u_{\rm L}).
\end{equation}
\end{corollary}

For an invariant compact set $K$ of the dynamical system
$\big(\{\varphi^t\}_{t\geq0},(U\subseteq \mathbb{R}^n,||\cdot||) \big)$
one may consider for a fixed $t$ the evolutionary operator $\varphi^t(u)$, then
\[
  \varphi^t(K)=K
\]
and the corresponding Lyapunov dimension (\emph{finite time Lyapunov dimension})
\begin{equation}\label{DOmapt}
\begin{aligned}
  &
  d_{\rm L}(\varphi^t,K) &=& \sup\limits_{u \in K} d_{\rm L}(\varphi^t,u)=
  \inf\{d \in [0,n]: \max\limits_{u \in K}\omega_{d}(D\varphi^t(u))<1\}.
\end{aligned}
\end{equation}

\begin{example}
If for a nonempty compact set $K \subset U \subseteq \mathbb{R}^n$
it is considered the identical map $\varphi={\rm id}$,
then $D\varphi(u)= I$ and by the definition of the Lyapunov dimension we have
\(
  d_{\rm L}({\rm id},K)=n.
\)
Remark that for $t=0$  we have $\varphi^0={\rm id}$
and $d_{\rm L}(\varphi^0,K)=n$, thus we further consider $t>0 $.
\end{example}

\begin{remark}
   For the numerical estimations of dimension, the following remark is important:
   for any $t>0$ the equality \eqref{winfliminf}
   for a compact invariant set $K$ implies
   the existence of $s=s(t)>0$ such that
   \begin{equation}\label{DOinctT}
     d_{\rm L}(\varphi^{t+s},K) \leq d_{\rm L}(\varphi^t,K).
   \end{equation}
\end{remark}

While in the computations we can consider only finite time $t \leq T$
and evolutionary operator $\varphi^T(u)$,
from a theoretical point of view,
it is interesting to study
the limit behavior of dynamical system $\{\varphi^t\}_{t\geq0}$ as $t \to +\infty$.
Next, unless otherwise stated,
$K \subset U \subseteq \mathbb{R}^n$ denotes a \emph{compact invariant set}
with respect to
the dynamical system
$\big(\{\varphi^t\}_{t\geq0},(U\subseteq \mathbb{R}^n,||\cdot||) \big)$:
$\varphi^t(K)\!=\!K,\ t\!>\!0.$

\section{Lyapunov dimensions of dynamical system}
According to the Douady-Oesterl\'{e} theorem
it is natural to give the following generalization of Definition~\ref{def:DO-map}
for dynamical systems.
\begin{definition}\label{defdimLt}
The Lyapunov dimension of the dynamical system $\{\varphi^t\}_{t\geq0}$
with respect to a compact invariant set $K$ is defined as
\begin{equation}\label{DOds}
  d_{\rm L}(\{\varphi^t\}_{t\geq0},K) :=
  \inf_{t > 0}d_{\rm L}(\varphi^t,K) =
  \inf_{t > 0}\sup\{d \in [0,n]: \max\limits_{u \in K}\omega_{d}(D\varphi^t(u)) \geq 1 \}.
\end{equation}
\end{definition}

By Theorem~\ref{thm:DO} we have
\begin{equation}\label{dimHKlessDOK}
  \dim_{\rm H} K \leq d_{\rm L}(\{\varphi^t\}_{t\geq0},K) \leq d_{\rm L}(\varphi^t,K).
\end{equation}
By \eqref{DOmapt} and Lemma~\ref{thm:fg}\,(property~\eqref{wp5})
 we have\footnote{
While inf and sup give the same values
for $\omega_{d}(D\varphi^t(u))$ in \eqref{supUsup=supsupU} and \eqref{mapDOsupinf},
for $\inf_{t > 0}\max\limits_{u \in K}\omega_{d}(D\varphi^t(u))$ we need consider

$\sup\{d \in [0,n]:
\forall \widetilde{d} \in [0,d] \
\inf_{t > 0}\max\limits_{u \in K}\omega_{\widetilde{d}}(D\varphi^t(u))
\geq 1
\}
=
\inf\{d \in [0,n]:
\inf_{t > 0}\max\limits_{u \in K}\omega_{d}(D\varphi^t(u))
< 1
\}.
$
}
\begin{equation}\label{DOdsDOmapinf}
  d_{\rm L}(\{\varphi^t\}_{t\geq0},K)=
  \inf_{t > 0}d_{\rm L}(\varphi^t,K) =
  \inf\{d \in [0,n]: \inf_{t > 0}\max\limits_{u \in K}\omega_{d}(D\varphi^t(u)) < 1 \}
\end{equation}
and, finally, by \eqref{winfliminf} we have (see also \cite[p.65]{Eden-1989-PhD})
\begin{equation}\label{DOdsDOmap}
  d_{\rm L}(\{\varphi^t\}_{t\geq0},K)=
  \inf_{t > 0}d_{\rm L}(\varphi^t,K) =
  \inf\{d \in [0,n]: \liminf_{t \to +\infty}\max\limits_{u \in K}\omega_{d}(D\varphi^t(u)) = 0\}.
\end{equation}

It is interesting to consider a critical point $u_{\rm L}(T) \in K$
such that the supremum of the local finite time Lyapunov dimension
$d_{\rm L}(\varphi^T,u)$
is achieved at this point\footnote{
  Let there exists
  $\lim_{t \to +\infty}\max\limits_{u \in K}\omega_{d}(D\varphi^t(u))=0$.
  It is interesting to study
  1) the existence of critical point $u_0 \in K$ such that
  $
    \lim_{t \to +\infty}\max\limits_{u \in K}\omega_{d}(D\varphi^t(u))
    =
    \limsup_{t \to +\infty}\omega_{d}(D\varphi^t(u_0))
  $
  and 2) the estimations
  $\dim_{H}K\, ?\!\leq
  \inf\{d \in [0,n]: \sup\limits_{u \in K}\limsup_{t \to +\infty}\omega_{d}(D\varphi^t(u)) < 1\}$
  or
  $\dim_{H}K\, ?\!\leq
  \sup\limits_{u \in K}\limsup_{t \to +\infty}\sup\{d \in [0,n]: \omega_{d}(D\varphi^t(u)) \geq 1\}$.
  Remark, it is clear that
  $
  \lim_{t \to +\infty}\max\limits_{u \in K}\omega_{d}(D\varphi^t(u))
  \geq
  \sup\limits_{u \in K}\limsup_{t \to +\infty}\omega_{d}(D\varphi^t(u))
  $.

  From \eqref{DOTmax} it follows the existence of a critical point $u_{\rm L}(t)$ such that
  $
   d_{\rm L}(\varphi^t,u_{\rm L}(t))
   =
   \max_{u \in K}d_{\rm L}(\varphi^t,u).
  $
  Taking into account \eqref{DOinctT}
  we can consider a sequence $t_k \to +\infty$ such that
  $d_{\rm L}(\varphi^{t_k},u_{\rm L}(t_k))$
  is monotonically decreasing to $\inf_{t\geq0}\max_{u \in K}d_{\rm L}(\varphi^t,u)$.
  Since $K$ is a compact set, we can consider
  a subsequence $t_m = t_{k_m} \to +\infty$ such that
  there exists a limit critical point $u^{cr}_{\rm L}$:
  $u_{\rm L}(t_m) \to u^{cr}_{\rm L} \in K$ as $t_m \to +\infty$.
  Thus we have
  $d_{\rm L}(\varphi^{t_m},u_{\rm L}(t_m)) \searrow
  d_{\rm L}(\{\varphi^t\}_{t\geq0},K)$ and $u_{\rm L}(t_m) \to u^{cr}_{\rm L} \in K$
  as $m \to +\infty$.
  One may guess that $u^{cr}_{\rm L}$ coincides with  
  a critical point $u^{\rm E}_{\rm L}$ from \eqref{ucrL}.
}


\begin{proposition}\label{thm:DOeqpoint}
Suppose that for a certain $t=T>0$
the supremum of the local finite time Lyapunov dimensions
$d_{\rm L}(\varphi^T,u)$
is achieved at one of the equilibria points:
\begin{equation}\label{DOeqpoint}
   d_{\rm L}(\varphi^T,u^{cr}_{eq}) =
   \sup_{u \in K}d_{\rm L}(\varphi^T,u),
   \quad
   \varphi^t(u^{cr}_{eq})\equiv u^{cr}_{eq}.
\end{equation}
Then
\begin{equation}\label{DOeqpoint}
   \dim_{\rm H} K \leq
   d_{\rm L}(\varphi^T,u^{cr}_{eq}) = d_{\rm L}(\{\varphi^t\}_{t\geq0},K)
   = d_{\rm L}(\varphi^{T},K).
\end{equation}
\end{proposition}
{\bf Proof.}
From \eqref{weqpoint} we have
\[
  \omega_{d}(D\varphi^{t}(u^{cr}_{eq}))
  = \big(\omega_{d}(D\varphi(u^{cr}_{eq}))\big)^{t}, \quad t \geq 0.
\]
Therefore
\[
  \big(\omega_{d}(D\varphi^{T}(u^{cr}_{eq})) < 1\big)
  \Leftrightarrow
  \big(\omega_{d}(D\varphi(u^{cr}_{eq})) < 1\big)
  \Leftrightarrow
  \big(\omega_{d}(D\varphi^{t}(u^{cr}_{eq})) < 1,\ \forall t>0 \big)
\]
and
\[
\]
\[
\big(\liminf\limits_{t \to +\infty}\max\limits_{u \in K}\omega_{d}(D\varphi^t(u)) < 1\big)
\Rightarrow
\big(\liminf\limits_{t \to +\infty}\omega_{d}(D\varphi^t(u^{cr}_{eq}))<1\big)
\Leftrightarrow
\big(\omega_{d}(D\varphi^T(u^{cr}_{eq}))<1\big).
\]
By Lemma~\ref{thm:fg}\,(property~\eqref{wp1}) we obtain
\[
  \inf\{d \in [0,n]: \omega_{d}(D\varphi^T(u^{cr}_{eq})) < 1\}
  = d_{\rm L}(\varphi^T,u^{cr}_{eq}) \leq d_{\rm L}(\{\varphi^t\}_{t\geq0},K).
\]
Finally, by from \eqref{dimHKlessDOK}
we get the assertion of the proposition.
$\blacksquare$

Further, to consider $\ln\omega_{d}(D\varphi^t(u))$,
we suppose that $\det J(u)\neq0$ $\forall u \in U$ and thus
\begin{equation}\label{sn>0}
  \sigma_i(t,u) > 0, \quad i=1,\ldots,n.
\end{equation}

The following definitions of Lyapunov dimension
are inspirited by Constantin, Foias, Temam \cite[p.31,Remark 3.1., ii)]{ConstantinFT-1985}
and Eden \cite[p.114]{Eden-1990}
\footnote{
In \cite[p.31,Remark 3.1., ii)]{ConstantinFT-1985}
Constantin, Foias, Temam stated that
if
$
    \sup_{u \in K}\limsup_{t\to+\infty}
    \big(\omega_{d}(D\varphi^t(u))\big)^{1/t} < 1
$
or
$
    \limsup_{t\to+\infty}
    \sup_{u \in K}\big(\omega_{d}(D\varphi^t(u))\big)^{1/t}<1,
$
then $\dim_{\rm H}(K)\leq d$.
In \cite[p.114]{Eden-1990} Eden considered the value
$d_{\rm O}(K) =
\inf\{ d>0:\, \sup_{u \in K} \omega_{d}(D\varphi^t(u))
\mbox{\ converges to zero exponentially as\ }
t \to \infty
\}
$
and called it \emph{the Douady-Oesterl\'{e} dimension of $K$}.}.

\begin{definition}\label{defdimLt}
The (global) Lyapunov dimension of the dynamical system $\{\varphi^t\}_{t\geq0}$
with respect to a compact invariant set $K$ is defined as\footnote{
Comparing the expressions in the definitions \eqref{DOds}
and \eqref{dimLt},
remark that we can change $\frac{1}{t}$ in \eqref{dimLt}
to another scalar positive monotonically decreasing function $q(t)$
such that
$\inf_{t>0}q(t)\max_{u \in K}\omega_d(D\varphi^{t}(u))=
\lim_{t\to+\infty}q(t)\max_{u \in K}\omega_d(D\varphi^{t}(u))$.
The last relation is important from a computational point of view.
}
\begin{equation}\label{dimLt}
  d_{\rm L}^{\rm E}(\{\varphi^t\}_{t\geq0},K):=
  \,\inf\{d \in [0,n]:
  \lim\limits_{t \to +\infty}\max\limits_{u \in K}\frac{1}{t}\ln\omega_{d}(D\varphi^t(u))<0\}.
\end{equation}
\end{definition}
Correctness of the definition follows from \eqref{1tlnwinflim}.

By \eqref{sn>0} we have
\(
   \big(
   \liminf\limits_{t \to +\infty}\max_{u \in K}
   (\omega_{d}(D\varphi^t(u)))<1
   \Leftrightarrow
   \liminf\limits_{t \to +\infty}\max_{u \in K}
   \ln(\omega_{d}(D\varphi^t(u)))<0
   \big)
\)
and
\(
  \bigg(
   \liminf\limits_{t \to +\infty}\max_{u \in K}
   \frac{1}{t} \ln
   \big(\omega_{d}(D\varphi^t(u))\big)
   <0
   \Rightarrow
   \liminf\limits_{t \to +\infty}\max_{u \in K}
   \ln(\omega_{d}(D\varphi^t(u)))<0
   \bigg).
\)
Thus, taking into account \eqref{DOdsDOmap} and \eqref{1tlnwinflim},
by Lemma~\ref{thm:fg}\,(property~\eqref{wp1})
we have
\begin{equation}\label{dimHKdDOKdL}
\begin{aligned}
  &
  d_{\rm L}(\{\varphi^t\}_{t\geq0},K)
  &=& \inf\{d \in [0,n]: \liminf\limits_{t \to +\infty}\max\limits_{u \in K}\ln\omega_{d}(D\varphi^t(u))<0\}
  \leq
  \\
  && \leq &
  \inf\{d \in [0,n]:
  \lim\limits_{t \to +\infty}\max\limits_{u \in K}\frac{1}{t}\ln\omega_{d}(D\varphi^t(u))<0\}
  =
  d_{\rm L}^{\rm E}(\{\varphi^t\}_{t\geq0},K).
\end{aligned}
\end{equation}
Since for fixed $t>0$ and $d \in [0,n]$
\[
  \bigg( 1>\max_{u \in K}\big(\omega_{d}(D\varphi^t(u))\big) \bigg)
   \Rightarrow
  \bigg( 0>\max_{u \in K} \frac{1}{t} \ln\omega_{d}(D\varphi^t(u))
   \geq
   \lim\limits_{t \to +\infty} \max_{u \in K}\frac{1}{t} \ln\omega_{d}(D\varphi^t(u)) \bigg),
\]
by Lemma~\ref{thm:fg}\,(property~\eqref{wp1}) and \eqref{dimHKlessDOK} we have
\begin{proposition}
\begin{equation}\label{allglobaldim}
\begin{aligned}
  \dim_{\rm H}(K) \leq& d_{\rm L}(\{\varphi^t\}_{t\geq0},K)
  \leq
  d_{\rm L}^{\rm E}(\{\varphi^t\}_{t\geq0},K)
  \leq
  d_{\rm L}(\varphi^{t},K)\ \ \ \forall t>0.
\end{aligned}
\end{equation}
\end{proposition}
\begin{corollary}
Taking $\inf_{t>0}$ in \eqref{allglobaldim}, we obtain
\begin{equation}\label{DO=DL}
  d_{\rm L}(\{\varphi^t\}_{t\geq0},K) =
  d_{\rm L}^{\rm E}(\{\varphi^t\}_{t\geq0},K).
\end{equation}
\end{corollary}

\begin{definition}\label{deflocalDL}
The local Lyapunov dimension of the dynamical system $\{\varphi^t\}_{t\geq0}$
at the point $u$ is defined as
\begin{equation}\label{dimL}
  d_{\rm L}^{\rm E}(\{\varphi^t\}_{t\geq0},u):=
  \inf\{d \in [0,n]: \limsup\limits_{t \to +\infty}\frac{1}{t}\ln\omega_{d}(D\varphi^t(u))<0 \}.
\end{equation}
\end{definition}

By \eqref{suplim=limsup} and Lemma~\ref{thm:fg}\,(property~\eqref{wp4}) we have
\begin{equation}\label{localLD=globalLD}
\begin{aligned}
  \sup\limits_{u \in K}\inf\{d \in [0,n]:
  \limsup\limits_{t \to +\infty}\frac{1}{t}\ln\omega_{d}(D\varphi^t(u))<0\}
  =
  \inf\{d \in [0,n]:
  \sup\limits_{u \in K}\limsup\limits_{t \to +\infty}\frac{1}{t}\ln\omega_{d}(D\varphi^t(u))<0\}
  \end{aligned}
\end{equation}
Therefore, by \eqref{suplim-limsup} and Lemma~\ref{thm:fg}\,(property~\eqref{wp2})
we get
\begin{equation}\label{localLD=globalLD}
\begin{aligned}
  &
  \sup\limits_{u \in K}d_{\rm L}^{\rm E}(\{\varphi^t\}_{t\geq0},u)
  &=&
  \sup\limits_{u \in K}\inf\{d \in [0,n]:
  \limsup\limits_{t \to +\infty}\frac{1}{t}\ln\omega_{d}(D\varphi^t(u))<0\}
  =
  \\
  &&=&
  \inf\{d \in [0,n]:
  \sup\limits_{u \in K}\limsup\limits_{t \to +\infty}\frac{1}{t}\ln\omega_{d}(D\varphi^t(u))<0\}
  \leq
  \\
  &&\leq&
  \inf\{d \in [0,n]:
  \limsup\limits_{t \to +\infty}\sup\limits_{u \in K}\frac{1}{t}\ln\omega_{d}(D\varphi^t(u))<0\}
  =
  d_{\rm L}^{\rm E}(\{\varphi^t\}_{t\geq0},K).
\end{aligned}
\end{equation}

\begin{proposition}\label{thm:DODLeqpoint}
If there is a critical equilibrium point $u^{cr}_{eq}$
such that \eqref{DOeqpoint} is valid, then
\[
  d_{\rm L}(\varphi^{T},u^{cr}_{eq}) =
  d_{\rm L}^{\rm E}(\{\varphi^t\}_{t\geq0},u^{cr}_{eq})
  =
  \sup_{u \in K}d_{\rm L}^{\rm E}(\{\varphi^t\}_{t\geq0},K)
\]
and from \eqref{allglobaldim} it follows
 \begin{equation}\label{alldimeq}
\begin{aligned}
  \dim_{\rm H}(K) \leq &
  d_{\rm L}(\{\varphi^t\}_{t\geq0},K) =
  d_{\rm L}^{\rm E}(\{\varphi^t\}_{t\geq0},K) = d_{\rm L}^{\rm E}(\{\varphi^t\}_{t\geq0},u^{cr}_{eq})
  =
  d_{\rm L}(\varphi^{T},K).
\end{aligned}
\end{equation}
\end{proposition}
In this case for the estimation of the Hausdorff dimension
by \eqref{alldimeq} we need only the Douady-Oesterl\'{e} theorem (see Theorem~\ref{thm:DO}).
In the general case the existence of a critical point $u^{\rm E}_{\rm L}$
(it may be not unique) such that
\begin{equation}\label{ucrL}
  d_{\rm L}^{\rm E}(\{\varphi^t\}_{t\geq0},u^{\rm E}_{\rm L})
  = \sup\limits_{u \in K}d_{\rm L}^{\rm E}(\{\varphi^t\}_{t\geq0},u)
  =  d_{\rm L}^{\rm E}(\{\varphi^t\}_{t\geq0},K)
\end{equation}
follows from \eqref{suplim=limsup}
and the so-called \emph{Eden conjecture} is that $u^{\rm E}_{\rm L}$
corresponds to an equilibrium point or to a periodic orbit
(\cite[p.98, Question 1.]{Eden-1989-PhD}).

Finally, from \eqref{allglobaldim} and \eqref{localLD=globalLD} we have
\begin{theorem}\label{thm:alldim}
\begin{equation}\label{alldim}
\begin{aligned}
  \dim_{\rm H}(K) \leq &
  d_{\rm L}(\{\varphi^t\}_{t\geq0},K)
  =
  d_{\rm L}^{\rm E}(\{\varphi^t\}_{t\geq0},K)= \sup_{u\in K} d_{\rm L}^{\rm E}(\{\varphi^t\}_{t\geq0},u)
  \leq
  d_{\rm L}(\varphi^{t},K)=\sup\limits_{u \in K}d_{\rm L}(\varphi^t,u).
\end{aligned}
\end{equation}
\end{theorem}

\subsection{Lyapunov exponents: various definitions}
\begin{definition}\label{defLE}
The Lyapunov exponent functions of singular values
(also called finite-time Lyapunov exponents \cite{AbarbanelBK-1991})
of the dynamical system
$\big(\{\varphi^t\}_{t\geq0},(U\subseteq \mathbb{R}^n,||\cdot||) \big)$
at the point $u \in U$ are denoted by
\[
\begin{aligned}
   & \LEs_i(t,u) = \LEs_i(D\varphi^t(u)),\quad \ i=1,2,...,n,
   \\ &
   \LEs_1(t,u) \geq\cdots\geq\LEs_n(t,u), \ \forall t> 0,
\end{aligned}
\]
and defined as
\[
  \LEs_i(t,u) := \frac{1}{t}\ln\sigma_i(t,u).
\]
\end{definition}

\begin{definition}\label{def:lLE}
The \emph{Lyapunov exponents (LEs) of singular values}\footnote{
We add ``of singular value''
to distinguish this definition
from other definitions of Lyapunov exponents;
if the differences in the definitions
are not significant for the presentation,
we use the term "Lyapunov exponents" or "LEs".
}
of the dynamical system $\{\varphi^t\}_{t\geq0}$
at the point $u$ are defined (\alarm{see, e.g. \cite{Oseledec-1968}},\cite[p.29,eq.3.26]{ConstantinFT-1985})
as
\[
   \LEs_i(u) :=
   \limsup\limits_{t \to +\infty} \LEs_i(t,u) =
   \limsup\limits_{t \to +\infty} \frac{1}{t}\ln(\sigma_i(t,u)), \quad i=1,2,..,n.
\]
\end{definition}

Often $\LEs_i(u)$ are called upper LEs and denoted as $\overline{\LEs}_i(u)$,
while $\underline{\LEs}_{\,i}(u) := \liminf\limits_{t \to +\infty}\LEs_{k}(t,u)$
are called lower LEs.
Remark that the Lyapunov exponents of singular values are
the same for any fundamental matrices
of the linearized systems \eqref{sflct} or \eqref{sfldt}
\begin{proposition}(see, e.g. \cite{KuznetsovAL-2014-arXiv-LE})\label{thm:LE-def-inv}
For the matrix $D\varphi^t(u)P$, where $P$ is a nonsingular
$n\times n$ matrix (i.e. $\det P \neq 0$), one has
 \[
   \lim\limits_{t \to +\infty} \bigg(\LEs_i(D\varphi^t(u)) - \LEs_i(D\varphi^t(u)P) \bigg) = 0,
   \quad i=1,2,...,n.
 \]
\end{proposition}

\begin{definition}\label{def:LCE}
The Lyapunov exponent functions of the fundamental matrix columns
$(y^1(t,u),...,y^n(t,u))=D\varphi^t(u)$
\[
   \LCE_i(t,u) = \LCE_i(D\varphi^t(u)), \ i=1,2,...,n, \ u \in U
\]
are defined
as
\[
  \LCE_i(t,u) := \frac{1}{t}\ln||y^i(t,u)||.
\]
The ordered Lyapunov exponent functions of the fundamental matrix columns
at the point $u$
(also called finite-time Lyapunov characteristic exponents)
are given by the ordered set (for all $t>0$) of $\LCE_i(t,u)$:
\[
   \LCE^{o}_1(t,u) \geq\cdots\geq\LCE^{o}_n(t,u), \ \forall t\geq 0.
\]
\end{definition}

\begin{definition}\label{def:LCE}
The Lyapunov exponents of the fundamental matrix columns\footnote{
Often they are called Lyapunov characteristic exponents (LCE) \cite{LeonovK-2007}.
In \cite{Lyapunov-1892} these values are defined
with the opposite sign and called \emph{characteristic exponents}
at the point $u$.
}
are defined (see \cite{Lyapunov-1892}) as
\[
   \LCE_i(u) :=
   \limsup\limits_{t \to +\infty} \LCE^{o}_i(t,u), \quad i=1,2,..,n.
\]
\end{definition}

\begin{remark}
  The Lyapunov exponents of fundamental matrix columns
  may be different for different fundamental matrices
  in contrast to the definition of Lyapunov exponents of singular values
  (see, e.g. Proposition~\ref{thm:LE-def-inv}).
  To get the set of all possible values of Lyapunov exponents of fundamental matrix columns
  (the set with the minimal sum of values),
  one has to consider the so-called normal fundamental matrices
  (see \cite{Lyapunov-1892},\cite{LeonovK-2007}).
\end{remark}

\begin{definition}\label{def:lvolLE}
The \emph{relative Lyapunov exponents of singular value functions}
of the dynamical system $\{\varphi^t\}_{t\geq0}$ at the point $u$
are defined (\alarm{see, e.g. \cite{Oseledec-1968}}) as
\[
 \begin{aligned}
 & \LEo_{1}(u) &:=& \limsup_{t \to +\infty}(\LEs_1(t,u)), \\
 &
 \LEo_{i+1}(u) &:=&
 \limsup_{t \to +\infty}(\LEs_1(t,u)+\cdots+\LEs_{i+1}(t,u))-
 \\&&-&\limsup_{t \to +\infty}(\LEs_1(t,u)+\cdots+\LEs_{i}(t,u)),\  i=1,...,n-1.
 \end{aligned}
\]
\end{definition}

For $k=1,2,..,n$ we have
\[
 \LEo_{1}(u)+\cdots+\LEo_{k}(u) = \limsup_{t \to +\infty}(\LEs_1(t,u)+\cdots+\LEs_{k}(t,u)),
\]
\[
 \begin{aligned}
 &
    \underline{\LEs}_{\,k}(u)
    \leq
    \LEo_{k}(u)
    &=&
    \limsup\limits_{t \to +\infty}\sum\limits_{i=1}^{k}\LEs_i(t,u)
    -
    \limsup\limits_{t \to +\infty}\sum\limits_{i=1}^{k-1}\LEs_i(t,u)
    \leq
  \\&&  \leq&
    \limsup\limits_{t \to +\infty}\LEs_{k}(t,u) = \LEs_{k}(u).
 \end{aligned}
\]

  From the Courant-Fischer theorem \cite{HornJ-1994-book} it follows (see, e.g. \cite{Barabanov-2005})\footnote{
 For example \cite{KuznetsovAL-2014-arXiv-LE},
 for the matrix
 \(
    u(t)=\left(
      \begin{array}{cc}
        1 & g(t)-g^{-1}(t) \\
        0 & 1 \\
      \end{array}
    \right)
 \) we have the following ordered values:
 $  \LCE_1 =
  {\rm max}\big(\limsup\limits_{t \to +\infty}\frac{1}{t}\ln|g(t)|,
  \limsup\limits_{t \to +\infty}\frac{1}{t}\ln|g^{-1}(t)|\big),
  \LCE_2 = 0$;
 $
  \LEs_{1,2} = {\rm max, min}
  \big(
     \limsup\limits_{t \to +\infty}\frac{1}{t}\ln|g(t)|,
     \limsup\limits_{t \to +\infty}\frac{1}{t}\ln|g^{-1}(t)|
  \big).
 $
} that
\begin{equation}\label{LE<=LCE}
    \LEs_i(t,u) \leq \LCE_i(t,u), \ \forall t\geq 0,\ \forall u \in U \quad i=1,2,..,n.
\end{equation}

\begin{definition} \cite{ConstantinF-1985,ConstantinFT-1985}
The relative global (or uniform) Lyapunov exponents of singular value functions
of the dynamical system $\{\varphi^t\}_{t\geq0}$
with respect to the compact invariant set $K \subset U$ are defined as
\[
 \begin{aligned}
 &
 \LEo_{1}(K) &:=& \limsup_{t \to +\infty}\sup_{u \in K} \LEs_1(t,u),
 \\ &
 \LEo_{i+1}(K) &:=&
 \limsup_{t \to +\infty}\sup_{u \in K} \big(\LEs_1(t,u)+\cdots+\LEs_{i+1}(t,u)\big) -
 \\&&
 -&\limsup_{t \to +\infty}\sup_{u \in K} \big(\LEs_1(t,u)+\cdots+\LEs_{i}(t,u)\big), \ i=1,...,n-1.
 \end{aligned}
\]
\end{definition}
For $i=1,2,...,n$ we have
\[
 \LEo_{1}(K)+\cdots+\LEo_{i}(K) = \limsup_{t \to +\infty}\sup_{u \in K}\big(\LEs_1(t,u)+\cdots+\LEs_{i}(t,u)\big).
\]

By \eqref{1twinflim} and \eqref{wsupmax} we get (see, e.g. \cite[pp.360-361]{Temam-1997})
\[
 \begin{aligned}
 &
 \LEo_{1}(K) &=& \lim_{t \to +\infty}\max_{u \in K} \LEs_1(t,u),
 \\ &
 \LEo_{i+1}(K) &=&
 \lim_{t \to +\infty}\max_{u \in K} \big(\LEs_1(t,u)+\cdots+\LEs_{i+1}(t,u)\big) -
 \\&&
 -&\lim_{t \to +\infty}\max_{u \in K} \big(\LEs_1(t,u)+\cdots+\LEs_{i}(t,u)\big), \ i=1,...,n-1.
 \end{aligned}
\]

From \eqref{suplim-limsup} (see, e.g. \cite[p.49]{Eden-1989-PhD}, \cite[p.146]{EdenFT-1991})
for $u \in K$ we obtain the following inequality
\begin{equation}\label{LEvol<LEK}
   \LEo_{1}(u)+\cdots+\LEo_{i}(u)
   \leq \LEo_{1}(K)+\cdots+\LEo_{i}(K), \ i=1,2,...,n.
\end{equation}
At the same time, according to \eqref{suplim=limsup},
there exists $u^{cr}(m) \in K$ (it may be not unique)
such that the above expressions in \eqref{LEvol<LEK} coincide \cite{Eden-1989,Eden-1989-PhD,EdenFT-1991}:
\begin{equation}\label{gLEgequalLE}
\begin{aligned}
   & \LEo_{1}(K)+\cdots+\LEo_{m}(K) = \LEo_{1}(u^{cr}(m))+\cdots+\LEo_{m}(u^{cr}(m)) = \\
   & \qquad
   = \max_{u \in K} \big(\LEo_{1}(u)+\cdots+\LEo_{m}(u)\big).
\end{aligned}
\end{equation}

Various characteristics of chaotic behavior are based on Lyapunov exponents
(e.g., LEs are used in the Kaplan-Yorke formula of the Lyapunov dimension
and the sum of positive LEs may be used \cite{Millionschikov-1976,Pesin-1977}
as the characteristic of Kolmogorov-Sinai entropy rate \cite{Kolmogorov-1959,Sinai-1959}).
The properties of Lyapunov exponents and their various generalizations
are studied, e.g., in \cite{Lyapunov-1892,BylovVGN-1966,Oseledec-1968,Pesin-1977,Adrianova-1998,KunzeK-2001,KuznetsovL-2005,LeonovK-2007,BarreiraG-2011,Izobov-2012,CzornikNN-2013,KuznetsovAL-2014-arXiv-LE,LeonovAK-2015}.

The existence of different definitions of LEs, computational methods,
and related assumptions led to the appeal:
{\it ''Whatever you call your exponents, please state clearly how are they being computed''}
\cite{ChaosBook}.


\subsection{Kaplan-Yorke formula of the Lyapunov dimension}

\subsubsection{Kaplan-Yorke formula with respect to the finite time Lyapunov exponents}
Consider the dynamical system
$\big(\{\varphi^t\}_{t\geq0},(U\subseteq \mathbb{R}^n,||\cdot||) \big)$.
For $t>0$ we have
\begin{equation}\label{defomegat} 
  \frac{1}{t}\ln(\omega_d(D\varphi^t(u))) = \left\{
  \begin{aligned}
    & 0, && d=0, \\
    & \sum_{i=1}^{\lfloor d \rfloor}\LEs_i(t,u),   && d=\lfloor d \rfloor \in \{1,\ldots,n\}, \\ 
    & \sum_{i=1}^{\lfloor d \rfloor}\LEs_i(t,u) + ({d-\lfloor d \rfloor})\LEs_{\lfloor d \rfloor+1}(t,u), && d \in (0,n). \\
  \end{aligned}
  \right.
\end{equation}
If $\ln(\omega_{n}(D\varphi^t(u)))\leq0$
for fixed $t>0$ and point $u \in K$,
then by \eqref{maplocalDOinf} for $d(t,u):=d_{\rm L}(\varphi^t,u)$
we have
\begin{equation}\label{constr}
  \begin{aligned}
  \frac{1}{t}\ln(\omega_{d(t,u)}(D\varphi^t(u))) = 0,
  \quad
  \frac{1}{t}\ln(\omega_{d(t,u)+\delta}(D\varphi^t(u))) < 0
  \quad \forall \delta \in (0,n-d(t,u)].
  \end{aligned}
\end{equation}
Let for $t>0$
\[
  \begin{aligned}
  &
  j(t,u)=j\big(\{\LEs_i(t,u)\}_1^n\big) := \lfloor d(t,u) \rfloor \in \{0,..,n\},
  \\
  & s(t,u)=s\big(\{\LEs_i(t,u)\}_1^n\big) := d(t,u) - j(t,u) \in [0,1).
  \end{aligned}
\]

Then for $j(t,u)\leq n-1$
from \eqref{constr} it follows that
\begin{equation}\label{jut}
  \sum\limits_{i=1}^{j(t,u)}\LEs_i(t,u) \geq 0, \quad
  \sum\limits_{i=1}^{j(t,u)+1}\LEs_i(t,u)<0.
\end{equation}
We have
\[
\begin{aligned}
  &
  \frac{1}{t}\ln(\omega_{d(t,u)}(D\varphi^t(u))) =
  \frac{1}{t}\ln \big(
    (\omega_{j(t,u)}(D\varphi^t(u)))^{(1-s(t,u))}
    (\omega_{j(t,u)+1}(D\varphi^t(u)))^{s(t,u)}
  \big)=
  \\ &
  =
  (1-s(t,u))\sum\limits_{i=1}^{j(t,u)}\LEs_k(t,u)+s(t,u)\sum\limits_{i=1}^{j(t,u)+1}\LEs_{k}(t,u) =0.
\end{aligned}
\]
Therefore
\[
  s(t,u) = \left\{
  \begin{aligned}
     &\dfrac{\LEs_1(t,u)+\cdots+\LEs_{j(t,u)}(t,u)}{|\LEs_{j(t,u)+1}(t,u)|} < 1, & j(t,u) \in \{1,\ldots,n-1\}, \\
     & 0, & j(t,u)=0 \mbox{\ or\ } j(t,u)=n.
  \end{aligned}
  \right.
\]

The expression
\begin{equation}\label{lftKY}
  d_{\rm L}^{\rm KY}(\{\LEs_i(t,u)\}_{1}^{n}) :=j(t,u)+\dfrac{\LEs_1(t,u)+\cdots+\LEs_{j(t,u)}(t,u)}{|\LEs_{j(t,u)+1}(t,u)|}
\end{equation}
corresponds to the Kaplan-Yorke formula \cite{KaplanY-1979}
with respect to the finite time Lyapunov exponents, i.e.
the ordered set $\{\LEs_i(t,u)\}_{1}^{n}$.
The idea of $d_{\rm L}^{\rm KY}$ construction
may be used with other types of Lyapunov exponents (see below).

Further we assume that the relation $s(t,u)=0$
for $j(t,u)=0$ and $j(t,u)=n$ follows from the first expression for $s(t,u)$.
Since
$\frac{1}{t}\ln(\omega_{d(t,u)}(D\varphi^t(u))) \leq 0 \Leftrightarrow
\omega_{d(t,u)}(D\varphi^t(u)) \leq 1$ for $t>0$,
from \eqref{DOmapt} we have
\begin{proposition}\label{thm:ftKY}
\begin{equation}\label{ftKY}
 d_{\rm L}(\varphi^t,K) =
 \sup_{u\in K}d_{\rm L}(\varphi^t,u) =
 \sup_{u\in K}d_{\rm L}^{\rm KY}(\{\LEs_i(t,u)\}_{1}^{n})=
 \sup_{u\in K}
 \left(
   j(t,u)+\dfrac{\LEs_1(t,u)+\cdots+\LEs_{j(t,u)}(t,u)}{|\LEs_{j(t,u)+1}(t,u)|}
 \right).
\end{equation}
\end{proposition}

While in computing we can consider only finite time $t\leq T$,
from a theoretical point of view, it may be interesting to study
the limit behavior of $\sup_{u\in K}d_{\rm L}^{\rm KY}(\{\LEs_i(t,u)\}_{1}^{n})$
as $t \to +\infty$.

\subsubsection{Kaplan-Yorke formula with respect to
the relative global Lyapunov exponents of singular value functions}
Let
\[
  \begin{aligned}
  & j = j\big(\{\LEo_i(K)\}_1^n\big)
  := \max\{m \in \{0,\ldots,n\}:\ \sum\limits_{i=1}^{m}\LEo_i(K) \geq 0 \}, \quad
  \\ &
  s = s\big(\{\LEo_i(K)\}_1^n\big) := \left\{
  \begin{aligned}
     & 0 \leq \frac{\LEo_1(K)+\cdots+\LEo_{j}(K)}{|\LEo_{j+1}(K)|} < 1, & j \in \{1,\ldots,n-1\}, \\
     & 0, & j = 0 \mbox{\ or\ } j=n.
  \end{aligned}
  \right.
  \end{aligned}
\]
The expression
\(
  d_{\rm L}^{\rm KY}(\{\LEo_i(K)\}_1^n):=
  j+\frac{\LEo_1(K)+\cdots+\LEo_{j}(K)}{|\LEo_{j+1}(K)|}
\)
is the Kaplan-Yorke formula of Lyapunov dimension
with respect to the relative global Lyapunov exponents of singular value functions.
Then
\[
 \begin{aligned}
  &
  \lim_{t \to +\infty} \max_{u \in K}
  \frac{1}{t}\ln\big(\omega_{j+s}(D\varphi^t(u))\big)
  =
  \lim_{t \to +\infty} \max_{u \in K}
  \left(
    \sum\limits_{i=1}^{j}\LEs_i(t,u)+s\LEs_{j+1}(t,u)
  \right)
  =
  \\ &
  =
  \lim_{t \to +\infty} \max_{u \in K}
  \left(
    (1-s)\sum\limits_{i=1}^{j}\LEs_i(t,u)+s\sum\limits_{i=1}^{j+1}\LEs_i(t,u)
  \right) \leq
  \\ &
  (\mbox{since, in general, the maximums may be achieved at different points $u$})
  \\ &
  \leq
  \lim_{t \to +\infty} \max_{u \in K} (1-s)\sum\limits_{i=1}^{j}\LEs_i(t,u)
  +
  \lim_{t \to +\infty}\max_{u \in K}  s\sum\limits_{i=1}^{j+1}\LEs_i(t,u)
  = \\ &
  = (1-s)\lim\limits_{t \to +\infty} \max_{u \in K} \sum\limits_{i=1}^{j}\LEs_i(t,u)
  +
  s \lim\limits_{t \to +\infty} \max_{u \in K} \sum\limits_{i=1}^{j+1} \LEs_i(t,u) =0.
  \end{aligned}
\]
Thus, for any $\overline{s}:\ s<\overline{s}<1$,
$\lim\limits_{t \to +\infty} \max\limits_{u \in K}\frac{1}{t}\ln(\omega_{j+\overline{s}}(D\varphi^t(u)))<0$
and from Definition~\ref{defdimLt} we have
\begin{proposition}\label{thm:dimLdLKYglob} (see, e.g. \cite[pp.30-31]{ConstantinFT-1985})
  \[
    d_{\rm L}^{\rm E}(\{\varphi^t\}_{t\geq0},K) \leq
    d_{\rm L}^{\rm KY}(\{\LEo_i(K)\}_1^n)
  \]
\end{proposition}

Under some conditions we can obtain the equality.
\alarm{}
\begin{corollary}\label{thm:uo-sumk-k1}
If critical points $u^{cr}(j)$ and $u^{cr}(j+1)$ from \eqref{suplim=limsup}
coincide, i.e. $u^{cr}=u^{cr}(j)=u^{cr}(j+1)$,
then
\begin{equation}\label{uo-sumk-k1}
\begin{aligned}
  &
  \lim\limits_{t \to +\infty}\sum\limits_{k=1}^{j}\LEs_k(t,u^{cr}) =
  \lim\limits_{t \to +\infty}\max_{u \in K}\sum\limits_{k=1}^{j}\LEs_k(t,u),
  \\ &
  \lim\limits_{t \to +\infty}\sum\limits_{k=1}^{j+1}\LEs_k(t,u^{cr}) =
  \lim\limits_{t \to +\infty}\max_{u \in K}\sum\limits_{k=1}^{j+1}\LEs_k(t,u),
\end{aligned}
\end{equation}
and
\[
  d_{\rm L}(\{\varphi^t\}_{t\geq0},K) = d_{\rm L}^{\rm KY}(\{\LEo_i(K)\}_1^n).
\]
\end{corollary}
In \cite[p.565]{Gelfert-2003} the systems, having property \eqref{uo-sumk-k1},
are called ``typical systems''.

\subsubsection{Kaplan-Yorke formula with respect to
relative Lyapunov exponents of singular value functions}

Let
\[
  j(u) = j\big(\{\LEo_i(u)\}_1^n\big)
  := \max\{m \in \{0,\ldots,n\}: \sum\limits_{i=1}^{m}\LEo_i(u) \geq 0 \}
\]
\[
  s(u) =s\big(\{\LEo_i(u)\}_1^n\big) := \left\{
  \begin{aligned}
     & 0 \leq \frac{\LEo_1(u)+\cdots+\LEo_{j(u)}(u) }{|\LEo_{j(u)+1}(u)|} < 1, & j(u) \in \{1,\ldots,n-1\}, \\
     & 0, &   j(u)=0 \mbox{\ or\ } j(u)=n.
  \end{aligned}
  \right.
\]

The expression
\(
  d_{\rm L}^{\rm KY}(\{\LEo_i(u)\}_1^n):=
  j(u)+\frac{\LEo_1(u)+\cdots+\LEo_{j(u)}(u)}{|\LEo_{j(u)+1}(u)|}
\)
is the Kaplan-Yorke formula
of Lyapunov dimension
with respect to the relative Lyapunov exponents of singular value functions.
We have
\begin{equation}\label{itlnwdest}
 \begin{aligned}
  &
  \limsup\limits_{t \to +\infty}
  \frac{1}{t}\ln\big(\omega_{j(u)+s(u)}(D\varphi^t(u))\big)
  =
  \limsup_{t \to +\infty}
  \left(
    \sum\limits_{i=1}^{j(u)}\LEs_i(t,u)+s(u)\LEs_{j(u)+1}(t,u)
  \right)
  =\\ &
  =
  \limsup_{t \to +\infty}
  \left(
    (1-s(u))\sum\limits_{i=1}^{j(u)}\LEs_i(t,u)+s(u)\sum\limits_{i=1}^{j(u)+1}\LEs_i(t,u)
  \right) \leq
  \\ &
  \leq
  \limsup_{t \to +\infty} (1-s(u))\sum\limits_{i=1}^{j(u)}\LEs_i(t,u)
  +
  \limsup_{t \to +\infty} s(u)\sum\limits_{i=1}^{j(u)+1}\LEs_i(t,u)
  = \\ &
  =
  (1-s(u))\limsup\limits_{t \to +\infty} \sum\limits_{i=1}^{j(u)}\LEs_i(t,u)
  +
  s(u) \limsup\limits_{t \to +\infty} \sum\limits_{i=1}^{j(u)+1} \LEs_i(t,u)=0.
  \end{aligned}
\end{equation}
Thus, for any $j(u)<n$ and $\overline{s}:\ s(u)<\overline{s}<1$,
$\limsup\limits_{t \to +\infty}\frac{1}{t}\ln\big(\omega_{j(u)+\overline{s}}(D\varphi^t(u))\big)<0$
and from Definition~\ref{deflocalDL} we have
\begin{proposition}\label{thm:dimLdLKYloc}
\[
  \sup_{u\in K} d_{\rm L}(\{\varphi^t\}_{t\geq0},u)
  \leq \sup_{u\in K} d_{\rm L}^{\rm KY}(\{\LEo_i(u)\}_1^n)
  = \sup_{u\in K}
  \left(
     j(u)+\frac{\LEo_1(u)+\cdots+\LEo_{j(u)}(u)}{|\LEo_{j(u)+1}(u)|}
  \right).
\]
\end{proposition}

\begin{proposition} (see, e.g. \cite[p.60]{Eden-1989-PhD})
 \begin{equation}\label{DKYloc<DKYglob}
  \sup_{u\in K} (d_{\rm L}^{\rm KY}(\{\LEo_i(u)\}_1^n))
  \leq d_{\rm L}^{\rm KY}(\{\LEo_i(K)\}_1^n).
\end{equation}
\end{proposition}
{\bf Proof.}
The assertion follows from the relation (see \cite[p.60]{Eden-1989-PhD})
\begin{equation}\label{jLEcoinside}
  \sup_{u \in K} j\big(\{\LEo_i(u)\}_1^n\big)=
  j\big(\{\LEo_i(K)\}_1^n\big)
\end{equation}
and inequality \eqref{LEvol<LEK}.
$\blacksquare$

Remark that there are examples in which inequality \eqref{DKYloc<DKYglob} is strict
(see, e.g. \cite[pp.49-51,62-63]{Eden-1989-PhD}\footnote{
 Let $\LEs_1(t,u) = (e^u)^t$, $\LEs_2(t,u) = (\frac{1}{2}(1-u))^t$ for all $u \in K=[0,1]$.
 Thus $\LEs_1(u)=\LEo_1(u) = u$, $\LEs(u)=\LEo_2 = \ln(1-u)-\ln2$;
 $\LEo_1(K) = 1$, $\LEo_2 =-1-\ln2$;
 Here $u^{cr}(1)=1$:\ $\LEo_1(1)=\LEo_1(K)=1$;
 $u^{cr}(2)=0$:\ $\LEo_1(0)+\LEo_2(0)=\LEo_1(K)+\LEo_2(K)=-\ln2$.
 Then $\sup_{u \in [0,1]}d_{\rm L}^{\rm KY}(\{\LEo_i(u)\}_1^2)=\frac{u}{\ln2-\ln(1-u)} < 1+\frac{1}{1+\ln2}=
 d_{\rm L}^{\rm KY}(\{\LEo_i(K)\}_1^2)$.
}).

\subsubsection{Kaplan-Yorke formula with respect to
the Lyapunov exponents of singular values}

Let (see, e.g. \cite[pp.32-34]{ConstantinFT-1985})
\[
  \begin{aligned}
  & j(u) = j\big(\{\LEs_i(u)\}_1^n\big)
   := \max\{m \in \{0,\ldots,n\}: \sum\limits_{i=1}^{m}\LEs_i(u) \geq 0 \},
  \\ &
  s(u)= s\big(\{\LEs_i(u)\}_1^n\big) := \left\{
  \begin{aligned}
     & 0 \leq \frac{\LEs_1(u)+\cdots+\LEs_{j(u)}(u) }{|\LEs_{j(u)+1}(u)|} < 1, &  j(u) \in \{1,\ldots,n-1\}, \\
     & 0, & \!\!\!\!\!\!\!\!\!\!\!\!\!\!\!\!\!\!\!\!\!\!\!\!\!\!\!\!  j(u)=0 \mbox{\ or\ } j(u)=n.
  \end{aligned}
  \right.
  \end{aligned}
\]
The expression
\(
  d_{\rm L}^{\rm KY}(\{\LEs_i(u)\}_1^n):=
  j(u)+\frac{\LEs_1(u)+\cdots+\LEs_{j(u)}(u) }{|\LEs_{j(u)+1}(u)|}
\)
is the Kaplan-Yorke formula
of Lyapunov dimension
with respect to the Lyapunov exponents of singular values.

Then
\[
 \begin{aligned}
  &
  \limsup\limits_{t \to +\infty} \frac{1}{t}\ln\big(\omega_{j(u)+s(u)}(D\varphi^t(u))\big)
  \leq
  (1-s(u))\limsup\limits_{t \to +\infty} \sum\limits_{i=1}^{j(u)}\LEs_i(t,u)
  +
  s(u) \limsup\limits_{t \to +\infty} \sum\limits_{i=1}^{j(u)+1} \LEs_i(t,u)
  \leq
  \\ &
  \leq (1-s(u))\sum\limits_{i=1}^{j(u)}\LEs_i(u)
  +
  s(u) \sum\limits_{i=1}^{j(u)+1} \LEs_i(u)=0.
  \end{aligned}
\]
For $j(u)<n$ and any $\overline{s}:\ s(u)<\overline{s}<1$,
$\limsup\limits_{t \to +\infty}\frac{1}{t}\ln\big(\omega_{j(u)+\overline{s}}(D\varphi^t(u))\big)<0$
and from Definition~\ref{deflocalDL} we have

\begin{proposition}
  \[
    \sup_{u\in K} d_{\rm L}(\{\varphi^t\}_{t\geq0},u)
    \leq \sup_{u\in K} d_{\rm L}^{\rm KY}(\{\LEs_i(u)\}_1^n)
    =
    \sup_{u\in K}
    \left(
    j(u)+\frac{\LEs_1(u)+\cdots+\LEs_{j(u)}(u) }{|\LEs_{j(u)+1}(u)|}
    \right).
  \]
\end{proposition}

\subsubsection{Kaplan-Yorke formula with respect to
the Lyapunov exponents of fundamental matrix columns}
Let
\[
  \begin{aligned}
  & j(u) = j\big(\{\LCE_i(u)\}_1^n\big)
  := \max\{m \in \{0,\ldots,n\}: \sum\limits_{k=1}^{m}\LCE_k(u) \geq 0 \},
  \\ &
  s(u) = s \big(\{\LCE_i(u)\}_1^n\big)
  := \left\{
  \begin{aligned}
     & 0 \leq \frac{\LCE_1(u)+\cdots+\LCE_{j(u)}(u) }{|\LCE_{j(u)+1}(u)|} < 1, &  j(u) \in \{1,\ldots,n-1\}, \\
     & 0, & \!\!\!\!\!\!\!\!\!\!\!\!\!\!\!\!\!\!\!\!\!\!\!\!\!\!\!\!  j(u)=0 \mbox{\ or\ }   j(u)=n.
     \end{aligned}
  \right.
  \end{aligned}
\]
The expression
\(
  d_{\rm L}^{\rm KY}(\{\LCE_i(u)\}_1^n):=
  j(u)+\frac{\LCE_1(u)+\cdots+\LCE_{j(u)}(u) }{|\LCE_{j(u)+1}(u)|}
\)
is the Kaplan-Yorke formula
of Lyapunov dimension
with respect to the Lyapunov exponents of fundamental matrix columns.

Then, similar to \eqref{itlnwdest}, by \eqref{LE<=LCE} we obtain
\[
 \begin{aligned}
  &
  \limsup\limits_{t \to +\infty} \frac{1}{t}\ln\big(\omega_{j(u)+s(u)}(D\varphi^t(u))\big)
  \leq
  & \\
 & \leq
  (1-s(u))\sum\limits_{i=1}^{j(u)}\LCE_i(u)
  +
  s(u) \sum\limits_{i=1}^{j(u)+1} \LCE_i(u)=0.
  \end{aligned}
\]

Thus, for $j(u)<n$ and any $\overline{s}:\ s(u)<\overline{s}<1$,
$\limsup\limits_{t \to +\infty}\frac{1}{t}\ln\big(\omega_{j(u)+\overline{s}}(D\varphi^t(u))\big)<0$
and from Definition~\ref{deflocalDL} we get
\begin{proposition}
 \[
    \sup_{u\in K} d_{\rm L}(\{\varphi^t\}_{t\geq0},u)
    \leq \sup_{u\in K} d_{\rm L}^{\rm KY}(\{\LCE_i(u)\}_1^n)
    =
    \sup_{u\in K}
    \left(
       j(u)+\frac{\LCE_1(u)+\cdots+\LCE_{j(u)}(u) }{|\LCE_{j(u)+1}(u)|}
    \right).
 \]
 \end{proposition}

\subsubsection{Computation by the Kaplan-Yorke formulas}

  For a given invariant set $K$ and a given point $u_0 \in K$
  there are two essential questions related
  to the computation of Lyapunov exponents and
  the use of the Kaplan-Yorke formulas of local Lyapunov dimension
  $\sup_{u \in K}d_{\rm L}^{\rm KY}(\{\LEs_i(u)\}_1^n)$ and
  $\sup_{u \in K}d_{\rm L}^{\rm KY}(\{\LEo_i(u)\}_1^n)$:
  \medskip

  (a)
  $\limsup\limits_{t\to+\infty}\LEs_i(t,u_0) \,?\!=\, \lim\limits_{t\to+\infty}\LEs_i(t,u_0)$
  \ or \
  $\limsup\limits_{t \to +\infty}(\sum\limits_1^i\LEs_i(t,u))
  \,?\!=\,
  \lim\limits_{t \to +\infty}(\sum\limits_1^m\LEs_i(t,u))$

  (b) if the above limits do not exist, then

  $
   \sup_{u \in K} d_{\rm L}^{\rm KY}(\{\LEs_m(u)\}_1^n)
   \,?\!=\,
   \sup_{u \in K \backslash \{\varphi^t(u_0), t \geq 0 \}} d_{\rm L}^{\rm KY}(\{\LEs_i(u)\}_1^n)
   $

   or

  $
   \sup_{u \in K} d_{\rm L}^{\rm KY}(\{\LEo_i(u)\}_1^n)
   \,?\!=\,
   \sup_{u \in K \backslash \{\varphi^t(u_0), t \geq 0 \}} d_{\rm L}^{\rm KY}(\{\LEo_i(u)\}_1^n)
   $.

  \medskip

  In order to get rigorously the positive answer to these questions,
  from a theoretical point of view, one may use
  various ergodic properties of the dynamical system $\{\varphi^t\}_{t\geq0}$
  (see, Oseledec~\cite{Oseledec-1968}, Ledrappier~\cite{Ledrappier-1981},
  and some auxiliary results in \cite{BogoliubovK-1937,DellnitzJ-2002}).
  However, from a practical point of view,
  the rigorous use of the above results is a challenging task
  (e.g. even for the well-studied Lorenz system)
  and hardly can be done effectively in the general case
  (see, e.g. the corresponding discussions
  in \cite{BarreiraS-2000},\cite[p.118]{ChaosBook},\cite{OttY-2008},\cite[p.9]{Young-2013}
   and the works on the Perron effects of the largest Lyapunov exponent sign reversals \cite{LeonovK-2007,KuznetsovL-2005}).
  For an example of the effective rigorous use of the ergodic theory
  for the estimation of the Hausdorff and Lyapunov dimensions see, e.g. \cite{Schmeling-1998}.

  Thus, in the general case, from a practical point of view,
  one cannot rely on the above relations (a) and (b)
  and shall use $\limsup_{t\to+\infty}$ in the definitions
  of local Lyapunov exponents
  and the corresponding formulas for the Lyapunov dimension
  (see, e.g. Temam~\cite{Temam-1997}).


  \medskip

  However,
  if $u_0$ is an equilibrium point, 
  then 
  the expression "$\limsup_{t \to +\infty}$"
  in Definitions~\ref{def:lLE},~\ref{def:LCE}, and \ref{def:lvolLE}
  can be replaced by "$\lim_{t \to +\infty}$"
  and we have
\begin{lemma}\label{thm:allLE}
  Let $\varphi^t(u_0)$ be a stationary point, i.e. $\varphi^t(u_0)\equiv u_0$. 
  Then for $i=1,2,...,n$ we have
  \[
    \lim\limits_{t \to +\infty} \LEs_i(t,u_0) =\LEs_i(u_0) =\LEo_i(u_0)=\LCE_i(u_0).
  \]
\end{lemma}
Thus, for $j = j\big(\{\LEo_i(K)\}_1^n\big)$, we get
\begin{proposition}\label{thm:dimu0}
If critical points in \eqref{ucrL} and \eqref{uo-sumk-k1} coincide with
an equilibrium point $u^{cr}_{eq}$,
i.e  $\varphi^t(u^{cr}_{eq})\equiv u^{cr}_{eq}=u^{cr}_{\rm L}\equiv u^{cr}(j)=u^{cr}(j+1)$, then
\[
  d_{\rm L}(\{\varphi^t\}_{t\geq0},u^{cr}_{eq})
  = d_{\rm L}(\{\varphi^t\}_{t\geq0},K) = d_{\rm L}^{\rm KY}(\{\LEo_i(K)\}_1^n)
\]
and
\[
\begin{aligned}
    &
    d_{\rm L}(\{\varphi^t\}_{t\geq0},u^{cr}_{eq})
    =
    \sup_{u\in K} d_{\rm L}^{\rm KY}(\{\LEo_i(u)\}_1^n)
    =
    \sup_{u\in K} d_{\rm L}^{\rm KY}(\{\LEs_i(u)\}_1^n)
    =
    \sup_{u\in K} d_{\rm L}^{\rm KY}(\{\LCE_i(u)\}_1^n)=
    \\ &
    =
    \lim\limits_{t \to +\infty}\max_{u \in K}d_{\rm L}^{\rm KY}(\{\LEs_i(t,u)\}_{1}^{n})
\end{aligned}
\]
\end{proposition}
If $u^{cr}_{\rm L}= u^{cr}(j)=u^{cr}(j+1)$ belongs to a periodic orbit with period $T$,
then the same reasoning can be applied for $(\varphi^{T})^t$.

The last section of this survey is devoted to the examples in which
the maximum of the local Lyapunov dimension achieves at an equilibrium point.

Taking into account the existence of different definitions of Lyapunov dimension and related formulas
and following \cite{ChaosBook},
we recommend that {\it whatever you call your Lyapunov dimension, please state clearly how
is it being computed}.

\section{Analytical estimates of the 
Lyapunov dimension and its invariance with respect to diffeomorphisms}

Along with widely used numerical methods for estimating and computing
the Lyapunov dimension
(see, e.g. MATLAB realizations of the methods based on QR and SVD decompositions
in \cite{KuznetsovMV-2014-CNSNS,LeonovKM-2015-EPJST})
there is an effective analytical approach, proposed by G.A.Leonov in 1991 \cite{Leonov-1991-Vest}
(see also \cite{LeonovB-1992,Leonov-2002,BoichenkoLR-2005,Leonov-2008,Leonov-2012-PMM,LeonovK-2015-AMC,LeonovKM-2015-EPJST}).
The Leonov method is based on the direct Lyapunov method
with special Lyapunov-like functions.
The advantage of this method is that it allows one to estimate the Lyapunov dimension
of invariant set without localization of the set in the phase space
and in many cases get
effectively
exact Lyapunov dimension formula
\cite{Leonov-2002,LeonovP-2005,LeonovPS-2013-PLA,LeonovK-2015-AMC,LeonovKKK-2015-arXiv-YangTigan,LeonovAK-2015,LeonovKKK-2015-arXiv-Lorenz,LeonovKM-2015-ArXiv}.

Following \cite{Kuznetsov-2016}, next the invariance of Lyapunov dimension
with respect to diffeomorphisms and its relation with the Leonov method are discussed.
An analog of Leonov method for discrete time dynamical systems is suggested.

While topological dimensions are invariant with respect to Lipschitz homeomorphisms,
the Hausdorff dimension is invariant with respect to Lipschitz diffeomorphisms
and noninteger Hausdorff dimension is not invariant with respect to
homeomor\-phisms \cite{HurewiczW-1941}.
Since the Lyapunov dimension is used as an upper estimate of Hausdorff dimension,
the question arises whether the Lyapunov dimension
is invariant under diffeomorphisms
(see, e.g. \cite{OttWY-1984,Kuznetsov-2016-ArXiv}).

Consider the dynamical system
$\big(\{\varphi^t\}_{t\geq0},(U\subseteq \mathbb{R}^n,||\cdot||) \big)$
under the change of coordinates $w = h(u)$,
where $h: U \subseteq \mathbb{R}^n \to \mathbb{R}^n$ is a diffeomorphism.
In this case
the semi-orbit
$\gamma^{+}(u) = \{\varphi^t(u), t \geq 0 \}$
is mapped to the semi-orbit defined by
$\varphi_h^t(w)=\varphi_h^t(h(u))=h(\varphi^t(u))$,
the dynamical system
$\big(\{\varphi^t\}_{t\geq0},(U\subseteq \mathbb{R}^n,||\cdot||) \big)$
is transformed to
the dynamical system
$\big(\{\varphi_h^t\}_{t\geq0},(h(U)\subseteq \mathbb{R}^n,||\cdot||) \big)$,
and a compact set $K \subset U$ invariant with respect to $\{\varphi^t\}_{t\geq0}$
is mapped to the compact set $h(K) \subset h(U)$
invariant with respect to $\{\varphi_h^t\}_{t\geq0}$.
Here
\[
  D_w\varphi_h^t(w)=
  D_w\big(h(\varphi^t(h^{-1}(w)))\big)=
  D_uh(\varphi^t(h^{-1}(w)))
  D_u\varphi^t(h^{-1}(w))
  D_w h^{-1}(w),
\]
\[
  D_u \big( \varphi_h^t(h(u)) \big)=
  D_w\varphi_h^t(h(u)) D_uh(u)=
  D_u\big(h(\varphi^t(u))\big)=
  D_u h(\varphi^t(u))D_u\varphi^t(u).
\]
Therefore
\[
  D_w h^{-1}(w) = \big(D_uh(u)\big)^{-1}
\]
and\begin{equation}\label{Dphih}
  D\varphi_h^t(w)=
  Dh(\varphi^t(u))
  D\varphi^t(u)
  \big(Dh(u)\big)^{-1}.
\end{equation}

If $u \in K$, then $\varphi^t(u)$ and $\varphi_h^t(h(u))$
define bounded semi-orbits.
Remark that $Dh$ and $(Dh)^{-1}$ are continuous and, thus,
$Dh(\varphi^t(u))$ and $(Dh(\varphi^t(u)))^{-1}$ are bounded in $t$.
From \eqref{wsupmax} it follows that
for any $d \in [0,n]$ there is a constant $c=c(d)\geq 1$
such that for any $t\geq 0$
\begin{equation}\label{Dhest}
\begin{aligned}
  &
  \max_{u \in K}\omega_d\big(Dh(u)\big)\leq c,
  \quad
  \max_{u \in K}\omega_d\big((Dh(u))^{-1}\big)
  \leq c, \quad t\geq 0.
\end{aligned}
\end{equation}

\begin{lemma}\label{thm:hdiff}
If for a fixed $t>0$  there exist diffeomorphism $h: U \subseteq \mathbb{R}^n \to \mathbb{R}^n$
and $d \in [0,n]$ such that the estimation
\begin{equation}\label{wDphiht<1}
  \max_{w \in h(K)}\omega_d\big(D\varphi_h^t(w)\big)
  =
  \max_{u \in K}
  \omega_d\bigg(
  Dh(\varphi^t(u))
  D\varphi^t(u)
  \big(Dh(u)\big)^{-1}
  \bigg) <1
\end{equation}
is valid\footnote{The expression in \eqref{wDphiht<1}
corresponds to the expressions considered in
\cite[eq.(1)]{Leonov-1991-Vest} for $p(u)=Dh(u)$,
\cite[eq.(1)]{Leonov-2002} and \cite[p.99, eq.10.1]{Leonov-2008} for $Q(u)=Dh(u)$.},
then for $u \in K$
 \[
    \liminf\limits_{t \to +\infty}
      \bigg(
      \omega_d\big(D\varphi^t(u)\big)
      -
      \omega_d\big(D\varphi_h^t(h(u))\big)
      \bigg)=0
 \]
 and
 \[
  \liminf\limits_{t \to +\infty}
  \omega_d\big( D\varphi_h^t(h(u)) \big)
  =
  \liminf\limits_{t \to +\infty}
  \omega_d\big(D\varphi^t(u)\big)=0.
 \]
\end{lemma}
{\bf Proof.}
Applying \eqref{horn} to \eqref{Dphih}, we get
\[
  \omega_d\big(D\varphi_h^t(h(u))\big) \leq
  \omega_d\big(Dh(\varphi^t(u))\big)
  \omega_d\big(D\varphi^t(u)\big)
  \omega_d\big(\big(Dh(u)\big)^{-1}\big).
\]
By \eqref{Dhest} we obtain
\[
  \omega_d\big(D\varphi_h^t(h(u))\big) \leq
  c^2
  \omega_d\big(D\varphi^t(u)\big).
\]
Similarly
\[
  \omega_d\big(D\varphi^t(u)\big) \leq
  \omega_d\big(\big(Dh(\varphi^t(u))\big)^{-1}\big)
  \omega_d\big(D\varphi_h^t(h(u))\big)
  \omega_d\big(Dh(u)\big)
\]
and
\[
  \omega_d\big(D\varphi^t(u)\big) \leq
  c^2
  \omega_d\big(D\varphi_h^t(h(u))\big).
\]
Therefore for any $d \in [0,n]$, $t\geq0$, and $u \in K$
\begin{equation}\label{Omega2Horn}
  c^{-2}
  \omega_d\big(D\varphi_h^t(h(u))\big)
  \leq
  \omega_d\big(D\varphi^t(u)\big)
  \leq
  c^2
  \omega_d\big(D\varphi_h^t(h(u))\big)
\end{equation}
and
\[
  (c^{-2}-1)
  \omega_d\big(D\varphi_h^t(h(u))\big)
  \leq
  \omega_d\big(D\varphi^t(u)\big)
  -
  \omega_d\big(D\varphi_h^t(h(u))\big)
  \leq
  (c^2-1)
  \omega_d\big(D\varphi_h^t(h(u))\big).
\]
If for a fixed $t\geq0$ there is $d \in [0,n]$ such that
$\sup_{u \in K}\omega_d\big(D\varphi_h^t(h(u))\big) <1$,
then by \eqref{liminf} we have
\[
  \liminf\limits_{t \to +\infty}\omega_d\big(D\varphi_h^t(h(u))\big)=0
\]
and
\[
  0 \leq
  \liminf\limits_{t \to +\infty}
  \left(
  \omega_d\big(D\varphi^t(u)\big)
  -
  \omega_d\big(D\varphi_h^t(h(u))\big)
  \right)
  \leq
  0.
\]
$\blacksquare$
\begin{corollary}\label{thm:hdiffLE} (see, e.g. \cite{KuznetsovAL-2016})
For $u \in K$ we have
\[
  \lim\limits_{t \to +\infty}
  \bigg( \LEs_i\big( D\varphi_h^t(h(u)) \big) - \LEs_i\big(D\varphi^t(u)\big)\bigg)=0,
  \quad \quad i=1,2,..,n
 \]
 and, therefore,
 \[
  \limsup\limits_{t \to +\infty} \LEs_i\big( D\varphi_h^t(h(u)) \big)
  =
  \limsup\limits_{t \to +\infty} \LEs_i\big( D\varphi^t(u) \big),
  \quad \quad i=1,2,..,n.
 \]
\end{corollary}
{\bf Proof.}
For $t>0$ from \eqref{Omega2Horn} we get
\begin{equation}\label{LE2Horn}
  \frac{1}{t}\ln c^{-2} +
  \frac{1}{t}\ln\omega_d\big(D\varphi_h^t(h(u))\big)
  \leq
  \frac{1}{t}\ln\omega_d\big(D\varphi^t(u)\big)
  \leq
  \frac{1}{t}\ln c^2 +
  \frac{1}{t}\ln\omega_d\big(D\varphi_h^t(h(u))\big).
\end{equation}
Thus for the integer $d=m$ we have
\[
  \lim\limits_{t \to +\infty}
  \left(
  \frac{1}{t}\ln\omega_m\big(D\varphi^t(u)\big)
  -
  \frac{1}{t}\ln\omega_m\big(D\varphi_h^t(h(u))\big)
  \right)
  =
  \lim\limits_{t \to +\infty}
  \left(
  \sum_{i=1}^{m} \LEs_i\big( D\varphi^t(u) \big)
  -
  \sum_{i=1}^{m} \LEs_i\big( D\varphi_h^t(h(u)) \big)
  \right)
  = 0.
\]
$\blacksquare$

The above statements are rigorous reformulation from \cite{KuznetsovAL-2016,LeonovAK-2015}
and implies the following

\begin{proposition}\label{thm:dDOunderdiff}
The Lyapunov dimension of the dynamical system $\{\varphi^t\}_{t\geq0}$
with respect to the compact invariant set $K$
is invariant with respect to any diffeomorphism
$h: U \subseteq \mathbb{R}^n \to \mathbb{R}^n$, i.e.
\begin{equation}\label{dDOunderdiff}
\begin{aligned}
  &
  d_{\rm L}(\{\varphi^t\}_{t\geq0},K)
  =
  d_{\rm L}(\{\varphi_h^t\}_{t\geq0},h(K)).
\end{aligned}
\end{equation}
\end{proposition}
{\bf Proof.}
Lemma~\ref{thm:hdiff} implies that
if
$
  \max_{w \in h(K)}\omega_d\big(D\varphi_h^{t}(w)\big) < 1
$
for a fixed $t>0$ and $d \in [0,n]$,
then there exists $T>t$ such that
\begin{equation}\label{defT}
  \max_{u \in K}\omega_d\big(D\varphi^{T}(u)\big) < 1
\end{equation}
and vice verse.
Thus the set of $d$,
over which $\inf_{t>0}$ is taken in \eqref{DOmapt},
is the same for $D\varphi^t(u)$ and $D\varphi_h^t(w)$ and, therefore,
\[
  \inf_{t>0}\inf\{d\in[0,n]: \max\limits_{u\in K}\omega_{d}(D\varphi^t(u))<1\}
  =
  \inf_{t>0}\inf\{d\in[0,n]: \max\limits_{w\in h(K)}\omega_{d}(D\varphi_h^t(w))<1\}.
\]
$\blacksquare$

\begin{corollary}\label{thm:hdiffLE}
Suppose $H(u)$ is a $\,n\times n $ matrix, the elements of which are scalar continuous functions of $u$ and $\det H(u) \neq 0$ for $u \in K$.
If for a fixed $t>0$  there is $d \in (0,n]$ such that
\begin{equation}\label{wDphiH<1}
  \max_{w \in h(K)}\omega_d\big(D\varphi_h^t(w)\big)
  =
  \max_{u \in K}
  \omega_d\bigg(
  H(\varphi^t(u))
  D\varphi^t(u)
  \big(H(u)\big)^{-1}
  \bigg) <1,
\end{equation}
then by \eqref{wDphiht<1} with $H(u)$ instead of $Dh(u)$,
\eqref{dDOunderdiff} and \eqref{defT} for sufficiently large $t=T>0$
we have
\[
   \dim_{\rm H}K \leq d_{\rm L}(\{\varphi^t\}_{t\geq0},K)
   \leq d_{\rm L}(\varphi^T,K) \leq d.
\]
If it is considered $H(u)=p(u)S$,
where $p(u): U \subseteq \mathbb{R}^n \to \mathbb{R}^1$
is a continuous positive scalar function
and $S$ is a nonsingular $n \times n$ matrix,
then condition \eqref{wDphiH<1}
takes the form
\begin{equation}\label{wDphip<1}
  \sup_{u \in K}
  \omega_d\bigg(
  H(\varphi^t(u))
  D\varphi^t(u)
  \big(H(u)\big)^{-1}
  \bigg)
  =
  \sup_{u \in K}
  \bigg(
  \big(p(\varphi^t(u))p(u)^{-1}\big)^d
  \omega_d\big(SD\varphi^t(u)S^{-1}\big)
  \bigg)
  <1.
\end{equation}
\end{corollary}

\bigskip
Consider now the Leonov method of analytical estimation of the Lyapunov dimension
and its relation with the invariance of Lyapunov dimension with respect to diffeomorphisms.
Following \cite{Leonov-1991-Vest, Leonov-2002}, we consider the special class of diffeomorphisms
such that $Dh(u)=p(u)S$,
where $p(u): U \subseteq \mathbb{R}^n \to \mathbb{R}^1$
is a continuous scalar function and $S$ is a nonsingular $n\times n$ matrix.
As it is shown below the multiplier of the type $p(\varphi^t(u))(p(u))^{-1}$ in \eqref{wDphip<1}
plays the role of Lyapunov-like functions\footnote{
In \cite{NoackR-1996} it is interpreted as changes of Riemannian metrics.}.

Let us apply the linear change of variables
$w=h(u)=Su$ with a nonsingular $n \times n$ matrix $S$.
Then
$\varphi^t(u_0) = u(t,u_0)$
is transformed into $\varphi_S^t(w_0)$:
\[
   \varphi_S^t(w_0) = w(t,w_0)=S\varphi^t(u_0)=Su(t,S^{-1}w_0).
\]
Consider the transformed
systems \eqref{eq:ode} and \eqref{eq:odife}
\[
  \dot w = Sf(S^{-1}w)
  \mbox{\ or \ }
  w(t+1) = S\varphi(S^{-1}w(t))
\]
and their linearizations along the solution
$\varphi_S^t(w_0)=w(t,w_0)=S\varphi^t(u_0)$:
\begin{equation}\label{jacobian-new}
\begin{aligned}
  &
  \dot v = J_S(w(t,w_0))v
  \mbox{\ \ \it or \ \ }v(t+1) = J_S(w(t,w_0))v(t),
 \\ &
  J_S(w(t,w_0)) =S \, J(S^{-1}w(t,w_0)) \, S^{-1}
  = S \, J(u(t,u_0)) \, S^{-1}.
\end{aligned}
\end{equation}
For the corresponding fundamental matrices we have $D\varphi_S^t(w)=SD\varphi^t(u)S^{-1}$.

First we consider continuous time dynamical system.
Let $\lambda_i(u_0,S)=\lambda_i(S\varphi^t(u_0))$, $i=1,2,...,n$,
be eigenvalues of the  symmetrized Jacobian matrix 
\begin{equation}\label{SJS}
 \frac{1}{2} \left( S J(u(t,u_0)) S^{-1} +
 (S J(u(t,u_0)) S^{-1})^{*}\right)
 =
 \frac{1}{2} \left(
 J_S(w(t,w_0))+J_S(w(t,w_0))^{*}
 \right),
\end{equation}
ordered so that
$\lambda_1(u_0,S)\geq \cdots \geq \lambda_n(u_0,S)$ for any $u_0 \in U$.
The following theorem is rigorous reformulation
of results from \cite{Leonov-2002,Leonov-2008,Leonov-2012-PMM}.

\begin{theorem}\label{thm:LD-estimate-V}
Let $d=(j+s) \in [1,n]$, where
integer $j=\lfloor d \rfloor \in \{1,\ldots,n\}$
and real $s = (d - \lfloor d) \rfloor \in [0,1)$.
If there are a differentiable scalar function $V(u): U \subseteq \mathbb{R}^n \to \mathbb{R}^1$
and a nonsingular $n\times n$ matrix $S$
such that
\begin{equation}\label{ineq:weilSVct}
  \sup_{u \in K} \big(
  \lambda_1 (u,S) + \cdots + \lambda_j (u,S)
  + s\lambda_{j+1}(u,S) + \dot{V}(u)
  \big) < 0,
\end{equation}
where $\dot{V} (u) = ({\rm grad}(V))^{*}f(u)$,
then
\[
   \dim_{\rm H}K \leq
    d_{\rm L}(\{\varphi^t\}_{t\geq0},K)
    \leq
    d_{\rm L}(\varphi^T,K)
   \leq j+s
\]
for sufficiently large $T>0$.
\alarm{} 
\end{theorem}
{\bf Proof}.
From the following relations
(see Liouville's formula and, e.g., \cite[p.102]{Leonov-2008})
\begin{equation}\label{weil}
\begin{aligned}
  &
  \omega_{j+s}\big(SD\varphi^t(u)S^{-1}\big)
  =
   \exp\left(
     \int_0^t
     \lambda_1(S\varphi^\tau(u))+\cdots+\lambda_j(S\varphi^\tau(u))+s\lambda_{j+1}(S\varphi^\tau(u))
      d\tau
   \right)
\end{aligned}
\end{equation}
and
\[
  \big(p(\varphi^t(u))p(u)^{-1}\big)^{j+s}
  =
  \exp\big( V(\varphi^t(u))-V(u)\big)=
  \exp\left( \int_{0}^{t} \dot V(\varphi^\tau(u)) d\tau \right)
\]
we get
\begin{equation}
\begin{aligned}
  &
  \big(p(\varphi^t(u))p(u)^{-1}\big)^{j+s}
  \omega_{j+s}\big(SD\varphi^t(u)S^{-1}\big)
  \leq
  \\
  & \qquad
  \leq
   \exp\left(
     \int_0^t \big(
     \lambda_1(S\varphi^\tau(u))+\cdots+\lambda_j(S\varphi^\tau(u))+s\lambda_{j+1}(S\varphi^\tau(u))
     +\dot V(\varphi^\tau(u))
     \big) d\tau
   \right).
\end{aligned}
\end{equation}
Since $\varphi^t(u) \in K$ for any $u \in K$,
for $t>0$ by \eqref{ineq:weilSVct} we have
\[
  \max_{u \in K} \bigg(
  \big(p(\varphi^t(u))p(u)^{-1}\big)^{j+s}
  \omega_{j+s}\big(SD\varphi^t(u)S^{-1}\big)\bigg)
  < 1, \quad t>0.
\]
Therefore by Corollary~\ref{thm:hdiffLE} with
$H(u) = p(u)S$,
where
$p(u) = \left(e^{V(u)}\right)^{\frac{1}{d}}$,
we get the assertion of the theorem.
$\blacksquare$

Now consider discrete time dynamical system.
Let $\lambda_i(u_0,S)=\lambda_i(S\varphi^t(u_0))$, $i=1,2,...,n$,
be positive square roots of the eigenvalues of the  symmetrized Jacobian matrix 
\begin{equation}\label{SJS}
 \left( (S J(u(t,u_0)) S^{-1})^{*}S J(u(t,u_0)) S^{-1} \right)
 =
 \left(
 J_S(w(t,w_0))^{*}J_S(w(t,w_0))
 \right),
\end{equation}
ordered so that
$\lambda_1(u_0,S)\geq \cdots \geq \lambda_n(u_0,S)$ for any $u_0 \in U$.

\begin{theorem}\label{thm:LD-estimate-Vdt} \cite{Kuznetsov-2016-ArXiv}
Let $d=(j+s) \in [1,n]$, where
integer $j=\lfloor d \rfloor \in \{1,\ldots,n\}$
and real $s = (d - \lfloor d) \rfloor \in [0,1)$.
If there is a scalar continuous function $V(u): U \subseteq \mathbb{R}^n \to \mathbb{R}^1$
and a nonsingular $n\times n$ matrix $S$
such that
\begin{equation}\label{ineq:weilSV}
  \sup_{u \in K}\bigg(
  \ln\lambda_1 (u,S) + \cdots + \ln\lambda_j (u,S)
  + s\ln\lambda_{j+1}(u,S) + \big(V(\varphi(u))-V(u)\big)
  \bigg) < 0,
\end{equation}
then
\[
   \dim_{\rm H}K \leq
    d_{\rm L}(\{\varphi^t\}_{t\geq0},K)
    \leq
    d_{\rm L}(\varphi^T,K)
   \leq j+s
\]
for sufficiently large $T>0$.
\end{theorem}
{\bf Proof}.
By \eqref{horn}
for
$D\varphi_S^t(w)=SD\varphi^t(u)S^{-1} =
\prod\limits_{\tau=0}^{t-1}\big(S \, J(u(\tau,u_0)) \, S^{-1}\big)$
we have
\begin{equation}\label{weilds}
\begin{aligned}
  &
  \omega_{j+s}\big(SD\varphi^t(u)S^{-1}\big)
  \leq
  \prod_{\tau=0}^{t-1} \omega_{j+s}\big(S \, J(u(\tau,u_0)) \, S^{-1}\big).
\end{aligned}
\end{equation}
Therefore by the discrete analog of \eqref{weil} we have
\begin{equation}
  \omega_{j+s}\big(SD\varphi^t(u)S^{-1}\big)
  \leq
  \prod_{\tau=0}^{t-1} \lambda_1(S\varphi^\tau(u))\cdots\lambda_j(S\varphi^\tau(u))\big(\lambda_{j+1}(S\varphi^\tau(u))\big)^s.
\end{equation}
By the relation
\[
  \big(p(\varphi^t(u))p(u)^{-1}\big)^{j+s}
  =
  \exp\big( V(\varphi^t(u))-V(u)\big)
  =
  \exp\bigg( \sum_{\tau=0}^{t-1} V(\varphi^{\tau+1}(u))-V(\varphi^{\tau}(u)) \bigg)
\]
and we get
\[
\begin{aligned}
  &
  \ln\big(p(\varphi^t(u))p(u)^{-1}\big)^{j+s}
  +\ln\omega_{j+s}\big(SD\varphi^t(u)S^{-1}\big)
  \leq
  \\
  & \qquad
  \leq\sum_{\tau=0}^{t-1}\bigg(
  \ln\lambda_1(S\varphi^\tau(u))+\cdots+\ln\lambda_j(S\varphi^\tau(u))+s\ln\lambda_{j+1}(S\varphi^\tau(u))
  +V(\varphi(\varphi^{\tau}(u)))-V(\varphi^{\tau}(u))\bigg).
\end{aligned}
\]
Since $\varphi^t(u) \in K$ for any $u \in K$, by \eqref{ineq:weilSV}
and Corollary~\ref{thm:hdiffLE} with $H(u) = p(u)S$, where
$p(u) = \left(e^{V(u)}\right)^{\frac{1}{d}}$,
we get the assertion of the theorem.
$\blacksquare$

\medskip

 From \eqref{DOeqpoint} we have
\begin{corollary}
  If at an equilibrium point $u^{cr}_{eq} \equiv \varphi^{t}(u^{cr}_{eq})$
  for a certain $t>0$  the relation
  \[
    d_{\rm L}(\varphi^t,u^{cr}_{eq}) = j+s
  \]
  holds, then for any invariant set $K \ni u^{cr}_{eq}$
  we get analytical formula of exact Lyapunov dimension
  \[
  \begin{aligned}
   & \dim_{\rm H}K =
   d_{\rm L}(\{\varphi^t\}_{t\geq0},K)
   =d_{\rm L}(\{\varphi^t\}_{t\geq0},u^{cr}_{eq})=j+s.
  \end{aligned}
  \]
\end{corollary}

Remark that in the above approach we need only the Douady-Oesterl\'{e} theorem
(see Theorem~\ref{thm:DO}) and do not use the results on the Lyapunov dimension
developed by Eden, Constantin, Foias, Temam in \cite{ConstantinFT-1985,EdenFT-1991}
(see \eqref{dimLt},\eqref{dimL}, Propositions~\ref{thm:dimLdLKYglob} and \ref{thm:dimLdLKYloc}).

In \cite{BoichenkoL-1998,PogromskyM-2011} it is demonstrated, how a technique similar to the above can be effectively applied to derive constructive upper bounds of the topological entropy of dynamical systems.

\bigskip
\medskip
For the study of continuous time dynamical system in $\mathbb{R}^3$
the following result is useful.
Consider a certain open set $K_{\varepsilon} \subset U \subseteq \mathbb{R}^n$, which is diffeomorphic
to a ball, whose boundary $\partial \overline{K_{\varepsilon}}$
is transversal to the vectors $f(u)$, $u \in \partial \overline{K_{\varepsilon}}$.
Let the set $K_{\varepsilon}$ be a positively invariant for the solutions of system \eqref{eq:ode}.

\begin{theorem}(see, \cite{Leonov-1991-Vest,Leonov-2008}]
\label{thm:LD-zero}
Suppose a continuously differentiable function $V(u)$ and a non-degenerate matrix $S$
exist such that
\begin{equation}\label{pmm36}
  \lambda_1(u,S)+\lambda_2(u,S)+\dot V(u) <0,\,\,\forall\, u\in K_{\varepsilon}.
\end{equation}
Then any solution of system \eqref{eq:ode} with the initial data $u_0\in K_{\varepsilon}$
tends to the stationary set as $t\to+\infty$.
\end{theorem}

\section{Analytical formulas of exact Lyapunov dimension
for well-known dynamical systems}

Next we consider examples in which the critical point,
corresponding to the maximum of the local Lyapunov dimension,
is one of the equilibrium points (see \eqref{alldimeq}).
Let us consider several examples of smooth dynamical systems
generated by difference and differential equations
(for an example of PDE see, e.g. \cite{DoeringGHN-1987}).
In these examples we assume the existence of invariant set $K$
in which the corresponding dynamical system $\{\varphi^t\}_{t\geq0}$ is defined,
and use the compact notation $d_{\rm L}(K)$ for the Lyapunov dimension instead of \eqref{DOds}.

\subsection{Henon map}

Consider the Henon map $F \colon \mathbb{R}^2\to\mathbb{R}^2$
\begin{equation}\label{henonmap}
 \left(
 \begin{array}{l}
   x\\ y
 \end{array}
 \right)
 \to
 \left(
 \begin{array}{l}
   a+by-x^2 \\ x
 \end{array}
 \right),
\end{equation}
where $a>0$, $b\in(0,1)$ are the parameters of mapping.
The stationary points $(x_\pm,x_\pm)$ of this map are the following
$$
\begin{array}{l}
x_+=\frac{1}{2} \,\left[b-1+\sqrt{(b-1)^2+4a}\
\right],\\[10pt]
x_-=\frac{1}{2}
\,\left[b-1-\sqrt{(b-1)^2+4a}\ \right].
\end{array}
$$
\begin{theorem}\cite{Leonov-2002}
 For a bounded invariant set $K \ni (x_-,x_-)$ with respect to \eqref{henonmap} we have
 $$
  d_{\rm L}(K)
  =1+\frac{1}{1-\ln b/\ln\sigma_1(x_-)},
 $$
 where
 $$
  \sigma_1(x_-)=\sqrt{x_-^2+b}-x_-.
 $$
\end{theorem}

\subsection{Lorenz system}
Consider the classical Lorenz system suggested in \cite{Lorenz-1963}:
\begin{equation}\label{sys:Lorenz-classic}
 \left\{
 \begin{aligned}
 &\dot x= \sigma(y-x),\\
 &\dot y= r x-y-xz,\\
 &\dot z=-b z+xy,
 \end{aligned}
 \right.
\end{equation}
where
\[
 \sigma>0, \ r>0, \ b >0
\]
because of their physical meaning
(e.g., $b = 4(1+a^2)^{-1}$ is positive and bounded).

Since the system is dissipative and generates a dynamical system $\{\varphi^t\}_{t\geq 0}$
(to verify this, it is sufficiently to consider the Lyapunov function
$V(x,y,z) = \frac{1}{2}(x^2 + y^2 + (z - r - \sigma)^2)$;
see, e.g., \cite{Lorenz-1963,BoichenkoLR-2005}),
it possesses a global attractor
\cite{Chueshov-2002-book,BoichenkoLR-2005}.

\begin{theorem}\label{thm:LorenzLDformula}\cite{LeonovKKK-2015-arXiv-Lorenz}
 Assume  that the following inequalities
 \begin{equation} \label{cond:mainTheoremR0}
 r - 1 > 0,
 \end{equation}
 \begin{equation}\label{cond:mainTheoremR}
 r - 1 \ge \frac{b(b + \sigma - 1)^2 -
 4\sigma(b + \sigma b - b^2)}{3\sigma^2}
 \end{equation}
 are satisfied.
 Suppose that one of the following two conditions holds:
 \begin{itemize}
 \item[a.]
 \begin{equation}\label{cond:mainTheorem1}
 \sigma^2 (r - 1)(b - 4)
 \le 4\sigma (\sigma b + b - b^2) - b(b + \sigma - 1)^2;
 \end{equation}

 \item[b.]
 there are two distinct real roots of equations
 \begin{equation}\label{eq:mainTheoremGamma2}
 \begin{split}
 (2\sigma - b + \gamma)^2 \left(b(b+\sigma-1)^2
 - 4\sigma(\sigma b + b-b^2) + \sigma^2(r-1)(b-4)\right) + \\
 + 4b\gamma(\sigma+1)\left(b(b+\sigma-1)^2
 - 4\sigma(\sigma b + b-b^2) - 3 \sigma^2(r-1)\right) =0
 \end{split}
 \end{equation}
 and
 \begin{equation}\label{cond:mainTheorem2}
 \left\{
 \begin{gathered}
 \sigma^2 (r - 1)(b - 4) > 4\sigma
 (\sigma b + b - b^2) - b(b + \sigma - 1)^2 \hfill, \\
 \gamma^{(II)} > 0 \hfill \\
 \end{gathered}
 \right. \hfill\\
 \end{equation}
 where $\gamma^{(II)}$ is
 the greater root of equation \eqref{eq:mainTheoremGamma2}.
 \end{itemize}

 In this case we have:
 \begin{enumerate}
 \item If
 \begin{equation}\label{cond:mainTheoremR1}
 (b-\sigma)(b-1)< \sigma r < (b + 1)(b + \sigma),
 \end{equation}
 then any bounded on $[0; +\infty)$
 solution of system \eqref{sys:Lorenz-classic}
 tends to a certain equilibrium as $t \to +\infty$.

 \item If
 \begin{equation}\label{cond:mainTheoremR2}
 \sigma r > (b + 1)(b+\sigma),
 \end{equation}
 then for a bounded invariant set $K \ni (0,0,0)$
 \begin{equation}\label{LorenzLD-formula}
 d_{\rm L}(K)
 = 3 - \frac{2 (\sigma + b + 1)}{\sigma + 1
 + \sqrt{(\sigma-1)^2 + 4 \sigma r}}.
 \end{equation}
If $(0,0,0) \notin K$, then
the right-hand side of \eqref{LorenzLD-formula} is an upper bound of $d_{\rm L}(\{\varphi^t\}_{t\geq0},K)$.
 \end{enumerate}
\end{theorem}

The existence of analytical formula for the exact Lyapunov dimension
of the Lorenz system with classical parameters
is known (see, e.g. \cite{LeonovL-1993}) as the \emph{Eden conjecture} on the Lorenz system
(see \cite[p.411,Question 3.]{Eden-1989},\cite[p.98, Question 2.]{Eden-1989-PhD}, \cite{Eden-1990}).

\begin{remark}\label{r:unstableEquilibria}
 It can be easily checked numerically that if all three equilibria are
 hyperbolic (see the theorem in \cite{LeonovPS-2013-PLA}),
 then the conditions of Theorem~\ref{thm:LorenzLDformula} are satisfied.
 For example, for the standard parameters $\sigma = 10$ and $b=\frac{8}{3}$
 formula \eqref{LorenzLD-formula} is valid for $r > \frac{209}{45}$.
\end{remark}

\subsection{Glukhovsky-Dolzhansky system}
Consider a system, suggested by Glukhovsky and Dolghansky~\cite{GlukhovskyD-1980}
\begin{gather}
\begin{cases}
	\dot{x} $ = $ -\sigma x + z + a_0 yz, \\
	\dot{y} $ = $ R - y - xz, \\
    \dot{z} $ = $ -z + xy,
\end{cases} \label{sys:conv_fluid}
\end{gather}
where $\sigma$, $R$, $a_0$ are positive numbers (here $u=(x,y,x)$).
By the change of variables
\begin{equation}
(x,y,z) \rightarrow (x,R - \frac{\sigma}{a_0 R+1} z,\frac{\sigma}{a_0 R+1} y)
\label{sys:conv_fluid:change_var}
\end{equation}
system~\eqref{sys:conv_fluid} becomes
\begin{gather}
\begin{cases}
	\dot{x} $ = $ -\sigma x + \sigma y - \frac{a_0\sigma^2}{(a_0 R + 1)^2} y z,\\
	\dot{y} $ = $ \frac{R}{\sigma}(a_0 R + 1) x - y - x z, \\
	\dot{z} $ = $ -z + x y.
\end{cases} \label{sys:conv_fluid_lg}
\end{gather}

System \eqref{sys:conv_fluid_lg} is a generalization
of Lorenz system~\eqref{sys:Lorenz-classic}
and can be written as
\begin{gather}
\begin{cases}
	\dot{x} $ = $ \sigma(y-x) - A y z\\
	\dot{y} $ = $ r x - y - x z \\
    \dot{z} $ = $ -bz + x y,
\end{cases}
\label{sys:lorenz-general}
\end{gather}
where
\begin{equation}
	\quad A = \frac{a_0\sigma^2}{(a_0 R + 1)^2}, \quad
	r = \frac{R}{\sigma}(a_0 R + 1),\quad  b=1. \label{sys:conv_fluid:param}
\end{equation}

\begin{theorem}\cite{LeonovKM-2015-ArXiv} \label{conseq}
If
\begin{enumerate}
  \item $\sigma = Ar$, $4 \sigma r > (b+1) (b + \sigma)$ \\ or
  \item $b = 1$, $r > 2$, and
  \[
    \begin{cases}
      \sigma > \frac{-3 + 2\sqrt{3}}{3} \, A r,
      &\text{ if } \quad  2 < r \leq 4,\\
      \sigma \in \left(\frac{-3 + 2\sqrt{3}}{3} \, A r, \,
      \frac{3 r + 2 \sqrt{r (2 r + 1)}}{r - 4} \, A r \right) , &\text{ if } \quad  r > 4,
    \end{cases}
  \]
\end{enumerate}
then for a bounded invariant set $K \ni (0,0,0)$
of system \eqref{sys:lorenz-general} with $b=1$ or $\sigma = A r$
we have
\begin{equation}\label{GDLD-formula}
  d_{\rm L}(K)
  = 3 - \frac{2(\sigma + 2)}{\sigma + 1 + \sqrt{(\sigma - 1)^2 + 4 \sigma r}},
\end{equation}
If $(0,0,0) \notin K$, then
the right-hand side of the above relation
is an upper bound of $d_{\rm L}(\{\varphi^t\}_{t\geq0},K)$.

\end{theorem}
Note that this formula coincides with the formula for the
classical Lorenz system~\cite{Leonov-2012-PMM}.
Remark that system \eqref{sys:conv_fluid} is dissipative
and possesses a global attractor (see, e.g. \cite{LeonovKM-2015-EPJST}).

\subsection{Yang and Tigan systems}
Consider the Yang system \cite{Yangc-2008}:
\begin{equation}\label{sys:Yang}
	\left\{
	\begin{aligned}
		&\dot x= \sigma(y-x),\\
		&\dot y= r x - xz,\\
		&\dot z=-b z+xy,
	\end{aligned}
	\right.
\end{equation}
where $\sigma>0, b>0$, and $r$ is a real number.
Consider also the T-system (Tigan system) \cite{Tigan-2008}:
\begin{equation}\label{T-sys}
	\left\{
	\begin{aligned}
		&\dot x= a(y-x),\\
		&\dot y= (c-a) x - axz,\\
		&\dot z=-b z+xy.
	\end{aligned}
	\right.
\end{equation}
By the transformation
$(x,y,z) \rightarrow (\frac{x}{\sqrt{a}},\frac{y}{\sqrt{a}},\frac{z}{a})$
the Tigan system takes the form of the Yang system with parameters $\sigma = a, r = c-a$.

\begin{theorem}\label{Yang}\cite{LeonovKKK-2015-arXiv-YangTigan}

	\begin{enumerate}
		\item Assume $r = 0$ and the following inequalities $b(\sigma - b) > 0$, $\sigma - \frac{(\sigma + b)^2}{4(\sigma - b)} \ge 0$ are satisfied.
		Then any bounded on $[0; +\infty)$ solution of system \eqref{sys:Yang}
		tends to a certain equilibrium as $t \to +\infty$.
		
		\item Assume $r < 0$ and $ r \sigma + b(\sigma - b) > 0$.
		Then any bounded on $[0; +\infty)$ solution of system \eqref{sys:Yang}
		tends to a certain equilibrium as $t \to +\infty$.
		
		\item Assume $r > 0$ and there are two distinct real roots $\gamma^{(II)} > \gamma^{(I)}$ of equation
		\begin{equation}\label{eq:mainTheorem:rPos}
			4 b r \sigma^2 (\gamma + 2\sigma - b)^2 + 16 \sigma b \gamma (r \sigma^2 + b(\sigma + b)^2 - 4 \sigma (\sigma r + \sigma b - b^2)) =0
		\end{equation}
		such that $\gamma^{(II)} > 0$.
		
		In this case
		\begin{enumerate}
			\item if
			\begin{equation*}\label{cond:mainTheorem:rPosStable}
				b(b - \sigma) < r \sigma < b (\sigma + b),
			\end{equation*}
			then any bounded on $[0; +\infty)$ solution of system \eqref{sys:Yang}
			tends to a certain equilibrium as $t \to +\infty$.
			
			\item if
			\begin{equation}\label{cond:mainTheorem:sPosFormula}
				r \sigma > b (\sigma + b),
			\end{equation}
			then
			\begin{equation}\label{YangTiganLD-formula}
                d_{\rm L}(K)
                = 3 - \frac{2 (\sigma + b)}{\sigma + \sqrt{{\sigma}^2 + 4 \sigma r}},
			\end{equation}
			where $K \ni (0,0,0)$ is a bounded invariant set of system \eqref{sys:Yang}.
If $(0,0,0) \notin K$, then
the right-hand side of the above relation
is an upper bound of $d_{\rm L}(\{\varphi^t\}_{t\geq0},K)$.
		\end{enumerate}
	\end{enumerate}  	
\end{theorem}


	
\subsection{Shimizu-Morioka system}
Consider the Shimizu-Morioka system \cite{Shimizu1980201} of the form

\begin{equation}\label{eq18}\begin{aligned}
&\dot x=y,\\
&\dot y=x-\lambda y - xz,\\
&\dot z=-\alpha z+x^2,\end{aligned}
\end{equation}
where $ \alpha, \lambda $ are positive parameters.

Using the diffeomorphism
\begin{equation}\label{eq19}
 \left(
 \begin{array}{l}
   x\\ y \\z
 \end{array}
 \right)
 \to
 \left(
 \begin{array}{l}
   x \\ y \\ z-\frac{x^2}{2}
 \end{array}
 \right),
\end{equation}
system \eqref{eq18} can be reduced to the following system
\begin{equation}\label{eq20}\begin{aligned}
&\dot x=y,\\
&\dot y=x-\lambda y - xz + \frac{x^3}{2},\\
&\dot z=-\alpha z + xy + \Bigl (1+\frac{\alpha}{2}\Bigr)x^2,\end{aligned}
\end{equation}
where $ \alpha, \lambda $ are the positive parameters of system \eqref{eq18}.
We say that system \eqref{eq20} is a {\it transformed}
Shimizu-Morioka system.

\begin{theorem}\cite{LeonovAK-2015}
Suppose, $K$ is a bounded invariant set of system \eqref{eq20}:
$(0,0,0)\in K$, and the following relations
\begin{equation}\label{eq22}
\lambda-4 \le \sqrt{10+\frac{3}{\alpha}-13\alpha},
\lambda < \frac{1}{\alpha}-\alpha,
4-\lambda \le \sqrt{\frac{8+15\alpha-8{\alpha}^2-24{\alpha}^3}{2\alpha(\alpha+1)}}
\end{equation}
are satisfied.
Then
\begin{equation}\label{eq23}
   d_{\rm L}(K)
   = 3-\frac{2(\lambda+\alpha)}{\lambda+\sqrt{4+{\lambda}^2}} .
\end{equation}
If $(0,0,0) \notin K$, then
the right-hand side of relation \eqref{eq23}
is an upper bound of $d_{\rm L}(\{\varphi^t\}_{t\geq0},K)$.
\end{theorem}
In the proof there are used the Lyapunov function of the form
$$
  V(x,y,z)=\frac{1-s}{4\sqrt{4+\lambda^2}}\vartheta,
$$
where
$$
\vartheta=\mu_1(2y^2-2xy-x^4+2x^2z)+\mu_2x^2-\frac{4}{\alpha}z
+\mu_3(z^2-x^2z+\frac{x^4}{4}+xy)+\mu_4(z^2+y^2-\frac{x^4}{4}-x^2),
$$
and the nonsingular matrix
$$
S=\begin{pmatrix} -\frac{1}{k}& 0 & 0\\
\lambda-\alpha & 1 & 0\\
0& 0 & 1\end{pmatrix}.
$$

\section{Attractors of dynamical systems}

Compact invariant sets of dynamical systems are related with the notions of attractors
(see, e.g. \cite{Ladyzhenskaya-1987,Ladyzhenskaya-1991,BabinV-1992,Chueshov-1993,Temam-1997,Chueshov-2002-book,BoichenkoLR-2005,Leonov-2008}).
Consider dynamical system
$\big(\{\varphi^t\}_{t\geq0},(U\subseteq \mathbb{R}^n,||\cdot||) \big)$.

\begin{property}\label{property:local_attr_set}
An invariant set $K \subset U \subseteq \mathbb{R}^n$ is said to be {\it locally attractive}
if for a certain $\varepsilon$-neighborhood of the set $K$: $K_\varepsilon \subseteq U$,
$$
 \lim_{t \to +\infty} \rho (K, \varphi^t(u)) = 0, \quad \forall ~
 u \in K_\varepsilon.
$$
Here $\rho(K, {u})$ is the distance from the point ${u}$
to the set $K$, defined as
$$
 \rho(K, {u}) = \inf_{{w} \in K} || {w} - {u}
||,
$$
and $K_\varepsilon$ is the set of points ${u}$ for which
$\rho (K, {u}) < \varepsilon$.
\end{property}

\begin{property}\label{property:global_attr_set}
An invariant set $K \subset U \subseteq \mathbb{R}^n$ is said to be {\it globally attractive}
if
$$
 \lim_{t \to +\infty} \rho (K, \varphi^t(u)) = 0, \quad \forall ~
 u \in U. 
$$
\end{property}

\begin{property}\label{property:unif_loc_attr_set}
An invariant set $K \subset U \subseteq \mathbb{R}^n$ is said to be
{\it uniformly locally attractive} with respect
to the dynamical system $\{\varphi^t\}_{t\geq0}$
if for a certain $\varepsilon$-neighborhood
$K_\varepsilon \subseteq U$,
any number $\delta > 0$, and any bounded set
$B \subseteq U \subseteq \mathbb{R}^n$,
there exists a number $t(\delta, B) > 0$ such that
$$
 \varphi^t(B \cap K_\varepsilon) \subset K_\delta,
 \quad \forall ~ t \geq t(\delta, B).
$$
Here
$$
 \varphi^t(B \cap K_\varepsilon) =
 \left\{\varphi^t({u}_0) ~|~ {u}_0 \in B \cap K_\varepsilon \right\}.
$$
\end{property}

\begin{property}\label{property:unif_glob_attr_set}
Invariant set $K \subset U \subseteq \mathbb{R}^n$ is said to be
{\it uniformly globally attractive} with respect to the dynamical system $\{\varphi^t\}_{t\geq0}$
if for any number $\delta > 0$
and any bounded set $B \subseteq U \subseteq \mathbb{R}^n$
there exists a number $t(\delta, B) > 0$ such that
$$
 \varphi^t(B) \subset K_\delta, \quad \forall ~ t \geq t(\delta, B).
$$
\end{property}

\begin{definition}\label{def:attractor}
For a dynamical system, a bounded closed
invariant set $K$ is
\begin{enumerate}[label=(\arabic*)]
 \item an \emph{attractor} if it is a locally attractive set
       (i.e., it satisfies Property~\ref{property:local_attr_set});
 \item a \emph{global attractor} if it is a globally attractive set
       (i.e., it satisfies Property~\ref{property:global_attr_set});
 \item a \emph{B-attractor} if it is a uniformly locally attractive set
       (i.e., it satisfies Property~\ref{property:unif_loc_attr_set}); or
 \item a \emph{global B-attractor} if it is a uniformly globally attractive set
       (i.e., it satisfies Property~\ref{property:unif_glob_attr_set}).
\end{enumerate}
\end{definition}

\begin{remark}\label{remark:minimal-attr}
In the above definition we assume the closedness for the sake of uniqueness.
The reason is that the closure of a locally attractive invariant set $K$ is also a
locally attractive invariant set
(for example, consider an attractor with excluded one of the embedded unstable periodic orbits).
Note that if a dynamical system  is defined for negative $t$,
then a locally attractive invariant set contains only whole trajectories, i.e.
if $u_0 \in K$, then $\varphi^t(u_0) \in K$ for $t \in \mathbb{R}$ (see \cite{Chueshov-2002-book}).
\end{remark}
\smallskip

\begin{remark}\label{remark:minimal-attr}
The definition under consideration implies that a global B-attractor
is also a global attractor (and an attractor).
Consequently, it is rational to introduce the notion of a
{\it minimal global attractor} (and {\it a minimal attractor})
\cite{Chueshov-1993,Chueshov-2002-book}.
This is the minimal bounded closed invariant set that
possesses Property~\ref{property:global_attr_set}
(or Property~\ref{property:local_attr_set}, i.e.
minimal local attractor is an attractor, which cannot be represented as a union of local attractors).
Further,  "global attractor" means "minimal global attractor".
\end{remark}
\smallskip

\begin{definition}
For an attractor $K$, the {\it basin of attraction}
is the set $\beta(K) \subseteq U \subseteq \mathbb{R}^n$
of all $u_0 \in U$ such that
\[
 \lim_{t \to +\infty} \rho (K, \varphi^t(u_0)) = 0.
\]
\end{definition}

\subsection{Computation of attractors and Lyapunov dimension}

The study of a dynamical system typically begins with an analysis of the equilibria,
which are easily found numerically or analytically.
Therefore, from a computational perspective, it is natural to suggest the following classification
of attractors, which is based on the simplicity of finding their basins of attraction in the phase space:

\begin{definition}
\cite{LeonovKV-2011-PLA,LeonovKV-2012-PhysD,LeonovK-2013-IJBC,LeonovKM-2015-EPJST}
 An attractor is called a \emph{self-excited attractor}
 if its basin of attraction
 intersects with any open neighborhood of a stationary state (an equilibrium),
 otherwise it is called a \emph{hidden attractor}.
\end{definition}

\begin{figure}[!ht]
\centerline{
\begin{tabular}{c c c}
\includegraphics[width=0.33\textwidth]{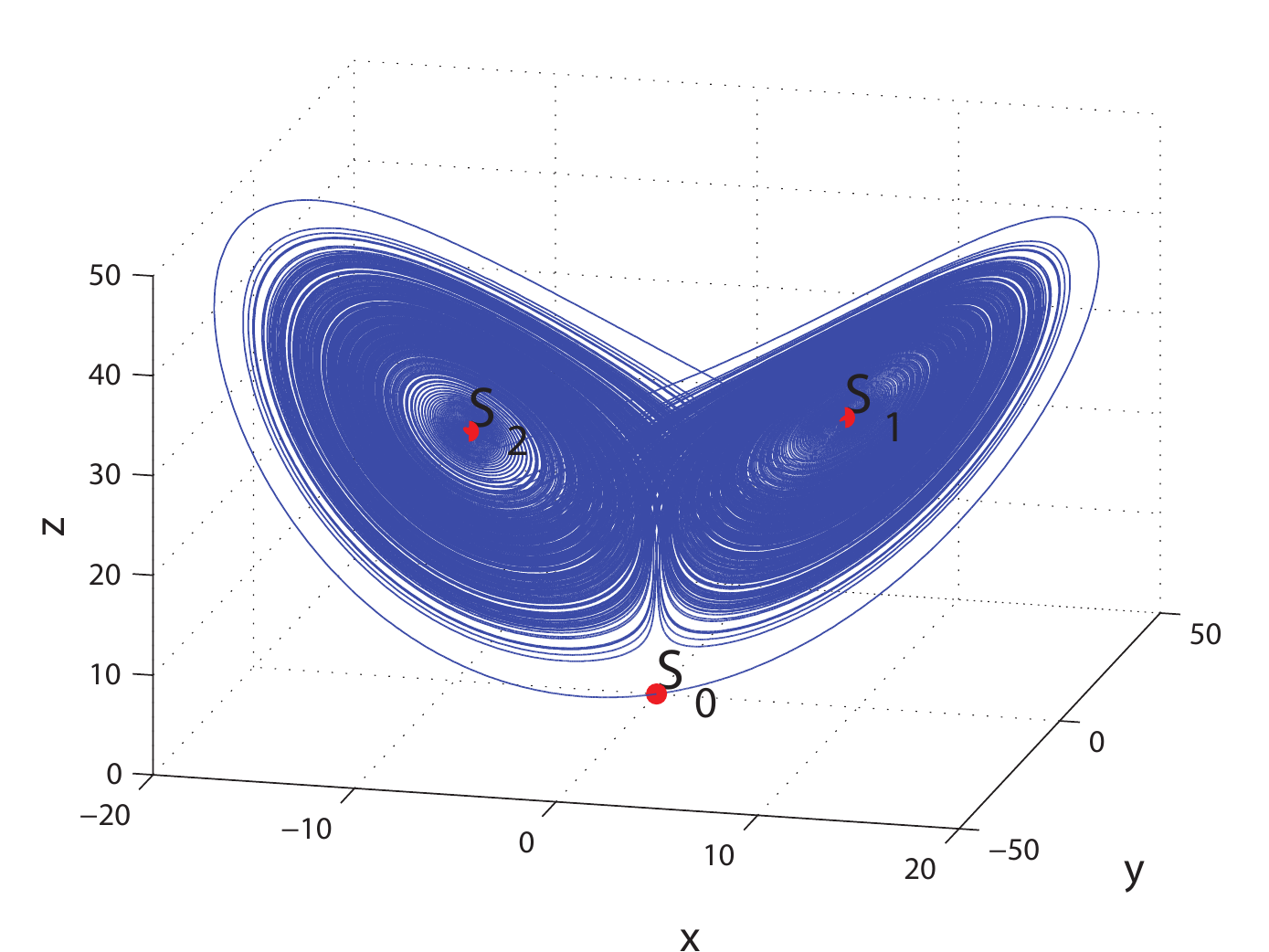} &
\includegraphics[width=0.33\textwidth]{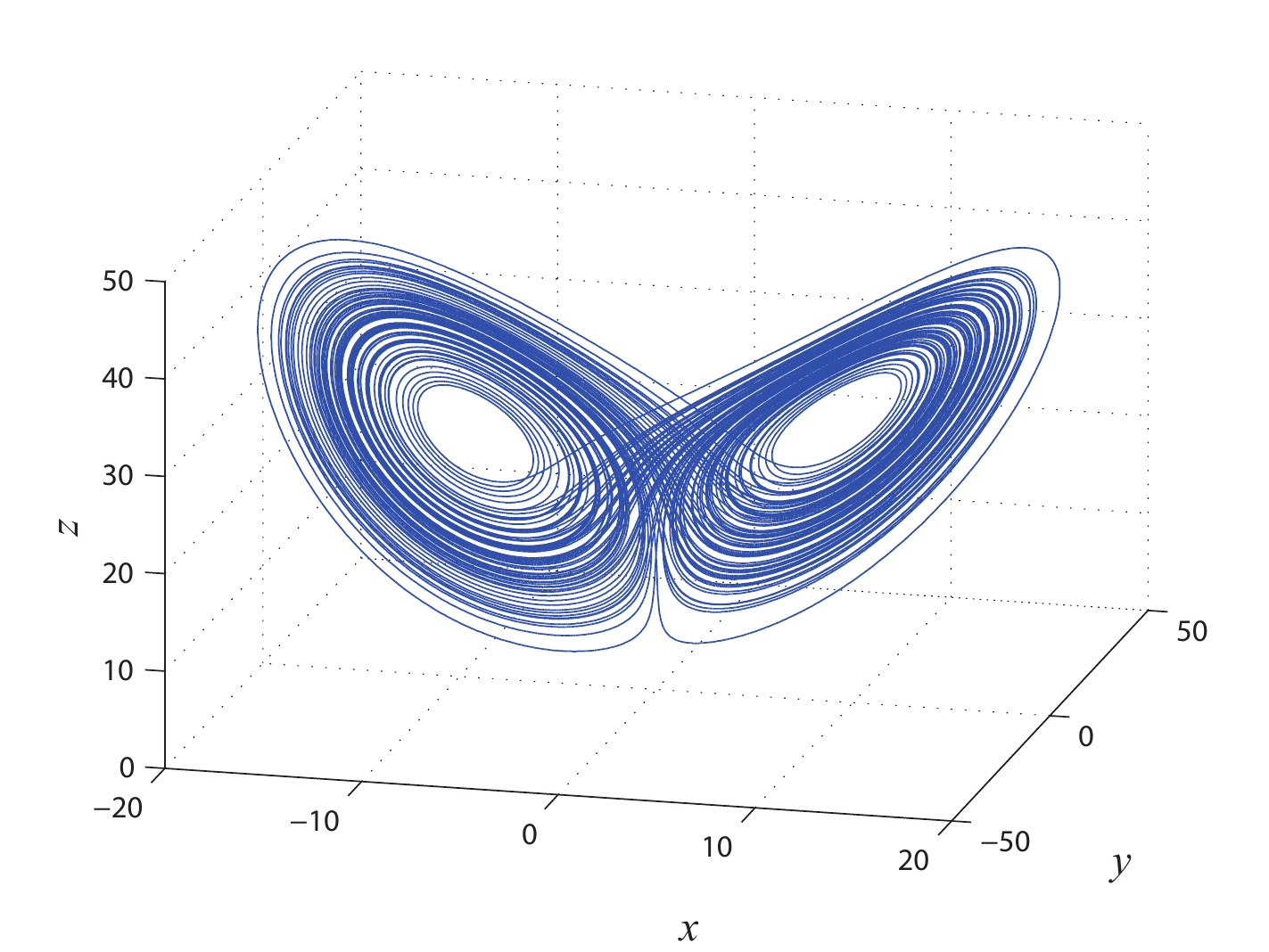}
\end{tabular}}
 \caption{
 Numerical visualization of
 self-excited chaotic attractor in the Lorenz system.
 Global B-attractor (left subfigure),
 $d_{\rm L}(K)=\sup_{u \in K}d_{\rm L}(u) = d_{\rm L}(S_0) = 2.4013$
 according to \eqref{LorenzLD-formula};
 global attractor (right subfigure),
 $d_{\rm L}(K) \approx 2.0565$ by numerical computation.
 Parameters: $r = 28$, $\sigma = 10$, $b = 8/3$.
 }
\label{attr-lorenz-3unstable}
\end{figure}

 Self-excited attractor in a system
 can be found using the standard computational procedure, i.e.
 by constructing a solution using initial data from a small neighborhood of the equilibrium,
 observing how it is attracted and, thus, visualizes the attractor.
 For example, in the Lorenz system~\eqref{sys:Lorenz-classic}
 with classical parameters $\sigma=10, \beta=8/3, \rho=28$
 there is a chaotic attractor, which is self-excited with respect to all three equilibria
 and could have been found using the standard computational procedure
 with initial data in vicinity of any of the equilibria
 (see Fig.~\ref{attr-lorenz-3unstable}).
Here it is possible to check numerically that for the considered parameters
the local attractor is a global attractor
(i.e. there are no other attractors in the phase space).
In this case the global B-attractor involves the chaotic local attractor,
three unstable equilibria and their unstable manifolds
attracted to the chaotic local attractor.

\begin{figure}[!ht]
\centerline{
\begin{tabular}{c c c}
\includegraphics[width=0.33\textwidth]{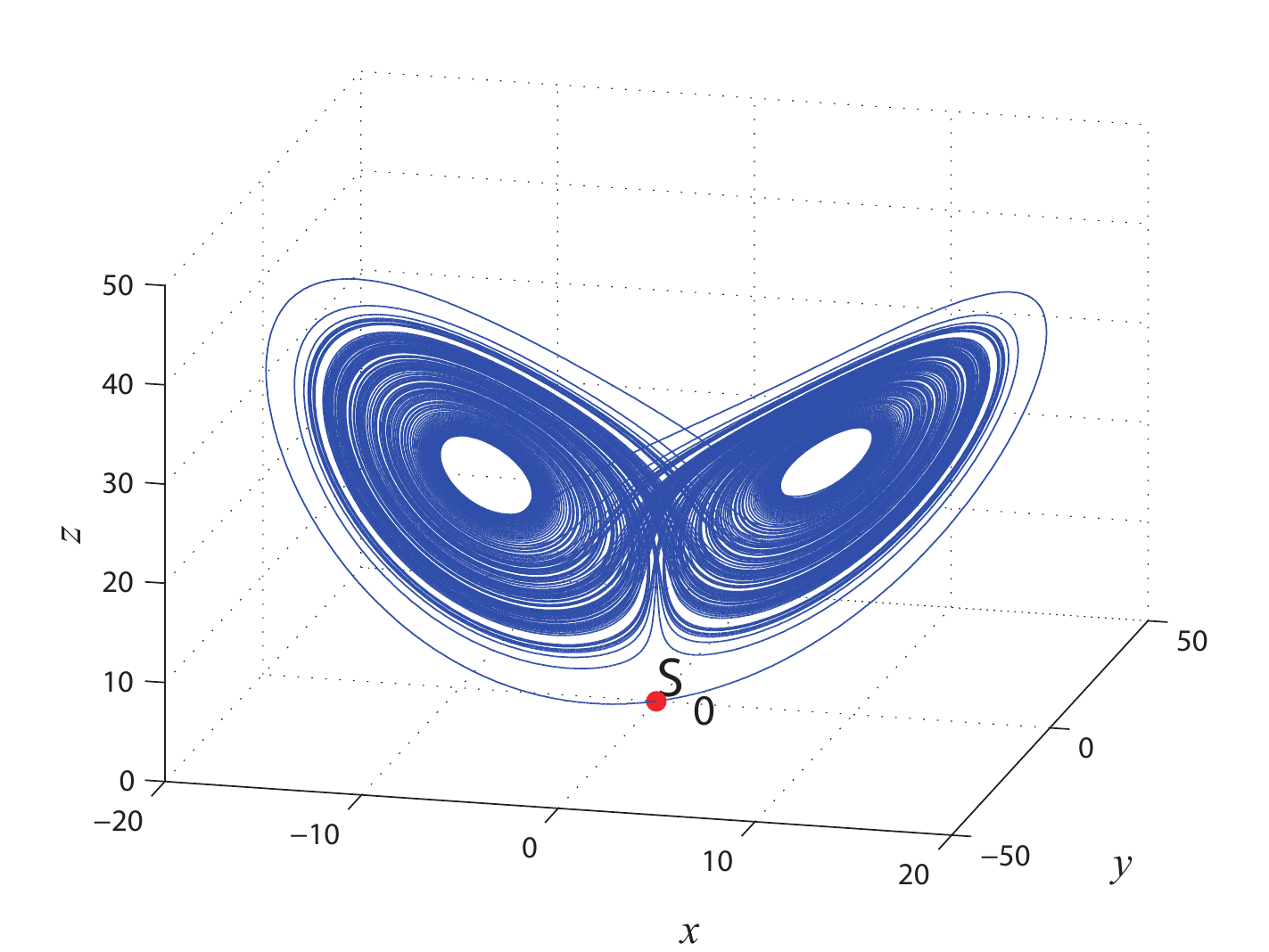} &
\includegraphics[width=0.33\textwidth]{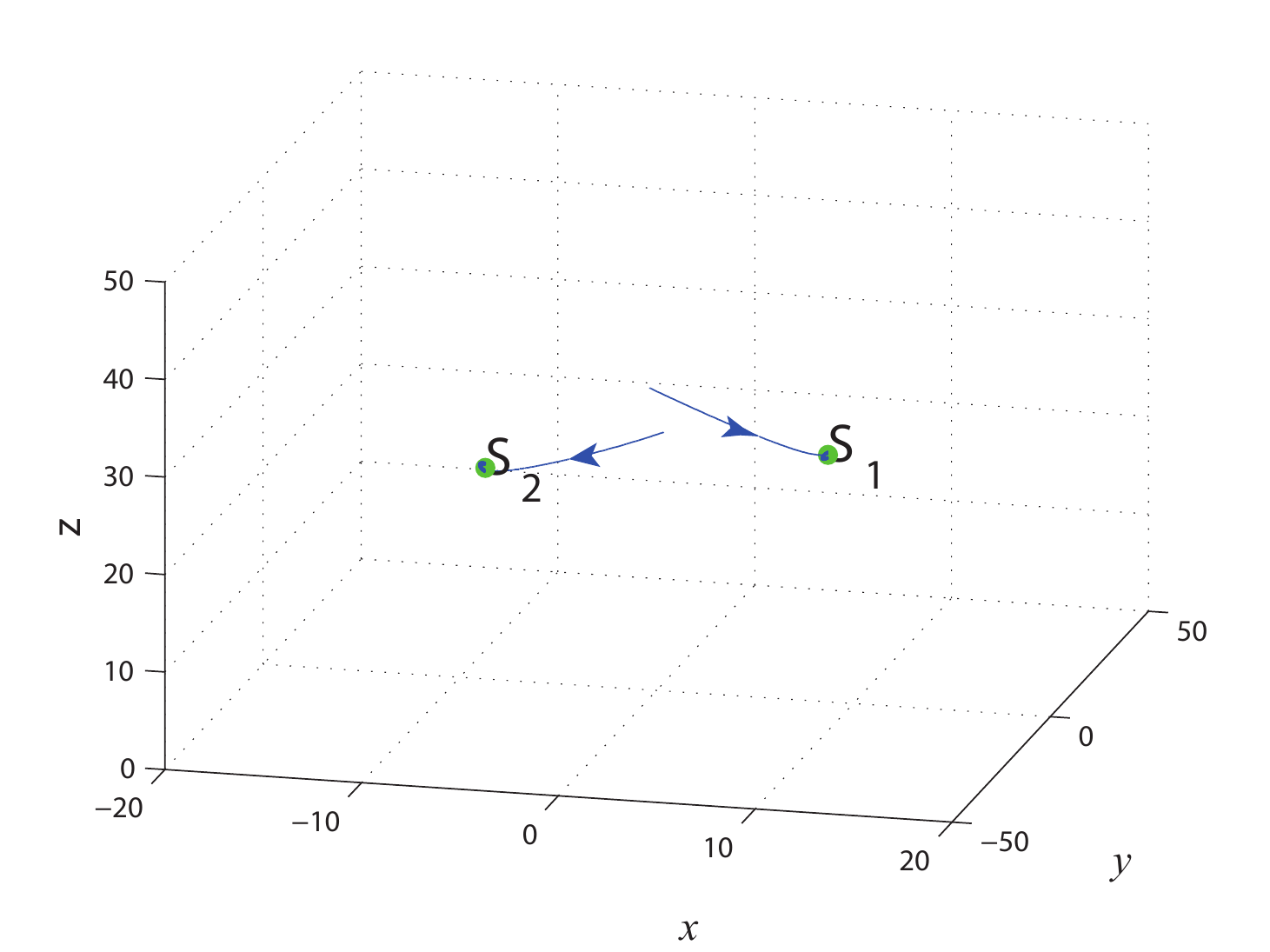} &
\includegraphics[width=0.33\textwidth]{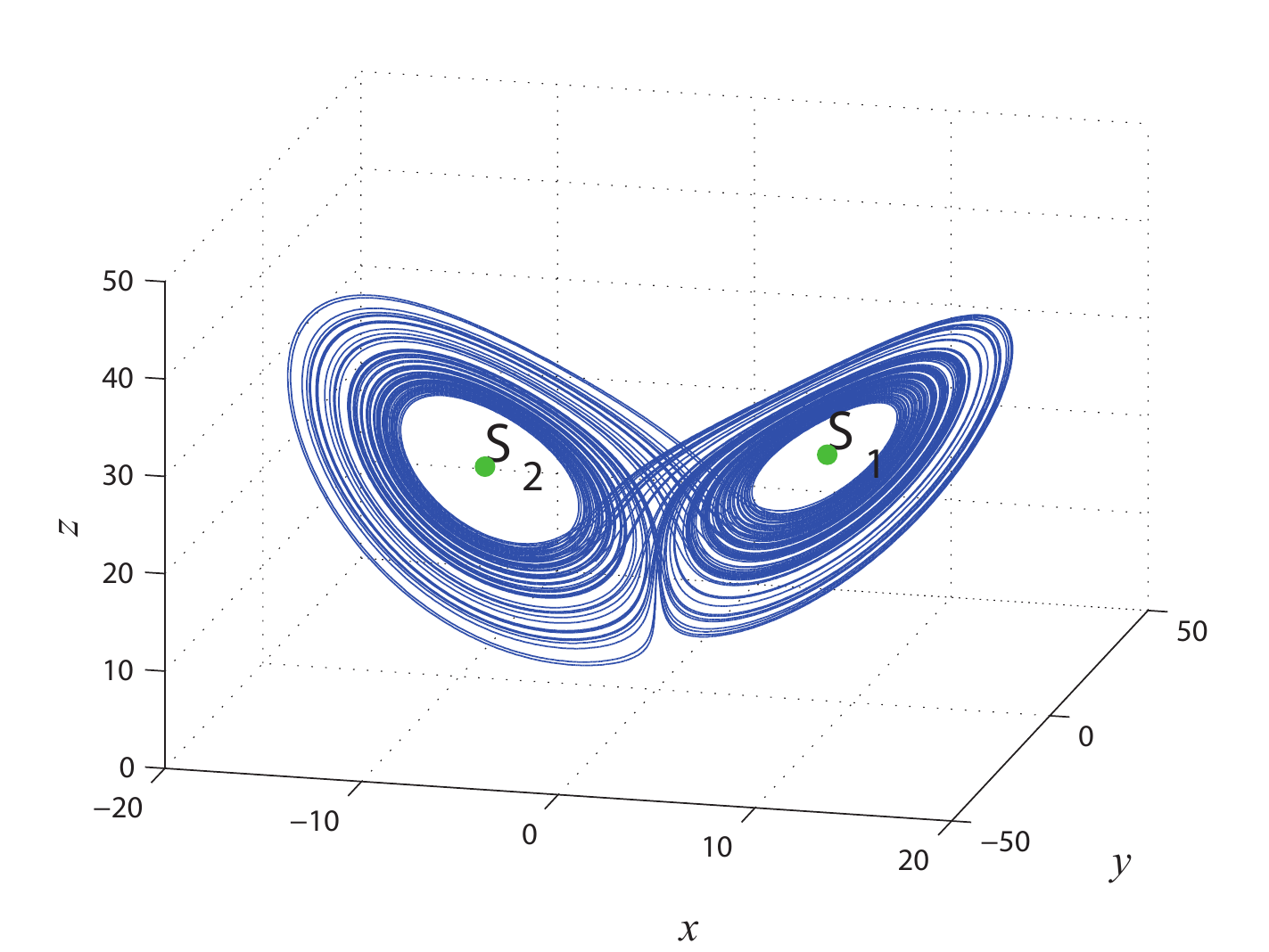}
\end{tabular}}
 \caption{
 Numerical visualization of
 self-excited chaotic local attractor in the Lorenz system.
 Local B-attractor involves self-excited chaotic local attractor,
 unstable zero equilibrium and its unstable manifold attracted
 to the chaotic local attractor (left subfigure),
 $d_{\rm L}(K)=\sup_{u \in K}d_{\rm L}(u)=d_{\rm L}(S_0)=2.3727$
 according to \eqref{LorenzLD-formula}.
 Trajectories with the initial data
 $(\pm 1.3276,\mp 9.7014, 28.7491)$
 tend to
 trivial local attractors --- equilibria $S_{2,1}$ (middle subfigure),
 $d_{\rm L}(S_{2,1})=1.9989$.
 Global attractor is the union of three coexisting local attractors:
 self-excited chaotic local attractor
 and two trivial local attractors (right subfigure),
 $d_{\rm L}(K)\approx2.0489$ by numerical computation.
 Parameters: $r = 24.5$, $\sigma = 10$, $b = 8/3$.
 }
\label{attr-lorenz-1unstable}
\end{figure}

However it is known that for other values of parameters,
e.g. $\sigma=10, \beta=8/3, \rho=24.5$ \cite{Sparrow-1982},
the chaotic local attractor
in the Lorenz system may be self-excited with respect to the zero unstable equilibrium only.
In this case there are three coexisting minimal local attractors (see Fig.~\ref{attr-lorenz-1unstable}):
chaotic local attractor and two trivial local attractors --- stable equilibria $S_{1,2}$.

Self-excited attractors in a multistable system can be found using
the standard computational procedure, whereas there is no standard way of predicting
the existence of hidden attractors in a system.

While the \emph{multistability} is a property of system,
the \emph{self-excited} and \emph{hidden} properties
are the properties of attractor and its basin.
For example, hidden attractors are attractors in systems
with no equilibria or with only one stable equilibrium
(a special case of multistability and coexistence of attractors).

\begin{figure}[!h]
 \centering
\begin{tabular}{c c}
 \includegraphics[width=0.4\textwidth]{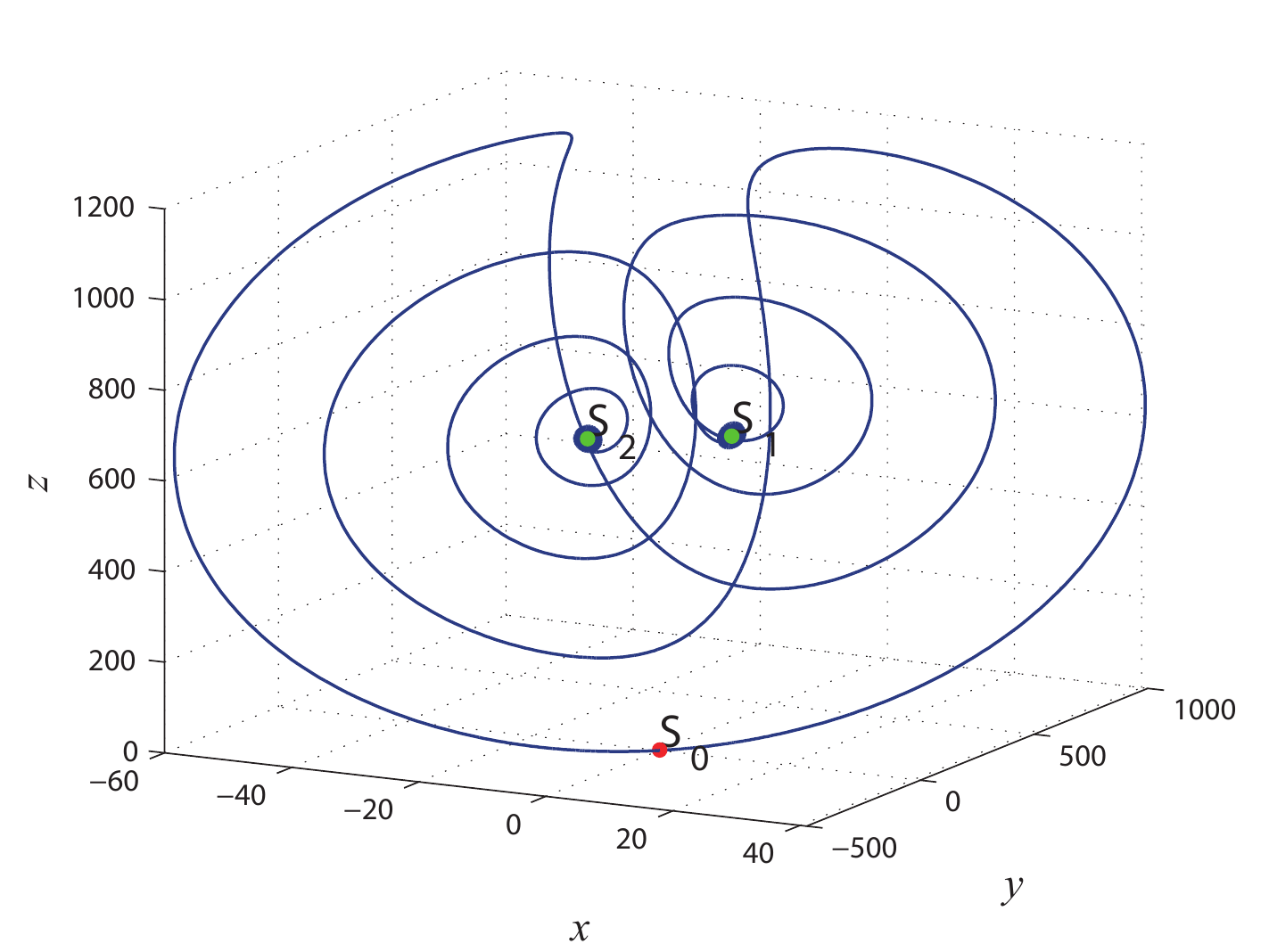} &
 \includegraphics[width=0.4\textwidth]{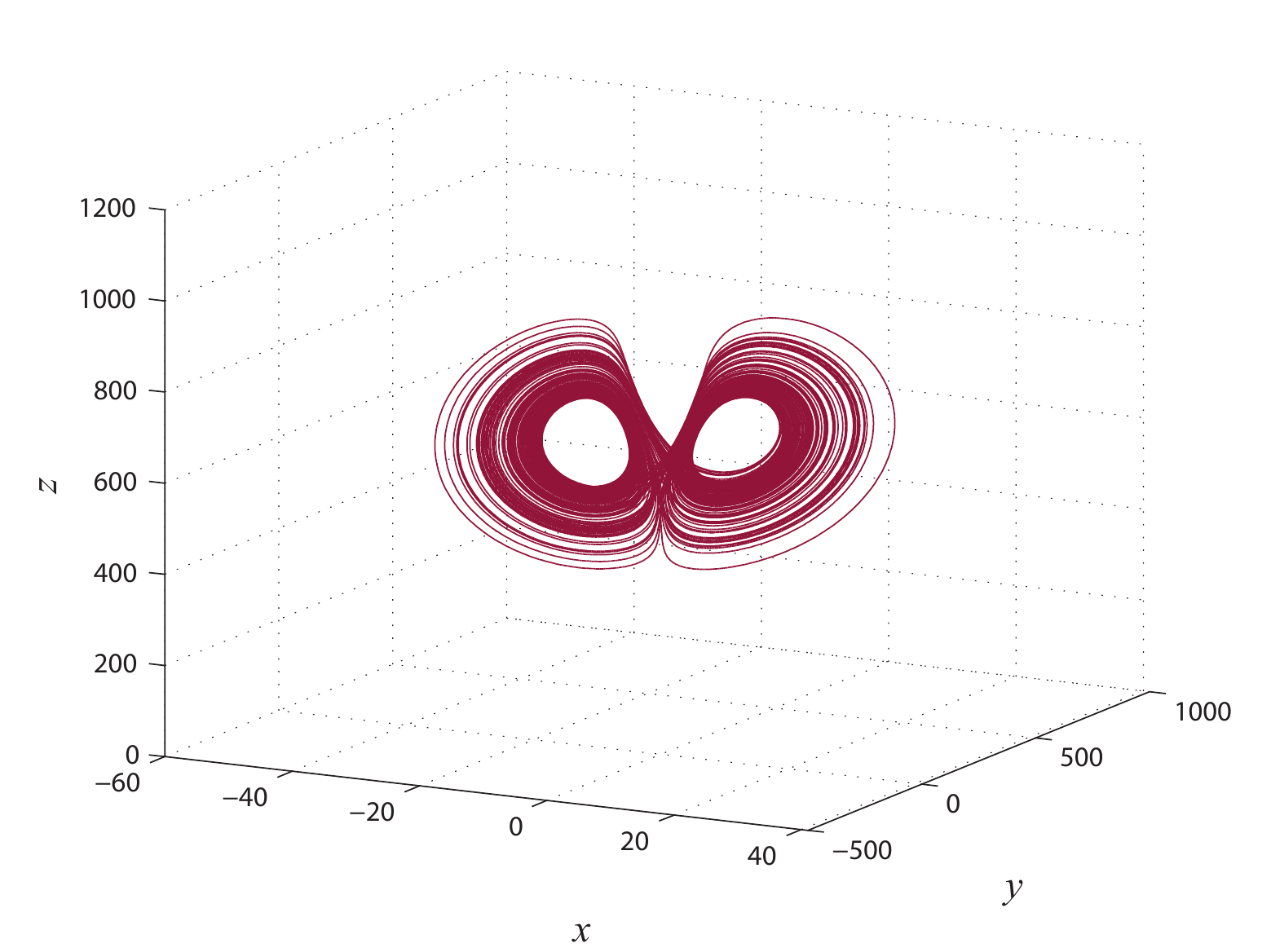}
\end{tabular}
 \caption{\label{fig:gd-hidden}
 Numerical visualization of local B-attractor and
 hidden local attractor in the Glukhovsky-Dolghansky system.
 Local B-attractor involves outgoing separatrix (blue) of the saddle $S_0$ (red)
 attracted to the stable equilibria $S_{1,2}$ (green) (left subfigure).
 Hidden local attractor
 (magenta,
 $d_{\rm L}(K) \approx 2.1322$
 by numerical computation)
 coexists with
 local B-attractor
 ($d_{\rm L}(K)=\sup_{u \in K}d_{\rm L}(u)=d_{\rm L}(S_0)=2.8917$
 by \eqref{GDLD-formula}).
 Global B-attractor involves the local B-attractor
 and the hidden local attractor.
 }
\end{figure}

In general, there is no straightforward way of predicting the existence
or coexistence of hidden attractors in a system
(see, e.g. \cite{KuznetsovLV-2010-IFAC,LeonovKV-2011-PLA,LeonovKV-2012-PhysD,LeonovK-2013-IJBC,KuznetsovL-2014-IFACWC,LeonovKM-2015-EPJST,LeonovKM-2015-CNSNS,KuznetsovLM-2015,DancaFKC-2016-IJBC,Kuznetsov-2016}).
A numerical search of hidden attractors by evolutionary algorithms
is discussed in \cite{Zelinka-2015,Zelinka-2016-HA}.
Recent examples of hidden attractors can be found in
\emph{The European Physical Journal Special Topics: Multistability: Uncovering Hidden Attractors}, 2015
(see \cite{ShahzadPAJH-2015-HA,BrezetskyiDK-2015-HA,JafariSN-2015-HA,ZhusubaliyevMCM-2015-HA,SahaSRC-2015-HA,Semenov20151553,FengW-2015-HA,Li20151493,FengPW-2015-HA,Sprott20151409,Pham20151507,VaidyanathanPV-2015-HA,SharmaSPKL-2015-EPJST}).

For example, in the Glukhovsky-Dolghansky system
and the corresponding generalized Lorenz system \eqref{sys:lorenz-general}
with parameters $r = 700, a = 0.0052, \sigma=ra, b=1$
a hidden chaotic local attractor can be found \cite{LeonovKM-2015-CNSNS,LeonovKM-2015-EPJST}
(see Fig.~\ref{fig:gd-hidden}).

Remark that if a system is proved to be dissipative
(i.e. it is possible to determine an absorbing bounded domain in the phase space
such that all trajectories enter this domain within a finite time),
then all self-excited or hidden local attractors of the system
are inside this absorbing bounded domain and can be found numerically.
However, in general,
\emph{the determination of the number and mutual disposition
of chaotic minimal local attractors in the phase space for a system}
may be a challenging problem \cite{LeonovK-2015-AMC}
(see, e.g. the corresponding well-known problem for two-dimensional polynomial systems
--- the second part of 16th Hilbert problem on the number and mutual disposition of limit cycles  \cite{Hilbert-1901})\footnote{
The numerical search of hidden attractors can be complicated by
the small size of the basin of attraction
with respect to the considered set of parameters $p \in P$
and subset of the phase space $U_0\subseteq U$:
following \cite{ZakrzhevskySY-2007,BrezetskyiDK-2015-HA},
the attractor may be called a \emph{rare attractor}
if the measure $\mu$ of the basin of attractors $\beta(K_p)$ for
the considered set of parameters $p \in P$ is small
with respect to the considered part of the phase space $U_0 \subseteq U$, i.e.
$\frac{\int_{p \in P} \mu(\beta(K_p)\cap U_0)}{\mu(U_0)} << 1$.
Also computational difficulties may be caused by the shape of basin of attraction,
e.g. by Wada and riddled basins.
}.
Thus the advantage of the analytical method for the Lyapunov dimension estimation,
suggested in Theorem~\ref{thm:LD-estimate-V},
is that it is useful not only for the dissipative systems
(see, e.g. estimation of the Lyapunov dimension for one of the Rossler systems \cite{LeonovB-1992})
but also allows one to estimate the Lyapunov dimension
of invariant set without localization of the set in the phase space.

Remark that, from a computational perspective,
it is not feasible to numerically check Property~\ref{property:local_attr_set}
for all initial states of the phase space of a dynamical system.
A natural generalization of the notion of an attractor is
the consideration of the weaker attraction requirements: almost everywhere or
on a set of positive measure (see, e.g., \cite{Milnor-2006}).
See also \emph{trajectory attractors} \cite{Sell-1996,ChepyzhovV-2002,ChueshovS-2005}.
In numerical computations,
to distinguish an artificial computer generated chaos
from a real behavior of the system, one can consider
the shadowing property of the system (see, e.g., the survey in \cite{Pilyugin-2011}).

We can typically see an attractor (or global attractor)
in numerical experiments.
The notion of a B-attractor is mostly used in the theory of dimensions,
where we consider invariant sets covered by balls.
The uniform attraction requirement
in Property~\ref{property:unif_loc_attr_set}
implies that a global B-attractor involves
a set of stationary points $\mathcal{S}$ and the corresponding unstable manifolds
$W^u (\mathcal{S}) = \left\{{u}_0 \in \mathbb{R}^n ~|~
\lim_{t \to -\infty} \rho (\mathcal{S}, \varphi^t(u_0)) = 0\right\}$
(see, e.g., \cite{Chueshov-1993,Chueshov-2002-book}).
The same is true for B-attractor
if the considered neighborhood $K_\varepsilon$
in Property~\ref{property:unif_loc_attr_set}
contains some of the stationary points from $\mathcal{S}$.
This allows one to get
analytical estimations
of the Lyapunov dimension for B-attractors
and even formulas
since the local Lyapunov dimension at a stationary point
can be easily obtained analytically
(but this does not help
for chaotic minimal local attractors,
hidden B-attractors
since they do not involve any stationary points).

From a computational perspective,
numerical check of Property~\ref{property:unif_loc_attr_set}
is also difficult.
Therefore, if the basin of attraction involves unstable manifolds of equilibria,
then computing the minimal attractor and the
unstable manifolds that are attracted to it may be regarded
as an approximation of minimal B-attractor.
For example, consider the visualization of the classical Lorenz attractor
from the neighborhood of the zero saddle equilibria.
Note that a minimal global attractor involves
the set $\mathcal{S}$ and
its basin of attraction involves the set $W^u (\mathcal{S})$.

For the computation of the Lyapunov dimension of an attractor $A$
we consider a sufficiently large time $T$
and a sufficiently dense grid of points $A_{\rm grid}$ on the attractor,
compute the local Lyapunov dimensions by
the corresponding Kaplan-Yorke formula  $d_{\rm L}^{\rm KY}(\{\LEs_i(T,u)\}_{1}^{n})$,
and take maximum on the grid: $\max_{u \in A_{\rm grid}}d_{\rm L}^{\rm KY}(\{\LEs_i(T,u)\}_{1}^{n})$.

Since numerically we can check only that all points of the grid belong
to the basin of attraction, the following remak is useful.
Let a point $u_0$ belongs to the basin of attraction of attractor $A$.
Consider the union of the semi-orbit
$\gamma^{+}(u_0) = \{\varphi^t(u_0), t \geq 0 \}$
and attractor $A$: $K(u_0)=A \cup \gamma^{+}(u_0)$.
According to the definition of the basin of attraction,
$\omega$-limit set of $\varphi^t(u_0)$ belong to $A$,
thus the set $K(u_0)$ is compact and invariant.
Since $A \supseteq K(\varphi^t(u_0)) \supseteq K(u_0)$, we have
\[
 d_{\rm L}(\varphi^t,A) =
 \max_{u \in A}d_{\rm L}(\varphi^t,u)
 \leq
 \max_{u \in K(\varphi^t(u_0))}d_{\rm L}(\varphi^t,u)
 \leq
 \max_{u \in K(u_0)}d_{\rm L}(\varphi^t,u).
\]
Since $\rho\big(K(u_0),K(\varphi^t(u_0))\big) \to 0$ for $t \to +\infty$,
from the properties of decreasing \eqref{DOinctT}
and continuity (Lemma~\ref{thm:wcontu}),
it follows that
\[
  d_{\rm L} = \liminf_{t\to+\infty}\max_{u \in K(\varphi^t(u_0))}d_{\rm L}(\varphi^t,u).
\]

\section{Computation of the finite-time Lyapunov exponents and dimension in MATLAB}
The singular value decomposition (\emph{SVD}) of a fundamental matrix $D\varphi^t(u_0)$ has the from
\[
D\varphi^t(u_0)={\rm U}(t,u_0){\rm \Sigma}(t,u_0){\rm V}^{T}(t,u_0): \quad
{\rm U}(t,u_0)^T{\rm U}(t,u_0) \equiv I \equiv {\rm V}(t,u_0)^T{\rm V}(t,u_0),
\]
where ${\rm \Sigma}(t)=\text{\rm diag}\{\sigma_1(t,u_0),...,\sigma_n(t,u_0)\}$
is a  diagonal matrix with positive real diagonal entries --- \emph{singular values}.
We now give a MATLAB implementation \cite{LeonovKM-2015-EPJST} of the discrete SVD method
for computing finite-time Lyapunov exponents $\{\LEs_i(t,u_0)\}_1^n$
based on the product SVD algorithm (see, e.g., \cite{Stewart-1997,DieciL-2008}).
For the computation of the Lyapunov dimension of an attractor by the considered code
one has to consider a sufficiently large time $T$
and a grid of points on the attractor $K_{\rm grid}$,
compute the local Lyapunov dimensions by
the corresponding Kaplan-Yorke formula  $d_{\rm L}^{\rm KY}(\{\LEs_i(T,u)\}_{1}^{n})$
(see, e.g. \cite{KuznetsovMV-2014-CNSNS}), and takes maximum on the grid:
$
  \max_{u \in K_{\rm grid}} d_{\rm L}^{\rm KY}(\{\LEs_i(T,u)\}_{1}^{n}).
$

{
\beginMatlab{\textbf{productSVD.m} -- product SVD algorithm}
\fontsize{8}{8}\selectfont
\begin{lstlisting}
function [U, R, V] = productSVD(initFactorization, nIterations)
% Parameters:
%   initFactorization - the array contains factor matrices of the
%                       fundamental matrix X, such that:
%          X = initFactorization(:,:,1) * ... * initFactorization(:,:,end);
%   nIterations - the number of iterations in the product SVD algorithm.

% dimOde - dimension of the ODEs, nFactors - the number of factor matrices
[~, dimOde, nFactors] = size(initFactorization);

% A - 2d array of matrices storing the factor matrices at each iteration
A = zeros(dimOde, dimOde, nFactors, nIterations);
A(:, :, :, 1) = initFactorization;

% Q - array of matrices storing orhogonal matrices of the QR decomposition
Q = zeros(dimOde, dimOde, nFactors+1);

% U, V - orthogonal matrices in the SVD decomposition
U = eye(dimOde); V = eye(dimOde);

% R - array of upper triangular factor matrices, such that after
% the last iteration \Sigma = R(:,:,1) * ... * R(:,:,end)
R = zeros(dimOde, dimOde, nFactors);

% Main loop
for iIteration = 1 : nIterations
    Q(:, :, nFactors + 1) = eye(dimOde, dimOde);
    for jFactor = nFactors : -1 : 1
        C = A(:, :, jFactor, iIteration) * Q(:, :, jFactor+1);
        [Q(:, :, jFactor), R(:, :, jFactor)] = qr(C);
        for kCoord = 1 : dimOde
            if R(kCoord, kCoord, jFactor) < 0
                R(kCoord, :, jFactor) = -1 * R(kCoord, :, jFactor);
                Q(:, kCoord, jFactor) = -1 * Q(:, kCoord, jFactor);
            end;
        end;
    end;

    if mod(iIteration, 2) == 1
        U = U * Q(:, :, 1);
    else
        V = V * Q(:, :, 1);
    end

    for jFactor = 1 : nFactors
        A(:, :, jFactor, iIteration + 1) = R(:, :, nFactors-jFactor+1)';
    end
end

end
\end{lstlisting}
}

{
\beginMatlab{\textbf{computeLEs.m} -- computation of the Lyapunov exponents}
\fontsize{8}{8}\selectfont
\begin{lstlisting}
function LEs = computeLEs(extOde, initPoint, tStep, ...
                                nFactors, nSvdIterations, odeSolverOptions)
% Parameters:
%   extOde - extended ODE system (system of ODEs + var. eq.);
%   initPoint -  initial point;
%   tStep - time-step in the factorization procedure;
%   nFactors - number of factor matrices in the factorization procedure;
%   nSvdIterations - number of iterations in the product SVD algoritm;
%   odeSolverOptions - solver options (sover = ode45);

% Dimension of the ODE :
dimOde = length(initPoint);

% Dimension of the extended ODE (ODE + Var. Eq.):
dimExtOde = dimOde * (dimOde + 1);

tBegin = 0; tEnd = tStep;
tSpan = [tBegin, tEnd];
initFundMatrix = eye(dimOde);
initCond = [initPoint(:); initFundMatrix(:)];

X = zeros(dimOde, dimOde, nFactors);

% Main loop : factorization of the fundamental matrix
for iFactor = 1 : nFactors
    [~, extOdeSolution] = ode45(extOde, tSpan, initCond, odeSolverOptions);

    X(:, :, iFactor) = reshape(...
                        extOdeSolution(end, (dimOde + 1) : dimExtOde), ...
                                                           dimOde, dimOde);
    currInitPoint = extOdeSolution(end, 1 : dimOde);
    currInitFundMatrix = eye(dimOde);

    tBegin = tBegin + tStep;
    tEnd = tEnd + tStep;
    tSpan = [tBegin, tEnd];
    initCond = [currInitPoint(:); currInitFundMatrix(:)];
end

% Product SVD of factorization X of the fundamental matrix
[~, R, ~] = productSVD(X, nSvdIterations);

% Computation of the Lyapunov exponents
LEs = zeros(1, dimOde);
for jFactor = 1 : nFactors
    LEs = LEs + log(diag(R(:, :, jFactor))');
end;
finalTime = tStep * nFactors;
LEs = LEs / finalTime;

end
\end{lstlisting}
}

{
\beginMatlab{\textbf{lyapunovDim.m} -- computation of the Lyapunov dimension}
\fontsize{8}{8}\selectfont
\begin{lstlisting}
function LD = lyapunovDim( LEs )
% For the given array of finite-time Lyapunov exponents at a point the function
% computes the local Lyapunov dimension by the Kaplan-Yorke formula.

% Parameters:
%   LEs - array of the finite-time Lyapunov exponents.

% Initialization of the local Lyapunov dimension:
LD = 0;

% Number of LEs :
nLEs = length(LEs);

% Sorted LEs :
sortedLEs = sort(LEs, 'descend');

% Main loop :
leSum = sortedLEs(1);
if ( sortedLEs(1) > 0 )
     for i = 1 : nLEs-1
        if sortedLEs(i+1) ~= 0
           LD = i + leSum / abs( sortedLEs(i+1) );
           leSum = leSum + sortedLEs(i+1);
           if leSum < 0
              break;
           end
        end
    end
end
end
\end{lstlisting}
}

{
\beginMatlab{\textbf{genLorenzSyst.m} -- generalized Lorenz system
\eqref{sys:lorenz-general} along with the variational equation}
\fontsize{8}{8}\selectfont
\begin{lstlisting}
function OUT = genLorenzSyst(t, x, r, sigma, b, a)

% Generalized Lorenz system with
% parameters: r sigma b a

OUT(1) = sigma*(x(2) - x(1)) - a*x(2)*x(3);
OUT(2) = r*x(1) - x(2) - x(1)*x(3);
OUT(3) = -b*x(3) + x(1)*x(2);

% Jacobian at the point [x(1), x(2), x(3)]
J = [-sigma, sigma-a*x(3),  -a*x(2);
    r-x(3),    -1,    -x(1);
    x(2),     x(1), -b];

X = [x(4), x(7), x(10);
    x(5), x(8), x(11);
    x(6), x(9), x(12)];

% Variational equation
OUT(4:12) = J*X;
\end{lstlisting}
}

{
\beginMatlab{\textbf{main.m} -- computation of the Lyapunov exponents and
local Lyapunov dimension for the hidden attractor of generalized Lorenz system
\eqref{sys:lorenz-general}}
\fontsize{8}{8}\selectfont
\begin{lstlisting}
function main

% Parameters of generalized Lorenz system
% that correspond to the hidden attractor
r = 700; sigma = 4; b = 1; a = 0.0052;

% Initial point for the trajectory which visualizes the hidden attractor
x0 = [-14.551336132013954 -173.86811769236883 718.92035664071227];

tStep = 0.1;
nFactors = 10000;
nSvdIterations = 3;

% ODE solver parameters
acc = 1e-8; RelTol = acc; AbsTol = acc; InitialStep = acc/10;
odeSolverOptions = odeset('RelTol', RelTol, 'AbsTol', AbsTol, ...
                        'InitialStep', InitialStep, 'NormControl', 'on');

LEs = computeLEs(@(t, x) genLorenzSyst(t, x, r, sigma, b, a), ...
                    x0, tStep, nFactors, nSvdIterations, odeSolverOptions);

fprintf('Lyapunov exponents: %6.4f, %6.4f, %6.4f\n', LEs);

LD = lyapunovDim(LEs);

fprintf('Lyapunov dimension: %6.4f\n', LD);

end
\end{lstlisting}
}

\bigskip
\section*{Conclusions}
 In this survey for finite dimensional dynamical systems in Euclidean space
 we have tried to discuss rigorously
 the connection between the works
 by Kaplan and Yorke (the concept of Lyapunov dimension, 1979),
 Douady and Oesterl\'{e}
 (estimation of Hausdorff dimension via the Lyapunov dimension of maps, 1980),
 Constantin, Eden, Foias, and Temam
 (estimation of Hausdorff dimension via the Lyapunov exponents and dimension of dynamical systems, 1985-90),
 Leonov (estimation of the Lyapunov dimension via the direct Lyapunov method, 1991),
 and numerical methods for the computation of Lyapunov exponents and Lyapunov dimension.
 Remark that in the numerical estimations
 we can consider only finite time and get finite-time Lyapunov exponents,
 thus we have discussed the justification of Kaplan-Yorke formula
 with respect to finite-time Lyapunov exponent,
 by the Douady--Oesterl\'{e} theorem for maps.
 For various self-excited and hidden attractors of
 well-known dynamical systems, the numerical values,
 analytical estimations and formulas of the Lyapunov dimension are given.

\bigskip
\section*{Acknowledgments}
The authors would like to thank
Luis~Barreira, Igor~Chueshov, Alp~Eden, Volker~Reitmann for valuable comments on this work.
The  work  was supported by the Russian Scientific Foundation (project 14-21-00041)
and Saint-Petersburg State University.


\begin{thebibliography}{}

\bibitem[Abarbanel et~al., 1991]{AbarbanelBK-1991}
Abarbanel, H., Brown, R., and Kennel, M. (1991).
\newblock Variation of {L}yapunov exponents on a strange attractor.
\newblock {\em Journal of Nonlinear Science}, 1(2):175--199.

\bibitem[Adrianova, 1998]{Adrianova-1998}
Adrianova, L.~Y. (1998).
\newblock {\em Introduction to Linear systems of Differential Equations}.
\newblock American Mathematical Society, Providence, Rhode Island.

\bibitem[Aleksandrov and Pasynkov, 1973]{AleksandrovP-1973}
Aleksandrov, P. and Pasynkov, B. (1973).
\newblock {\em Introduction to Dimension Theory (in Russian)}.
\newblock Nauka, Moscow.

\bibitem[Babin and Vishik, 1992]{BabinV-1992}
Babin, A.~V. and Vishik, M.~I. (1992).
\newblock {\em Attractors of Evolution Equations}.
\newblock North-Holland, Amsterdam.

\bibitem[Barabanov, 2005]{Barabanov-2005}
Barabanov, E. (2005).
\newblock Singular exponents and properness criteria for linear differential
  systems.
\newblock {\em Differential Equations}, 41:151--162.

\bibitem[Barreira and Gelfert, 2011]{BarreiraG-2011}
Barreira, L. and Gelfert, K. (2011).
\newblock Dimension estimates in smooth dynamics: a survey of recent results.
\newblock {\em Ergodic Theory and Dynamical Systems}, 31:641--671.

\bibitem[Barreira and Schmeling, 2000]{BarreiraS-2000}
Barreira, L. and Schmeling, J. (2000).
\newblock Sets of ``{N}on-typical'' points have full topological entropy and
  full {H}ausdorff dimension.
\newblock {\em Israel Journal of Mathematics}, 116(1):29--70.

\bibitem[Bogoliubov and Krylov, 1937]{BogoliubovK-1937}
Bogoliubov, N. and Krylov, N. (1937).
\newblock La theorie generalie de la mesure dans son application a l'etude de
  systemes dynamiques de la mecanique non-lineaire.
\newblock {\em Ann. Math. II (in French) (Annals of Mathematics)},
  38(1):65--113.

\bibitem[Boichenko and Leonov, 1998]{BoichenkoL-1998}
Boichenko, V. and Leonov, G. (1998).
\newblock Lyapunov's direct method in estimates of topological entropy.
\newblock {\em Journal of Mathematical Sciences}, 91(6):3370--3379.

\bibitem[Boichenko et~al., 2005]{BoichenkoLR-2005}
Boichenko, V.~A., Leonov, G.~A., and Reitmann, V. (2005).
\newblock {\em Dimension Theory for Ordinary Differential Equations}.
\newblock Teubner, Stuttgart.

\bibitem[Brezetskyi et~al., 2015]{BrezetskyiDK-2015-HA}
Brezetskyi, S., Dudkowski, D., and Kapitaniak, T. (2015).
\newblock Rare and hidden attractors in van der {P}ol-{D}uffing oscillators.
\newblock {\em European Physical Journal: Special Topics}, 224(8):1459--1467.

\bibitem[Bylov et~al., 1966]{BylovVGN-1966}
Bylov, B.~E., Vinograd, R.~E., Grobman, D.~M., and Nemytskii, V.~V. (1966).
\newblock {\em Theory of characteristic exponents and its applications to
  problems of stability (in Russian)}.
\newblock Nauka, Moscow.

\bibitem[Chepyzhov and Vishik, 2002]{ChepyzhovV-2002}
Chepyzhov, V. and Vishik, M. (2002).
\newblock {\em Attractors for equations of mathematical physics}.
\newblock American Mathematical Society, Providence, Rhode Island.

\bibitem[Choquet and Foias, 1975]{ChoquetF-1975}
Choquet, G. and Foias, C. (1975).
\newblock Solution d'un probleme sur les iteres d'un operateur positif sur
  {$C(K)$} et proprietes de moyennes associees.
\newblock {\em Annales de l'institut Fourier (in French)}, 25(3-4):109--129.

\bibitem[Chueshov, 2002]{Chueshov-2002-book}
Chueshov, I. (2002).
\newblock {\em Introduction to the Theory of Infinite-dimensional Dissipative
  Systems}.
\newblock Electronic library of mathematics. ACTA.

\bibitem[Chueshov, 2015]{Chueshov-2015}
Chueshov, I. (2015).
\newblock {\em Dynamics of Quasi-Stable Dissipative Systems}.
\newblock Springer.

\bibitem[Chueshov and Siegmund, 2005]{ChueshovS-2005}
Chueshov, I. and Siegmund, S. (2005).
\newblock On dimension and metric properties of trajectory attractors.
\newblock {\em Journal of Dynamics and Differential Equations}, 17(4):621--641.

\bibitem[Chueshov, 1993]{Chueshov-1993}
Chueshov, I.~D. (1993).
\newblock Global attractors in the nonlinear problems of mathematical physics.
\newblock {\em Russian Mathematical Surveys}, 48(3):135--162.

\bibitem[Constantin and Foias, 1985]{ConstantinF-1985}
Constantin, P. and Foias, C. (1985).
\newblock Global {L}yapunov exponents, {K}aplan-{Y}orke formulas and the
  dimension of the attractors for 2{D} {N}avier-{S}tokes equations.
\newblock {\em Communications on Pure and Applied Mathematics}, 38(1):1--27.

\bibitem[Constantin et~al., 1985]{ConstantinFT-1985}
Constantin, P., Foias, C., and Temam, R. (1985).
\newblock Attractors representing turbulent flows.
\newblock {\em Memoirs of the American Mathematical Society}, 53(314).

\bibitem[Cvitanovi\'c et~al., 2012]{ChaosBook}
Cvitanovi\'c, P., Artuso, R., Mainieri, R., Tanner, G., and Vattay, G. (2012).
\newblock {\em Chaos: Classical and Quantum}.
\newblock Niels Bohr Institute, Copenhagen.
\newblock {h}ttp://ChaosBook.org.

\bibitem[Czornik et~al., 2013]{CzornikNN-2013}
Czornik, A., Nawrat, A., and Niezabitowski, M. (2013).
\newblock {L}yapunov exponents for discrete time-varying systems.
\newblock {\em Studies in Computational Intelligence}, 440:29--44.

\bibitem[Danca et~al., 2016]{DancaFKC-2016-IJBC}
Danca, M.-F., Feckan, M., Kuznetsov, N., and Chen, G. (2016).
\newblock Looking more closely at the {R}abinovich-{F}abrikant system.
\newblock {\em International Journal of Bifurcation and Chaos}, 26(02).
\newblock art. num. 1650038.

\bibitem[Dellnitz and Junge, 2002]{DellnitzJ-2002}
Dellnitz, M. and Junge, O. (2002).
\newblock Set oriented numerical methods for dynamical systems.
\newblock In {\em Handbook of Dynamical Systems}, volume~2, pages 221--264.
  Elsevier Science.

\bibitem[Dieci and Elia, 2008]{DieciL-2008}
Dieci, L. and Elia, C. (2008).
\newblock {S}{V}{D} algorithms to approximate spectra of dynamical systems.
\newblock {\em Mathematics and Computers in Simulation}, 79(4):1235--1254.

\bibitem[Doering et~al., 1987]{DoeringGHN-1987}
Doering, C., Gibbon, J., Holm, D., and Nicolaenko, B. (1987).
\newblock Exact {L}yapunov dimension of the universal attractor for the complex
  {G}inzburg-{L}andau equation.
\newblock {\em Phys. Rev. Lett.}, 59:2911--2914.

\bibitem[Douady and Oesterle, 1980]{DouadyO-1980}
Douady, A. and Oesterle, J. (1980).
\newblock Dimension de {H}ausdorff des attracteurs.
\newblock {\em C.R. Acad. Sci. Paris, Ser. A. (in French)}, 290(24):1135--1138.

\bibitem[Eden, 1989a]{Eden-1989-PhD}
Eden, A. (1989a).
\newblock {\em An abstract theory of {L}-exponents with applications to
  dimension analysis ({PhD} thesis)}.
\newblock Indiana University.

\bibitem[Eden, 1989b]{Eden-1989}
Eden, A. (1989b).
\newblock Local {L}yapunov exponents and a local estimate of {H}ausdorff
  dimension.
\newblock {\em ESAIM: Mathematical Modelling and Numerical Analysis -
  Modelisation Mathematique et Analyse Numerique}, 23(3):405--413.

\bibitem[Eden, 1990]{Eden-1990}
Eden, A. (1990).
\newblock Local estimates for the {H}ausdorff dimension of an attractor.
\newblock {\em Journal of Mathematical Analysis and Applications},
  150(1):100--119.

\bibitem[Eden et~al., 1991]{EdenFT-1991}
Eden, A., Foias, C., and Temam, R. (1991).
\newblock Local and global {L}yapunov exponents.
\newblock {\em Journal of Dynamics and Differential Equations}, 3(1):133--177.
\newblock [Preprint No. 8804, The Institute for Applied Mathematics and
  Scientific Computing, Indiana University, 1988].

\bibitem[Feng et~al., 2015]{FengPW-2015-HA}
Feng, Y., Pu, J., and Wei, Z. (2015).
\newblock Switched generalized function projective synchronization of two
  hyperchaotic systems with hidden attractors.
\newblock {\em European Physical Journal: Special Topics}, 224(8):1593--1604.

\bibitem[Feng and Wei, 2015]{FengW-2015-HA}
Feng, Y. and Wei, Z. (2015).
\newblock Delayed feedback control and bifurcation analysis of the generalized
  {S}prott {B} system with hidden attractors.
\newblock {\em European Physical Journal: Special Topics}, 224(8):1619--1636.

\bibitem[Frederickson et~al., 1983]{FredericksonKYY-1983}
Frederickson, P., Kaplan, J., Yorke, E., and Yorke, J. (1983).
\newblock The {L}iapunov dimension of strange attractors.
\newblock {\em Journal of Differential Equations}, 49(2):185--207.

\bibitem[Gelfert, 2003]{Gelfert-2003}
Gelfert, K. (2003).
\newblock Maximum local {L}yapunov dimension bounds the box dimension. {D}irect
  proof for invariant sets on {R}iemannian manifolds.
\newblock {\em Z. Anal. Anwend.}, 22:553--568.

\bibitem[Glukhovskii and Dolzhanskii, 1980]{GlukhovskyD-1980}
Glukhovskii, A.~B. and Dolzhanskii, F.~V. (1980).
\newblock Three-component geostrophic model of convection in a rotating fluid.
\newblock {\em Academy of Sciences, USSR, Izvestiya, Atmospheric and Oceanic
  Physics (in Russian)}, 16:311--318.

\bibitem[Gundlach and Steinkamp, 2000]{GundlachS-2000}
Gundlach, V. and Steinkamp, O. (2000).
\newblock Products of random rectangular matrices.
\newblock {\em Mathematische Nachrichten}, 212(1):51--76.

\bibitem[Hilbert, 1902]{Hilbert-1901}
Hilbert, D. (1901-1902).
\newblock Mathematical problems.
\newblock {\em Bull. Amer. Math. Soc.}, (8):437--479.

\bibitem[Horn and Johnson, 1994]{HornJ-1994-book}
Horn, R. and Johnson, C. (1994).
\newblock {\em Topics in Matrix Analysis}.
\newblock Cambridge University Press, Cambridge.

\bibitem[Hunt, 1996]{Hunt-1996}
Hunt, B. (1996).
\newblock Maximum local {L}yapunov dimension bounds the box dimension of
  chaotic attractors.
\newblock {\em Nonlinearity}, 9(4):845--852.

\bibitem[Hurewicz and Wallman, 1941]{HurewiczW-1941}
Hurewicz, W. and Wallman, H. (1941).
\newblock {\em Dimension Theory}.
\newblock Princeton University Press, Princeton.

\bibitem[Ilyashenko and Li, 1999]{IlyashenkoW-1999-AMS}
Ilyashenko, Y. and Li, W. (1999).
\newblock {\em Nonlocal Bifurcations}.
\newblock American Mathematical Society, Rhode Island.

\bibitem[Izobov, 2012]{Izobov-2012}
Izobov, N.~A. (2012).
\newblock {\em Lyapunov exponents and stability}.
\newblock Cambridge Scientific Publischers, Cambridge.

\bibitem[Jafari et~al., 2015]{JafariSN-2015-HA}
Jafari, S., Sprott, J., and Nazarimehr, F. (2015).
\newblock Recent new examples of hidden attractors.
\newblock {\em European Physical Journal: Special Topics}, 224(8):1469--1476.

\bibitem[Kaplan and Yorke, 1979]{KaplanY-1979}
Kaplan, J.~L. and Yorke, J.~A. (1979).
\newblock Chaotic behavior of multidimensional difference equations.
\newblock In {\em Functional Differential Equations and Approximations of Fixed
  Points}, pages 204--227. Springer, Berlin.

\bibitem[Kolmogorov, 1959]{Kolmogorov-1959}
Kolmogorov, A. (1959).
\newblock On entropy per unit time as a metric invariant of automorphisms.
\newblock {\em Dokl. Akad. Nauk SSSR (In Russian)}, 124(4):754--755.

\bibitem[Kuczma and Gil{\'a}nyi, 2009]{Kuczma-2009}
Kuczma, M. and Gil{\'a}nyi, A. (2009).
\newblock {\em An Introduction to the Theory of Functional Equations and
  Inequalities: Cauchy's Equation and Jensen's Inequality}.
\newblock Birkh{\"a}user Basel.

\bibitem[Kunze and Kupper, 2001]{KunzeK-2001}
Kunze, M. and Kupper, T. (2001).
\newblock Non-smooth dynamical systems: An overview.
\newblock In {\em Ergodic Theory, Analysis, and Efficient Simulation of
  Dynamical Systems}, pages 431--452. Springer.

\bibitem[Kuratowski, 1966]{Kuratowski-1966}
Kuratowski, K. (1966).
\newblock {\em Topology}.
\newblock Academic press, New York.

\bibitem[Kuznetsov, 2016a]{Kuznetsov-2016-ArXiv}
Kuznetsov, N. (2016a).
\newblock Estimation of {L}yapunov dimension via the {L}eonov method.
\newblock {\em arXiv}.
\newblock http://arxiv.org/pdf/1602.05410v1.pdf.

\bibitem[Kuznetsov, 2016b]{Kuznetsov-2016}
Kuznetsov, N. (2016b).
\newblock Hidden attractors in fundamental problems and engineering models. {A}
  short survey.
\newblock {\em Lecture Notes in Electrical Engineering}, 371:13--25.
\newblock (plenary lecture at AETA 2015: Recent Advances in Electrical
  Engineering and Related Sciences).

\bibitem[Kuznetsov et~al., 2016]{KuznetsovAL-2016}
Kuznetsov, N., Alexeeva, T., and Leonov, G. (2016).
\newblock Invariance of {L}yapunov exponents and {L}yapunov dimension for
  regular and irregular linearizations.
\newblock {\em Nonlinear Dynamics}.
\newblock (arXiv e-prints 1410.2016v2, 2014).

\bibitem[Kuznetsov and Leonov, 2014]{KuznetsovL-2014-IFACWC}
Kuznetsov, N. and Leonov, G. (2014).
\newblock Hidden attractors in dynamical systems: systems with no equilibria,
  multistability and coexisting attractors.
\newblock {\em IFAC Proceedings Volumes (IFAC-PapersOnline)}, 19:5445--5454.

\bibitem[Kuznetsov et~al., 2015]{KuznetsovLM-2015}
Kuznetsov, N., Leonov, G.~A., and Mokaev, T.~N. (2015).
\newblock Hidden attractor in the {R}abinovich system.
\newblock {\em arXiv:1504.04723v1}.
\newblock http://arxiv.org/pdf/1504.04723v1.pdf.

\bibitem[Kuznetsov et~al., 2014a]{KuznetsovAL-2014-arXiv-LE}
Kuznetsov, N.~V., Alexeeva, T., and Leonov, G.~A. (2014a).
\newblock Invariance of {L}yapunov characteristic exponents, {L}yapunov
  exponents, and {L}yapunov dimension for regular and non-regular
  linearizations.
\newblock {\em arXiv:1410.2016v2}.
\newblock (accepted to Nonlinear Dynamics).

\bibitem[Kuznetsov and Leonov, 2005]{KuznetsovL-2005}
Kuznetsov, N.~V. and Leonov, G.~A. (2005).
\newblock On stability by the first approximation for discrete systems.
\newblock In {\em 2005 International Conference on Physics and Control, PhysCon
  2005}, volume Proceedings Volume 2005, pages 596--599. IEEE.

\bibitem[Kuznetsov et~al., 2010]{KuznetsovLV-2010-IFAC}
Kuznetsov, N.~V., Leonov, G.~A., and Vagaitsev, V.~I. (2010).
\newblock Analytical-numerical method for attractor localization of generalized
  {C}hua's system.
\newblock {\em IFAC Proceedings Volumes (IFAC-PapersOnline)}, 4(1):29--33.

\bibitem[Kuznetsov et~al., 2014b]{KuznetsovMV-2014-CNSNS}
Kuznetsov, N.~V., Mokaev, T.~N., and Vasilyev, P.~A. (2014b).
\newblock Numerical justification of {L}eonov conjecture on {L}yapunov
  dimension of {R}ossler attractor.
\newblock {\em Commun Nonlinear Sci Numer Simulat}, 19:1027--1034.

\bibitem[Ladyzhenskaya, 1987]{Ladyzhenskaya-1987}
Ladyzhenskaya, O. (1987).
\newblock Determination of minimal global attractors for the {N}avier-{S}tokes
  equations and other partial differential equations.
\newblock {\em Russian Mathematical Surveys}, 42(6):25--60.

\bibitem[Ladyzhenskaya, 1991]{Ladyzhenskaya-1991}
Ladyzhenskaya, O.~A. (1991).
\newblock {\em Attractors for semi-groups and evolution equations}.
\newblock Cambridge University Press.

\bibitem[Ledrappier, 1981]{Ledrappier-1981}
Ledrappier, F. (1981).
\newblock Some relations between dimension and {L}yapounov exponents.
\newblock {\em Communications in Mathematical Physics}, 81(2):229--238.

\bibitem[Leonov, 2002]{Leonov-2002}
Leonov, G. (2002).
\newblock Lyapunov dimension formulas for {H}enon and {L}orenz attractors.
\newblock {\em St.Petersburg Mathematical Journal}, 13(3):453--464.

\bibitem[Leonov et~al., 2015a]{LeonovAK-2015}
Leonov, G., Alexeeva, T., and Kuznetsov, N. (2015a).
\newblock Analytic exact upper bound for the {L}yapunov dimension of the
  {S}himizu-{M}orioka system.
\newblock {\em Entropy}, 17(7):5101.

\bibitem[Leonov et~al., 2015b]{LeonovKKK-2015-arXiv-Lorenz}
Leonov, G., Kuznetsov, N., Korzhemanova, N., and Kusakin, D. (2015b).
\newblock The {L}yapunov dimension formula for the global attractor of the
  {L}orenz system.
\newblock {\em arXiv}, http://arxiv.org/pdf/1508.07498v1.pdf.

\bibitem[Leonov et~al., 2015c]{LeonovKKK-2015-arXiv-YangTigan}
Leonov, G., Kuznetsov, N., Korzhemanova, N., and Kusakin, D. (2015c).
\newblock {L}yapunov dimension formula of attractors in the {T}igan and {Y}ang
  systems.
\newblock {\em arXiv:1510.01492v1}, http://arxiv.org/pdf/1510.01492v1.pdf.

\bibitem[Leonov et~al., 2015d]{LeonovKM-2015-CNSNS}
Leonov, G., Kuznetsov, N., and Mokaev, T. (2015d).
\newblock Hidden attractor and homoclinic orbit in {L}orenz-like system
  describing convective fluid motion in rotating cavity.
\newblock {\em Communications in Nonlinear Science and Numerical Simulation},
  28:166--174.

\bibitem[Leonov et~al., 2015e]{LeonovKM-2015-EPJST}
Leonov, G., Kuznetsov, N., and Mokaev, T. (2015e).
\newblock Homoclinic orbits, and self-excited and hidden attractors in a
  {L}orenz-like system describing convective fluid motion.
\newblock {\em Eur. Phys. J. Special Topics}, 224(8):1421--1458.

\bibitem[Leonov and Lyashko, 1993]{LeonovL-1993}
Leonov, G. and Lyashko, S. (1993).
\newblock {E}den's hypothesis for a {L}orentz system.
\newblock {\em Vestnik St. Petersburg University: Mathematics}, 26(3):15--18.
\newblock [Transl. from Russian. Vestnik Sankt-Peterburgskogo Universiteta. Ser
  1. Matematika, 26(3), 14-16].

\bibitem[Leonov et~al., 2012a]{LeonovPS-2013-PLA}
Leonov, G., Pogromsky, A., and Starkov, K. (2012a).
\newblock Erratum to ''{T}he dimension formula for the {L}orenz attractor''
  [{P}hys. {L}ett. {A} 375 (8) (2011) 1179].
\newblock {\em Physics Letters A}, 376(45):3472 -- 3474.

\bibitem[Leonov and Poltinnikova, 2005]{LeonovP-2005}
Leonov, G. and Poltinnikova, M. (2005).
\newblock On the {L}yapunov dimension of the attractor of {C}hirikov
  dissipative mapping.
\newblock {\em AMS Translations. Proceedings of St.Petersburg Mathematical
  Society. Vol. X}, 224:15--28.

\bibitem[Leonov, 1991]{Leonov-1991-Vest}
Leonov, G.~A. (1991).
\newblock On estimations of {H}ausdorff dimension of attractors.
\newblock {\em Vestnik St. Petersburg University: Mathematics}, 24(3):38--41.
\newblock [Transl from Russian. Vestnik Leningradskogo Universiteta.
  Mathematika, 24(3), 1991, pp.~41-44].

\bibitem[Leonov, 2008]{Leonov-2008}
Leonov, G.~A. (2008).
\newblock {\em Strange attractors and classical stability theory}.
\newblock St.Petersburg University Press, St.Petersburg.

\bibitem[Leonov, 2012]{Leonov-2012-PMM}
Leonov, G.~A. (2012).
\newblock Lyapunov functions in the attractors dimension theory.
\newblock {\em Journal of Applied Mathematics and Mechanics}, 76(2):129--141.

\bibitem[Leonov and Boichenko, 1992]{LeonovB-1992}
Leonov, G.~A. and Boichenko, V.~A. (1992).
\newblock Lyapunov's direct method in the estimation of the {H}ausdorff
  dimension of attractors.
\newblock {\em Acta Applicandae Mathematicae}, 26(1):1--60.

\bibitem[Leonov et~al., 2015f]{LeonovKM-2015-ArXiv}
Leonov, G.~A., Kuznetsov, N., and Mokaev, T.~N. (2015f).
\newblock The {L}yapunov dimension formula of self-excited and hidden
  attractors in the {G}lukhovsky-{D}olzhansky system.
\newblock {\em arXiv:1509.09161}.
\newblock http://arxiv.org/pdf/1509.09161v1.pdf.

\bibitem[Leonov and Kuznetsov, 2007]{LeonovK-2007}
Leonov, G.~A. and Kuznetsov, N.~V. (2007).
\newblock Time-varying linearization and the {P}erron effects.
\newblock {\em International Journal of Bifurcation and Chaos},
  17(4):1079--1107.

\bibitem[Leonov and Kuznetsov, 2013]{LeonovK-2013-IJBC}
Leonov, G.~A. and Kuznetsov, N.~V. (2013).
\newblock Hidden attractors in dynamical systems. {F}rom hidden oscillations in
  {H}ilbert-{K}olmogorov, {A}izerman, and {K}alman problems to hidden chaotic
  attractors in {C}hua circuits.
\newblock {\em International Journal of Bifurcation and Chaos}, 23(1).
\newblock {a}rt. no. 1330002.

\bibitem[Leonov and Kuznetsov, 2015]{LeonovK-2015-AMC}
Leonov, G.~A. and Kuznetsov, N.~V. (2015).
\newblock On differences and similarities in the analysis of {L}orenz, {C}hen,
  and {L}u systems.
\newblock {\em Applied Mathematics and Computation}, 256:334--343.

\bibitem[Leonov et~al., 2011]{LeonovKV-2011-PLA}
Leonov, G.~A., Kuznetsov, N.~V., and Vagaitsev, V.~I. (2011).
\newblock Localization of hidden {C}hua's attractors.
\newblock {\em Physics Letters A}, 375(23):2230--2233.

\bibitem[Leonov et~al., 2012b]{LeonovKV-2012-PhysD}
Leonov, G.~A., Kuznetsov, N.~V., and Vagaitsev, V.~I. (2012b).
\newblock Hidden attractor in smooth {C}hua systems.
\newblock {\em Physica D: Nonlinear Phenomena}, 241(18):1482--1486.

\bibitem[Li et~al., 2015]{Li20151493}
Li, C., Hu, W., Sprott, J., and Wang, X. (2015).
\newblock Multistability in symmetric chaotic systems.
\newblock {\em European Physical Journal: Special Topics}, 224(8):1493--1506.

\bibitem[Lorenz, 1963]{Lorenz-1963}
Lorenz, E.~N. (1963).
\newblock Deterministic nonperiodic flow.
\newblock {\em J. Atmos. Sci.}, 20(2):130--141.

\bibitem[Lyapunov, 1892]{Lyapunov-1892}
Lyapunov, A.~M. (1892).
\newblock {\em The General Problem of the Stability of Motion (in Russian)}.
\newblock Kharkov.
\newblock [English transl. Academic Press, NY, 1966].

\bibitem[Millionschikov, 1976]{Millionschikov-1976}
Millionschikov, V.~M. (1976).
\newblock A formula for the entropy of smooth dynamical systems.
\newblock {\em Differencial'nye Uravenija (in Russian)}, 12(12):2188--2192,
  2300.

\bibitem[Milnor, 2006]{Milnor-2006}
Milnor, J. (2006).
\newblock Attractor.
\newblock {\em Scholarpedia}, 1(11).
\newblock {d}oi:10.4249/scholarpedia.1815.

\bibitem[Noack and Reitmann, 1996]{NoackR-1996}
Noack, A. and Reitmann, V. (1996).
\newblock Hausdorff dimension estimates for invariant sets of time-dependent
  vector fields.
\newblock {\em Z. Anal. Anwend.}, 15:457--473.

\bibitem[Oseledec, 1968]{Oseledec-1968}
Oseledec, V. (1968).
\newblock Multiplicative ergodic theorem: Characteristic {L}yapunov exponents
  of dynamical systems.
\newblock In {\em Transactions of the Moscow Mathematical Society}, volume~19,
  pages 179--210.

\bibitem[Ott et~al., 1984]{OttWY-1984}
Ott, E., Withers, W., and Yorke, J. (1984).
\newblock Is the dimension of chaotic attractors invariant under coordinate
  changes?
\newblock {\em Journal of Statistical Physics}, 36(5-6):687--697.

\bibitem[Ott and Yorke, 2008]{OttY-2008}
Ott, W. and Yorke, J. (2008).
\newblock When {L}yapunov exponents fail to exist.
\newblock {\em Phys. Rev. E}, 78:056203.

\bibitem[Pesin, 1977]{Pesin-1977}
Pesin, Y. (1977).
\newblock Characteristic {L}yapunov exponents and smooth ergodic theory.
\newblock {\em Russian Mathematical Surveys}, 32(4):55--114.

\bibitem[Pham et~al., 2015]{Pham20151507}
Pham, V., Vaidyanathan, S., Volos, C., and Jafari, S. (2015).
\newblock Hidden attractors in a chaotic system with an exponential nonlinear
  term.
\newblock {\em European Physical Journal: Special Topics}, 224(8):1507--1517.

\bibitem[Pilyugin, 2011]{Pilyugin-2011}
Pilyugin, S. (2011).
\newblock Theory of pseudo-orbit shadowing in dynamical systems.
\newblock {\em Differential Equations}, 47(13):1929--1938.

\bibitem[Pogromsky and Matveev, 2011]{PogromskyM-2011}
Pogromsky, A.~Y. and Matveev, A.~S. (2011).
\newblock Estimation of topological entropy via the direct {L}yapunov method.
\newblock {\em Nonlinearity}, 24(7):1937.

\bibitem[Saha et~al., 2015]{SahaSRC-2015-HA}
Saha, P., Saha, D., Ray, A., and Chowdhury, A. (2015).
\newblock Memristive non-linear system and hidden attractor.
\newblock {\em European Physical Journal: Special Topics}, 224(8):1563--1574.

\bibitem[Schmeling, 1998]{Schmeling-1998}
Schmeling, J. (1998).
\newblock A dimension formula for endomorphisms -- the {B}elykh family.
\newblock {\em Ergodic Theory and Dynamical Systems}, 18:1283--1309.

\bibitem[Sell, 1996]{Sell-1996}
Sell, G.~R. (1996).
\newblock Global attractors for the three-dimensional {N}avier-{S}tokes
  equations.
\newblock {\em Journal of Dynamics and Differential Equations}, 8(1):1--33.

\bibitem[Semenov et~al., 2015]{Semenov20151553}
Semenov, V., Korneev, I., Arinushkin, P., Strelkova, G., Vadivasova, T., and
  Anishchenko, V. (2015).
\newblock Numerical and experimental studies of attractors in memristor-based
  {C}hua's oscillator with a line of equilibria. {N}oise-induced effects.
\newblock {\em European Physical Journal: Special Topics}, 224(8):1553--1561.

\bibitem[Shahzad et~al., 2015]{ShahzadPAJH-2015-HA}
Shahzad, M., Pham, V.-T., Ahmad, M., Jafari, S., and Hadaeghi, F. (2015).
\newblock Synchronization and circuit design of a chaotic system with
  coexisting hidden attractors.
\newblock {\em European Physical Journal: Special Topics}, 224(8):1637--1652.

\bibitem[Sharma et~al., 2015]{SharmaSPKL-2015-EPJST}
Sharma, P., Shrimali, M., Prasad, A., Kuznetsov, N., and Leonov, G. (2015).
\newblock Control of multistability in hidden attractors.
\newblock {\em Eur. Phys. J. Special Topics}, 224(8):1485--1491.

\bibitem[Shimizu and Morioka, 1980]{Shimizu1980201}
Shimizu, T. and Morioka, N. (1980).
\newblock On the bifurcation of a symmetric limit cycle to an asymmetric one in
  a simple model.
\newblock {\em Physics Letters A}, 76(3-4):201 -- 204.

\bibitem[Sinai, 1959]{Sinai-1959}
Sinai, Y. (1959).
\newblock On the notion of entropy of dynamical systems.
\newblock {\em Dokl. Akad. Nauk SSSR (In Russian)}, 124(4):768--771.

\bibitem[Smith, 1986]{Smith-1986}
Smith, R. (1986).
\newblock Some application of {H}ausdorff dimension inequalities for ordinary
  differential equation.
\newblock {\em Proc. Royal Society Edinburg}, 104A:235--259.

\bibitem[Sparrow, 1982]{Sparrow-1982}
Sparrow, C. (1982).
\newblock {\em The {L}orenz Equations: Bifurcations, Chaos, and Strange
  Attractors}.
\newblock Applied Mathematical Sciences. Springer New York.

\bibitem[Sprott, 2015]{Sprott20151409}
Sprott, J. (2015).
\newblock Strange attractors with various equilibrium types.
\newblock {\em European Physical Journal: Special Topics}, 224(8):1409--1419.

\bibitem[Stewart, 1997]{Stewart-1997}
Stewart, D.~E. (1997).
\newblock A new algorithm for the {S}{V}{D} of a long product of matrices and
  the stability of products.
\newblock {\em Electronic Transactions on Numerical Analysis}, 5:29--47.

\bibitem[Temam, 1997]{Temam-1997}
Temam, R. (1997).
\newblock {\em Infinite-dimensional Dynamical Systems in Mechanics and
  Physics}.
\newblock Springer-Verlag, New York, 2nd edition.

\bibitem[Tigan and Opris, 2008]{Tigan-2008}
Tigan, G. and Opris, D. (2008).
\newblock Analysis of a 3d chaotic system.
\newblock {\em Chaos, Solitons \& Fractals}, 36(5):1315--1319.

\bibitem[Vaidyanathan et~al., 2015]{VaidyanathanPV-2015-HA}
Vaidyanathan, S., Pham, V.-T., and Volos, C. (2015).
\newblock A 5-{D} hyperchaotic {R}ikitake dynamo system with hidden attractors.
\newblock {\em European Physical Journal: Special Topics}, 224(8):1575--1592.

\bibitem[Yang and Chen, 2008]{Yangc-2008}
Yang, Q. and Chen, G. (2008).
\newblock A chaotic system with one saddle and two stable node-foci.
\newblock {\em International Journal of Bifurcation and Chaos}, 18:1393--1414.

\bibitem[Young, 2013]{Young-2013}
Young, L.-S. (2013).
\newblock Mathematical theory of {L}yapunov exponents.
\newblock {\em Journal of Physics A: Mathematical and Theoretical},
  46(25):254001.

\bibitem[Zakrzhevsky et~al., 2007]{ZakrzhevskySY-2007}
Zakrzhevsky, M., Schukin, I., and Yevstignejev, V. (2007).
\newblock {\em Scientific Proc. Riga Technical Univ. Transp. Engin.}, 6:79.

\bibitem[Zelinka, 2015]{Zelinka-2015}
Zelinka, I. (2015).
\newblock A survey on evolutionary algorithms dynamics and its complexity --
  {M}utual relations, past, present and future.
\newblock {\em Swarm and Evolutionary Computation}, 25:2--14.

\bibitem[Zelinka, 2016]{Zelinka-2016-HA}
Zelinka, I. (2016).
\newblock Evolutionary identification of hidden chaotic attractors.
\newblock {\em Engineering Applications of Artificial Intelligence},
  50:159--167.
\newblock 10.1016/j.engappai.2015.12.002.

\bibitem[Zhusubaliyev et~al., 2015]{ZhusubaliyevMCM-2015-HA}
Zhusubaliyev, Z., Mosekilde, E., Churilov, A., and Medvedev, A. (2015).
\newblock Multistability and hidden attractors in an impulsive {G}oodwin
  oscillator with time delay.
\newblock {\em European Physical Journal: Special Topics}, 224(8):1519--1539.

\end{thebibliography}

\end{document}